\documentclass[twoside,11pt]{article}

\usepackage{lastpage}

\usepackage{natbib}%
\RequirePackage{epsfig}
\RequirePackage{amssymb}
\RequirePackage{natbib}
\RequirePackage{graphicx}
\if@usehyper
\RequirePackage[colorlinks=false,allbordercolors={1 1 1}]{hyperref}
\fi
\if@abbrvbib
\else
\fi
\bibpunct{(}{)}{;}{a}{,}{,}

\usepackage{import}
\usepackage[utf8]{inputenc}
\usepackage[USenglish]{babel}
\usepackage{geometry}
\usepackage{xpatch}
\usepackage{amssymb,amsmath,amsfonts,bm}
\usepackage{subcaption}
\usepackage{xcolor}
\usepackage{enumitem}
\usepackage{array}
\usepackage{xcolor}
\usepackage{xparse}
\usepackage{etoolbox}
\usepackage{calc}
\usepackage{dsfont}
\usepackage{mathtools}
\usepackage{pict2e}%

\usepackage{tikz}
\usetikzlibrary{arrows.meta}
\usetikzlibrary{calc}
\usetikzlibrary{arrows}
\usetikzlibrary{decorations.text,decorations.pathreplacing,matrix,calc,positioning}
\usetikzlibrary{calligraphy}%
\usetikzlibrary{patterns}
\usetikzlibrary{patterns.meta}
\pgfdeclarepatternformonly{hatch1}{\pgfqpoint{-2pt}{-2pt}}{\pgfqpoint{8pt}{8pt}}{\pgfqpoint{5pt}{5pt}}%
{
\pgfsetlinewidth{0mm} %
\pgfpathmoveto{\pgfqpoint{-1pt}{-1pt}}
\pgfpathlineto{\pgfqpoint{5pt}{5pt}}
\pgfsetstrokecolor{gray}
\pgfusepath{stroke}
\pgfpathmoveto{\pgfqpoint{5pt}{-1pt}}
\pgfpathlineto{\pgfqpoint{-1pt}{5pt}}
\pgfsetstrokecolor{gray}
\pgfusepath{stroke}
}
\tikzset{ 
    table/.style={
        matrix of math nodes,
        row sep=-\pgflinewidth,
        column sep=-\pgflinewidth,
        nodes={rectangle},
        text depth=1ex,
        text height=1ex,
        nodes in empty cells,
        left delimiter=[,
        right delimiter={]},
        ampersand replacement=\&
    }
}

\usepackage{booktabs}
\usepackage{rotating}%
\usepackage{adjustbox}%
\usepackage{wrapfig}%

\usepackage{graphicx}
\newcommand{\trimlr}{1mm} %
\newcommand{\trimbottom}{5mm} %
\newcommand{\trimtop}{15mm} %

\newenvironment{subsubfigure}[2][]{%
  \begin{subfigure}[#1]{#2}%
    \stepcounter{subsubfigure}%
}{%
    \addtocounter{subfigure}{-1}%
  \end{subfigure}%
}
\newcounter{subsubfigure}

\usepackage{hyperref}%
\usepackage[capitalise]{cleveref}
\crefname{proposition}{Proposition}{Propositions}
\Crefname{proposition}{Proposition}{Propositions}

\crefname{lemma}{Lemma}{Lemmas}
\Crefname{lemma}{Lemma}{Lemmas}
\newtheorem{lemma}{Lemma}

\usepackage{pgfplots}
\pgfplotsset{compat=newest}

\usepackage{algpseudocode}
\usepackage{algorithm}
\makeatletter
\renewcommand{\ALG@beginalgorithmic}{\small}%
\makeatother
\algnewcommand{\Initialize}[1]{%
  \Statex \textbf{Initialize:}
  #1
}

\algblockdefx{Nest}{EndNest}{}{}%
\algtext*{Nest}{}%
\algtext*{EndNest}{}%

\newcounter{algsubstate}

\newcounter{algsubsubstate}
\renewcommand{\thealgsubsubstate}{\roman{algsubsubstate}}

\algnewcommand{\lFor}[1]{\State\algorithmicfor\ #1\ \algorithmicdo} %
\algnewcommand{\EndlFor}{\unskip\ \algorithmicend\ \algorithmicfor} %
\algnewcommand{\IIf}[1]{\State\algorithmicif\ #1\ \algorithmicthen}%
\algnewcommand{\EndIIf}{\State\algorithmicend\ \algorithmicif}%
\algnewcommand{\EndIIfi}{\algorithmicend\ \algorithmicif}%
\algnewcommand{\ElseIIf}[1]{\State\algorithmicelse\ \algorithmicif\ #1\ \algorithmicthen} %

\makeatletter
\renewenvironment{subequations}{%
  \refstepcounter{equation}%
  \protected@edef\theparentequation{\theequation}%
  \setcounter{parentequation}{\value{equation}}%
  \setcounter{equation}{0}%
  \def\theequation{\theparentequation.\arabic{equation}}%
  \ignorespaces
}{%
  \setcounter{equation}{\value{parentequation}}%
  \ignorespacesafterend
}
\makeatother

\let\texdisplaystyle\displaystyle
\def\displaytotextstyle{\textstyle\let\displaystyle\texdisplaystyle}

\newenvironment{talign*}
 {\let\displaystyle\displaytotextstyle\csname align*\endcsname}
 {\endalign}

\newcommand{\mathopfont}{\mathsf}

\DeclareMathOperator*{\argmin}{arg\,min}
\DeclareMathOperator*{\Argmin}{arg\,min}

\DeclareMathOperator*{\projection}{\mathopfont{proj}}

\newcommand{\R}{\mathbb{R}}
\newcommand{\ind}{\mathds{1}}
\DeclarePairedDelimiter\abs{\lvert}{\rvert}%
\DeclarePairedDelimiter\ceil{\lceil}{\rceil}

\DeclarePairedDelimiterXPP\sign[1]{\operatorname{\mathopfont{sign}}}[]{}{#1}%

\DeclarePairedDelimiter\norm{\lVert}{\rVert}%
\DeclarePairedDelimiterX{\inner}[2]{\langle}{\rangle}{#1, #2}

\DeclarePairedDelimiterXPP\mydet[1]{\operatorname{\mathopfont{det}}}[]{}{#1}%

\let\oldvec\vec
\newcommand{\sorth}[1]{\reflectbox{\ensuremath{\oldvec{\reflectbox{\ensuremath{#1}}}}}}

\DeclarePairedDelimiterXPP\myvec[1]{\operatorname{\mathopfont{vec}}}[]{}{#1}%
\renewcommand{\vec}[1][]{\myvec{#1}}
\DeclarePairedDelimiterXPP\tangentcone[2]{\mathcal{T}_{#1}}(){}{#2}
\DeclarePairedDelimiterXPP\lineal[1]{\operatorname{\mathopfont{lin}}}[]{}{#1}%
\DeclarePairedDelimiterXPP\rank[1]{\operatorname{\mathopfont{rank}}}[]{}{#1}%
\DeclarePairedDelimiterXPP\nullspace[1]{\operatorname{\mathopfont{null}}}[]{}{#1}%
\DeclarePairedDelimiterX{\mathargument}[1]{(}{)}{\mathargumentA{#1}}
\NewDocumentCommand{\mathargumentA}{>{\SplitArgument{1}{|}}m}{%
  \mathargumentB#1%
}
\NewDocumentCommand{\mathargumentB}{mm}{%
  \IfNoValueTF{#2}{%
    #1%
  }{%
    #1\;\delimsize|\;#2%
  }%
}

\NewDocumentCommand{\indicator}{e{_}}{%
  \ind\IfValueT{#1}{_{\{#1\}}}\mathargument
}

\DeclarePairedDelimiter\rdbracket{(}{)}
\NewDocumentCommand\maxsum{t{'}e{_^}}{%
\mathopfont{T}%
\IfBooleanT{#1}{'}%
\IfNoValueF{#2}{_{#2}}%
\IfNoValueF{#3}{^{#3}}%
\rdbracket}
\NewDocumentCommand\proj{t{'}e{_^}}{%
\projection%
\textstyle%
\IfBooleanT{#1}{'}%
\IfNoValueF{#2}{_{\!#2\!}}%
\IfNoValueF{#3}{^{\!#3\!}}%
\rdbracket}
\NewDocumentCommand\prox{t{'}e{_^}}{%
\mathopfont{prox}%
\textstyle%
\IfBooleanT{#1}{'}%
\IfNoValueF{#2}{_{#2}}%
\IfNoValueF{#3}{^{#3}}%
\rdbracket}

\DeclarePairedDelimiterX{\probargument}[1]{(}{)}{\probargumentA{#1}}
\NewDocumentCommand{\probargumentA}{>{\SplitArgument{1}{|}}m}{%
  \probargumentB#1%
}
\NewDocumentCommand{\probargumentB}{mm}{%
  \IfNoValueTF{#2}{%
    #1%
  }{%
    #1\;\delimsize|\;#2%
  }%
}
\NewDocumentCommand{\prob}{e{_}}{%
  \mathopfont{P}\IfValueT{#1}{_{#1}}\probargument
}
\NewDocumentCommand{\expect}{e{_}}{%
  \mathopfont{E}\IfValueT{#1}{_{#1}}\probargument
}
\NewDocumentCommand{\cvar}{e{_}}{%
  \mathopfont{CVaR}\IfValueT{#1}{_{#1}}\probargument
}
\NewDocumentCommand{\quantile}{e{_}}{%
  \mathopfont{Q}\IfValueT{#1}{_{#1}}\probargument
}
\NewDocumentCommand{\superquantile}{e{_}}{%
  \mathopfont{CVaR}\IfValueT{#1}{_{#1}}\probargument
}

\newcommand{\ie}{\emph{i.e.,\space}}
\newcommand{\eg}{\emph{e.g.,\space}}

\allowdisplaybreaks

\hypersetup{
	linktoc=all,
	colorlinks=true,
	linkcolor={black},%
	filecolor=black,
	urlcolor=black,
	citecolor={black},
}

\usepackage{optidef}

\begin{document}

\title{Fast Computation of Superquantile-Constrained Optimization Through Implicit Scenario Reduction}

\author{Jake Roth \thanks{Department of Industrial and Systems Engineering, 
University of Minnesota, Minneapolis, MN 55414 (roth0674@umn.edu).}
\and Ying Cui 
\thanks{Department of Industrial Engineering and Operations Research, University of California, Berkeley, Berkeley, CA 94720 (yingcui@berkeley.edu).}
}
\date{}

\maketitle

\begin{abstract}%
Superquantiles have recently gained significant interest as a risk-aware metric for addressing fairness and distribution shifts in statistical learning and decision making problems. This paper introduces a fast, scalable and robust second-order computational framework to solve large-scale optimization problems with superquantile-based constraints. Unlike empirical risk minimization, superquantile-based optimization requires ranking random functions evaluated across all scenarios to compute the tail conditional expectation. While this tail-based feature might seem computationally unfriendly, it provides an advantageous setting for a semismooth-Newton-based augmented Lagrangian method. The superquantile operator effectively reduces the dimensions of the Newton systems since the tail expectation involves considerably fewer scenarios. Notably, the extra cost of obtaining relevant second-order information and performing matrix inversions is often comparable to, and sometimes even less than, the effort required for gradient computation. Our developed solver is particularly effective when the number of scenarios substantially exceeds the number of decision variables. In synthetic problems with linear and convex diagonal quadratic objectives, numerical experiments demonstrate that our method outperforms existing approaches by a large margin: It achieves speeds more than 750 times faster for linear and quadratic objectives than the alternating direction method of multipliers as implemented by OSQP for computing low-accuracy solutions. Additionally, it is up to 25 times faster for linear objectives and 70 times faster for quadratic objectives than the commercial solver Gurobi, and 20 times faster for linear objectives and 30 times faster for quadratic objectives than the Portfolio Safeguard optimization suite for high-accuracy solution computations. For the quantile regression problem involving over 30 million scenarios, our method computes solution paths up to 20 times faster than Gurobi. The Julia implementation of the solver is available at \url{https://github.com/jacob-roth/superquantile-opt}.
\end{abstract}

\section{Introduction}
\label{sec:intro}

In this paper, we develop a fast, scalable and robust second-order computational solver for the following superquantile-constrained convex optimization problem:
\begin{equation}\label{eq:source_problem}
\begin{array}{rl}
\displaystyle\operatornamewithlimits{minimize}_{x\in X\subseteq \R^n} \quad & \displaystyle f(x) + \superquantile_{\tau_0}[\big]{G^{\,0}(x;\omega^{[m]})} \\[0.1in]
\mbox{subject to}
& \superquantile_{\tau_{\ell}}[\big]{G^{\,\ell}(x;\omega^{[m]})} \leq 0,\; \ell \in\{1,2,\ldots,L\},
\end{array}
\end{equation}
where $f:\mathbb{R}^n \to \mathbb{R}$ is a convex and continuously differentiable function, $X\subseteq \mathbb{R}^n$ is a closed convex set with a strongly semismooth and easily computable Euclidean projection (such as a box), and for each $\ell \in \{0, 1, \ldots, L\}$,  $G^\ell(x;\omega^{[m]}) \coloneqq \left\{g^\ell(x;\omega^j)\right\}_{j=1}^m$ is a collection  of random convex continuously differentiable functions $g^\ell(\,\cdot\,; \omega)$ evaluated at scenarios $\{w^j\}_{j=1}^{m}$~\footnote{For simplicity of notation, we assume that  scenarios $\{\omega_j\}_{j=1}^m$ are the same for all  $\ell\in\{0,1,\ldots,L\}$.}, and 
\begin{equation}\label{cvar}
\superquantile_{\tau_\ell}{G^\ell(x; \omega^{[m]})}\coloneqq \displaystyle\inf_t\left\{t + \displaystyle\frac{1}{(1-\tau_\ell)m}\sum_{j=1}^m \max\left(g^\ell(x; \omega^j) - t, 0\right) \right\}
\end{equation}
denotes the empirical superquantile (also known as the conditional value-at-risk) of $G^\ell(x; \omega^{[m]})$ at the confidence level $\tau_\ell\in (0,1)$ \citep{rockafellar1999cvar}.
In the remainder of this paper, we assume that $\tau_\ell$ is chosen so that $k(\tau_\ell)\coloneqq(1-\tau_\ell) m$ is an integer.
Since the superquantile operator preserves the convexity, one can easily see that \eqref{eq:source_problem} is a convex optimization problem.

The optimization model \eqref{eq:source_problem} is flexible and  broadly applicable.
The minimization of a superquantile objective can be interpreted as distributionally robust optimization with a $\varphi$-divergence ambiguity set \citep{bertsimas2009constructing}, and varying the confidence level $\tau_0$ from $0$ to $1$ in \eqref{eq:source_problem} provides a mechanism to interpolate between (risk-neutral) stochastic programs and robust optimization formulations.
There is extensive literature on using the superquantile objectives and constraints to model risk-averse preferences  and ensure decisions against adverse uncertainty; for an early example, see \citet{pflug200some}.
In engineering design, superquantiles help identify potential vulnerabilities \citep{airaudo2023risk}, while in statistics, they can be used to formulate the quantile regression problem \citep{royset2021primer}.
Recently, risk-aware models have attracted attention in the machine learning community as alternatives to empirical risk minimization.
These models are especially valuable in situations where decision-makers aim to mitigate distribution shifts between training and test datasets \citep{laguel2021superquantiles,robey2022probabilistic}, or prioritize performance on certain regions of the distribution, \eg losses over the left tail to prevent over-fitting or the right tail to exclude outliers  \citep{yuan2020group,peng2022imbalanced}.
Meanwhile, superquantile constraints have been used in engineering to certify solutions against adverse events as well as adverse computational outcomes \citep{roysetOA2011,chaudhuri2022certifiable},
in operations research as the convex, conservative approximation of chance constraints \citep{nemirovski2006convex,chen2010cvar}, and in machine learning to promote fairness \citep{williamson2019fairness,liu2019human,frohlich2022risk}. 
Interested readers are refereed to the recent survey of \citet{royset2022risk} for a comprehensive overview of risk-aware optimization.
Despite the promising numerical results demonstrated by superquantile models in various applications, the absence of a highly stable, efficient, and scalable numerical solver has somewhat limited their widespread adoption.

In fact, by introducing auxiliary variables $\{u_{\ell j}\}_{j=1, \ldots, m}^{\ell=0, \ldots, L}\subseteq \mathbb{R}$, one can rewrite \eqref{eq:source_problem} in the form of the  standard nonlinear optimization problem:
\begin{equation}
    \label{eq:nlp}
    \begin{array}{rl}
\displaystyle\operatornamewithlimits{minimize}_{\substack{x\in X\subseteq\mathbb{R}^n \\ \{t_\ell\}_{\ell = 0, \ldots, L}, \{u_{\ell j}\}_{j=1, \ldots, m}^{\ell=0, \ldots, L}}} & f(x) + t_0 + \displaystyle\frac{1}{(1-\tau_0)m} \displaystyle\sum_{j=1}^m u_{0 j} \\[0.1in]
\mbox{subject to} & t_\ell + \displaystyle\frac{1}{(1-\tau_\ell)m} \sum_{j=1}^m u_{\ell j} \leq 0, \; \ell \in \{0, \ldots, L\} \\[0.2in]
& u_{\ell j} \geq g^\ell(x;\omega^j) - t_\ell, \; u_{\ell j}\geq 0, \; j\in \{1, \ldots, m\}, \; \ell \in \{0, \ldots, L\}.
\end{array}
\end{equation}
However, when the number of scenarios $m$ is large, such as in the millions, the above problem is of large-scale with $O(n+mL)$ number of variables and $O(mL)$ number of constraints. Consequently, standard off-the-shelf interior-point-based  solvers may not be able to handle it. 
In such cases, one may directly tackle the compact form of problem \eqref{eq:source_problem} with $n$ variables through (stochastic) subgradient-based first-order methods. In this approach, only scenarios from the $(1-\tau_\ell)$-tail of $\left\{g^\ell(x;\omega^j)\right\}_{j=1}^m$ contribute   a non-zero term to the subgradient of $\superquantile_{\tau_\ell}{G^\ell(\,\cdot\,; \omega^{[m]})}$. 
This differs from the risk-neutral expected value of $\left\{g^\ell(x;\omega^j)\right\}_{j=1}^m$, where any batch of scenarios can be used to construct an unbiased estimate of the subgradient.
To address these computational hurdles, recent research has suggested 
to take  larger mini-batches to increase the likelihood of including a tail scenario in the batch estimate \citep{yang2021robust}, or to adopt a biased sampling strategy \citep{tamar2015sample,levy2020large}.
While these methods maintain scalability, they do not fully exploit the inherent sparsity structure offered by superquantile operators, which suggest that only a small subset of scenarios is essential for computation. This issue similarly applies to the smoothing approach of the superquantile operators, where the nonsmooth term $\max(\,\cdot\,, 0)$ in \eqref{cvar}   is approximated by a smooth function such as the log-sum-exponential \citep{laguel2021superquantiles, laguel2022superquantile}.
The smoothing approach makes all scenarios relevant in the subgradient computation of $\superquantile_{\tau_\ell}{G^{\,\ell}(\,\cdot\,;\omega^{[m]})}$, detracting from the sparsity that could otherwise be leveraged.

Another line of research for solving problem \eqref{eq:source_problem} may be classified as \emph{scenario reduction} strategies. Recall that $m$ represents the total number of collected scenarios. 
The concept of scenario reduction emerges from the observation that \eqref{eq:source_problem} should have an equivalent formulation with a reduced number of scenarios $m_\ell^*\leq m$ for each constraint. Here, the $m_\ell^*$ number of ``effective" scenarios—those that make up the $(1-\tau_\ell)$-tail of ${G^\ell(x^*;\omega^{[m]})}$ at an optimal solution  $x^*$ of  \eqref{eq:source_problem} \cite[\eg see][]{rahimian2018effective}, are not known prior to solving problem \eqref{eq:source_problem}. 
This observation leads to a fixed-point approach for iteratively refining the set of active scenarios by excluding those outside the tail, thus leveraging the sparsity inherent in the superquantile operator. 
Methods that eliminate inactive scenarios based on a candidate solution are known as  scenario reduction methods and resemble active-set methods in  the optimization literature \citep{roysetOA2011,arpon2018scenario}.
The \emph{outer approximation} (OA) scheme, as outlined in \citet[Chapter 6.C]{royset2021primer}, provides a structured framework for these strategies, allowing  the dynamic expansion and pruning of the candidate active scenario set.
An advantage of these methods is that the subproblems can be solved in a significantly reduced space, potentially small enough to facilitate the use of off-the-shelf interior-point-based solvers. 
However, a major challenge in implementing effective scenario strategies lies in the iterative update of the active scenario set,  which, based on our experience, often requires manual tuning to achieve good performance.

The goal of this paper is to introduce a scalable and robust second-order computational solver specifically designed to efficiently solve problem \eqref{eq:source_problem} at both low and high accuracy levels. Our solver implements a semismooth-Newton-based  augmented Lagrangian method (ALM) with a carefully designed linear system solver. We have opted for this method for two primary reasons.
(i) Our algorithm is designed to calculate the superquantile terms $\{\superquantile_{\tau_\ell}{G^\ell(x; \omega^{[m]})}\}_{\ell=0}^L$  at every iteration, instead of relying on a subset of scenarios for estimation. This precise computation is especially critical when 
$\tau$ is close to $1$ to prevent biased estimations. However, this computation, which involves the evaluation of all $\{g^\ell(x;\omega^j)\}_{j=1}^m$,    could be expensive when $m$ is large. To reduce the number of such evaluations, we adopt a second-order method so that the total number of iterations is usually small.
(ii) To identify critical scenarios that belong to the $(1-\tau_\ell)$-tail of $G^\ell(x; \omega^{[m]})$, it is sufficient to compute the function values of $\{g^\ell(x;\omega^j)\}_{j=1}^m$ rather than their more computationally demanding gradients or Hessians. 
In each subproblem, the generalized Jacobian matrix in the semismooth Newton method use these function values to  identify scenarios likely to influence the optimal $(1-\tau_\ell)$-tail. This number may be slightly larger than $(1-\tau_\ell)m$ since there may exist ties for the $(1-\tau_\ell)m$-th largest value of $\{g^\ell(x;\omega^j)\}_{j=1}^m$,  but it is usually much smaller than the total number of scenarios $m$.  The inversion of the (generalized) Hessian matrix can thus be done in a small dimensional space using the Sherman-Morrison-Woodbury formula.
{\bf As a consequence,  our method can be viewed as an implicit scenario reduction approach that uses second-order variational analysis to perform the tuning automatically}.

Our method draws inspiration from recent advances in second-order methods for large-scale optimization, where similar frameworks have been utilized to exploit various forms of structured sparsity, such as in risk-neutral and risk-averse linear regression using $\ell_1$-norm or $\ell_1$-$\ell_2$ norm regularization
\citep{li2018highly, li2018fused, zhang2020hessian,wu2022convex}.
These studies employ the ALM on the dual problem, leveraging the structured sparsity of the solution to reduce computational costs. In contrast, our method applies the ALM to the primal problem, utilizing the sparsity introduced by superquantiles to achieve efficiency.

Numerical experiments on randomly generated problems with affine random functions $G^{\,\ell}(\,\cdot\,;\omega^{[m]})$ indicate that the ALM method outperforms general-purpose solvers such as Gurobi and OSQP, and a specialized OA method for small instances. In larger scenarios with multiple superquantile constraints or high risk aversion, the performance improvements are substantial, sometimes exceeding 50 times  of the aforementioned solvers.
In addition, we find that the ALM method compares favorably to the specialized Portfolio Safeguard suite in terms of both reliability and speed.
A significant advantage of the ALM framework is its ability to efficiently utilize warm-start solutions, enabling a highly efficient process for solving a series of similar problems, such as computing the solution path of \eqref{eq:source_problem} across various confidence levels $\{\tau_\ell\}_{\ell=0}^L$.
In the context of large-scale quantile regression, ALM proves adept at predicting the full distribution of a response variable across a spectrum of confidence levels. This approach effectively reuses prior solutions, with many  effective scenarios in one problem remaining so in subsequent ones, leading to improvements of up to 20 times over existing techniques.

The rest of the paper is organized as follows.
In \cref{sec:preliminaries}, we introduce some concepts from nonsmooth analysis and review known results related to the projection onto the top-$k$-sum constraint.
In \cref{sec:methodology}, we outline our computational framework to solve \eqref{eq:source_problem} using the augmented Lagrangian method, and discuss how to solve the subproblems via the semismooth Newton method. In \cref{sec:scenario_reduction},  we  describe the structure of the  resulting linear systems  from the semismooth Newton method and examine the computational cost associated with solving these systems.
Section \ref{sec:experiments} presents extensive numerical results from both synthetic data and real datasets. The paper ends with a concluding section that summarizes our findings and discusses  future research directions.

\subsubsection*{Notation}

For a matrix $A$, the submatrix formed by the rows in an index set $\mathcal{I}$ and the columns in an index set $\mathcal{J}$ is denoted $A_{\mathcal{I},\mathcal{J}}$. We use ``$:$'' to denote MATLAB notation for index sets, \eg $A_{\mathcal{I},:}$ represents the submatrix of $A$ formed by rows $\mathcal{I}$ and all columns.
The cardinality of a set $\mathcal{I}$ is denoted by $\abs{\mathcal{I}}$.
The vector of all ones in dimension $m$ is denoted by $\ind_m$; for an index set $\mathcal{I}\subseteq\{1,\ldots,m\}$, $\ind_{\mathcal{I}}$ denotes the vector with ones in the indices corresponding to $\mathcal{I}$ and zeros otherwise when the ambient dimension is clear; and $e^i$ denotes the $i$-th standard basis vector.
For a vector $x\in\R^m$, $x_i$ denotes the $i$-th element of $x$, $x_{(i)}$ denotes the $i$-th largest element of $x$, and $\sorth{x}=x_\pi$ denotes a nonincreasing rearrangement of the elements of $x$ according to permutation $\pi$ with the convention that $\sorth{x}_0\coloneqq +\infty$ and $\sorth{x}_{m+1}\coloneqq -\infty$.
For any $x\in\R^m$ and positive integer $1\leq k\leq m$, there exist integers $k_0,k_1$ satisfying $0\leq k_0\leq k-1$ and $k\leq k_1\leq m$ such that
\begin{equation}
  \label{eq:order_structure}
  \sorth x_1 \geq\cdots\geq \sorth x_{k_0} > \sorth x_{k_0+1} =\cdots= \sorth x_k =\cdots= \sorth x_{k_1} > \sorth x_{k_1+1} \geq\cdots\geq \sorth x_{m},
\end{equation}
with the convention that $k_0=0$ if $\sorth x_1=\sorth x_k$ and $k_1=m$ if $\sorth x_m = \sorth x_k$.
We call $(k_0,k_1)$ satisfying \eqref{eq:order_structure} the \emph{index-pair of $x$ associated with $k$} and define the corresponding partition of $\{1, \ldots, m\}$ as
\begin{equation}\label{eq:index}
\alpha\coloneqq\{1,\ldots,k_0\},\quad  \beta\coloneqq\{k_0+1,\ldots,k_1\},  
\quad \gamma\coloneqq\{k_1+1\ldots,m\}.
\end{equation}
For any positive integer $k(\leq m)$, we denote
the \emph{top-$k$-sum operator} of $x\in \mathbb{R}^m$ as 
\[
\maxsum_k{x}\coloneqq x_{(1)}+\cdots+x_{(k)}.
\]
In addition, we denote the set of all non-positive top-$k$-sum vectors as
\begin{equation}\label{eq:defn set b}
    \mathcal{B}_k\coloneqq\{x\in \mathbb{R}^m:\maxsum_k{x}\leq 0\}.
\end{equation}

For a closed convex set $S\subseteq\R^n$, its  indicator function at a point $x\in \mathbb{R}^n$ is denoted as $\delta_{S}(x)$, which takes the value 0 if $x\in S$ and $+\infty$ otherwise; the Euclidean projection of $x$ onto $S$ is written as $\proj_{S\,}{x} = \argmin_{y\in S} \|x-y\|_2$.
The normal cone of a convex set $S$ at $x\in S$ is denoted by ${\cal N}_S(x)\coloneqq\{d:\inner{d}{x-y}\geq0,\;\forall y\in S\}$.
The Fenchel conjugate of a function $f$ is denoted by $f^*$.

\section{Preliminaries}
\label{sec:preliminaries}
The purpose of this section is to review two known results that will be used in our subsequent development of the ALM framework.
The first result addresses the computation of the projection onto the set ${\cal B}_k$ defined in \eqref{eq:defn set b}. The second result further characterizes the Fr{\'e}chet differentiability of this projection, based on which we are able to derive elements from its generalized Jacobian.

We begin by reviewing some concepts from variational analysis. Interested readers are referred to \citet[Chapters 4 \& 7]{facchinei2007finite} for more detail. Let $\mathcal{O}\subseteq\mathbb{R}^m$ be an open set. 
For any locally Lipschitz function $F:\mathcal{O}\to\mathbb{R}^n$, it follows from the
Rademacher theorem  that $F$ is Fr{\'e}chet-differentiable almost everywhere.
We denote the set of Fr{\'e}chet-differentiable points of $F$ as $\mathcal{D}_F \coloneqq \{x\in\mathcal{O} : \text{$F$ is Fr{\'e}chet-differentiable at $x$}\}$, and  the Jacobian of $F$ at $x\in\mathcal{D}_F$ by $JF(x)$. In addition, we write the Bouligand subdifferential of $F$ at $x\in\mathcal{O}$ as 
\[
\partial_B F(x)\coloneqq\bigl\{ V\in \mathbb{R}^{n\times m} : V = \lim_{k\to\infty} JF(x^k),\; x^k\in\mathcal{D}_F,\;x^k\to x\bigr\},
\] 
and the Clarke generalized Jacobian of $F$ at $x\in\mathcal{O}$ as $\partial F(x) \coloneqq \mathsf{conv}\{\partial_B F(x)\}$, \ie the convex hull of the Bouligand subdifferential at $x$. The function $F$ is said to be directionally differentiable at $x\in\mathcal{O}$ along the direction $d$ if 
\[
F^\prime(x;d) \coloneqq \lim_{t\downarrow 0} \frac{F(x+td) - F(x)}{t}
\]
exists. We say $F$ is directionally differentiable at $x$ if it is directionally differentiable along any direction. Under this condition, 
 if it further holds that  for any $\Delta x\to0$ and $V\in\partial F(x+\Delta x)$, $F(x+\Delta x) - F(x) - V\Delta x =o(\|{\Delta x}\|_2)$ ($O(\|{\Delta x}\|_2^2)$), we say $F$ is semismooth (strongly semismooth) at $x\in\mathcal{O}$.
Finally, $F$ is (strongly) semismooth  on $\mathcal{O}$ if it is (strongly) semismooth  for any $x\in \mathcal{O}$. %

\subsection{Projection onto the set \texorpdfstring{$\mathcal{B}_k$}{Bk}}
Recall the definition of the set $\mathcal{B}_k$  in \eqref{eq:defn set b}.
It is easy to see that when $(1-\tau_\ell)m$ is an integer, %
the constraint $\superquantile_{\tau_\ell}[\big]{G^{\,\ell}(x;\omega^{[m]})} \leq 0$ in \eqref{eq:source_problem} can be equivalently written as $G^\ell(x;\omega^{[m]})\in \mathcal{B}_{k_\ell}$ where $k_\ell=(1-\tau_\ell)m$.
Given a vector $y\in \mathbb{R}^m$, its projection onto $\mathcal{B}_k$, denoted by 
\begin{align}\label{eq:maxksum_projection}
    \tilde{y} \coloneqq \proj_{\mathcal{B}_k}{y} = \displaystyle\operatornamewithlimits{argmin}_{y^\prime\in \mathbb{R}^m}\bigl\{\tfrac12 \norm{y^\prime-y}_2^2 : y^\prime \in \mathcal{B}_k\bigr\}, 
\end{align}
can be computed in $O(m\log m)$ operations \citep{wu2014moreau,davis2015algorithm,li2021fast,roth2023algorithms}.
For special cases when $k=1$ or $k=m$,  the projection has a closed form and can be computed in $O(m)$ operations.
Excluding the two special cases, all of these computations use the fact that the top-$k$-sum operator $\mathopfont{T}_k$ is permutation-invariant. Therefore,  one may instead first sort $y$ to $y_\pi$ such that $(y_\pi)_1 \geq\cdots\geq (y_\pi)_{m}$,  and then solve the following problem:
\begin{align}\label{eq:maxksum_projection_sorted}
    \bar{y} \coloneqq \Argmin_{y\in\R^m}\bigl\{\tfrac12 \norm{y-y_\pi}_2^2 : y \in \mathcal{B}_k\bigr\} = \Argmin_{y\in\R^m}\left\{\tfrac12 \norm{y-y_\pi}_2^2 : y_1\geq y_2\geq \cdots \geq y_m, \sum_{i=1}^k y_i \leq 0\right\}. 
\end{align}
In this case, the solution $\bar{y}$ must satisfy the following Karush–Kuhn–Tucker (KKT) optimality conditions for some unique indices $\bar k_0$ and $\bar k_1$ and  some multipliers $(\bar\mu, \bar\lambda,  \bar{\theta})\in \mathbb{R}^{m}\times \mathbb{R}\times \mathbb{R}$:
\begin{align}
  \label{eq:kkt_maxksum_projection_sort_2}
  \begin{array}{r*{1}{wc{0mm}}lcr*{1}{wc{0mm}}ll}
    \bar y_{\bar\alpha} &=& \sorth y_{\bar\alpha} - \bar\lambda\ind_{\bar\alpha} &\quad & \bar \mu_{\bar\alpha} &=& \ind_{\bar\alpha}\\*[0.05in]
   \bar y_{\bar\beta} &=& \sorth y_{\bar\beta} - \bar\lambda\bar{\mu}_{\bar\beta}  = \bar\theta\ind_{\bar\beta}  & \quad & \bar \mu_{\bar\beta} &\in& \{\mu \in \R^{\abs{\bar\beta}}: \mu\in[0,\ind],\;\ind^\top\mu = k-\bar{k}_0\}\\*[0.05in]
    \bar y_{\bar\gamma} &=& \sorth y_{\bar\gamma}  &\quad & \bar \mu_{\bar\gamma} &=& 0\\*[0.05in]
    0 &=& \ind_{\bar\alpha}^\top \sorth y_{\bar\alpha} - \bar k_0\bar \lambda + (k-\bar k_0)\bar\theta &\quad & \bar\lambda &>& 0  \qquad  \bar y_{\bar k_0} > \bar\theta > \bar y_{\bar k_1+1},
  \end{array}
\end{align}
where the sets $\{\bar{\alpha}, \bar{\beta}, \bar{\gamma}\}$ are defined according to \eqref{eq:index} based on the index-pair $(\bar{k}_0, \bar{k}_1)$. 
The pivoting procedures in \cite{roth2023algorithms} compute a (sorted) solution  to \eqref{eq:kkt_maxksum_projection_sort_2} by identifying the optimal index-pair $(\bar{k}_0,\bar{k}_1)$ satisfying the order structure \eqref{eq:order_structure}.
It was shown there that the total cost to compute $\bar{y}$ is $O(m+\bar{k}_1\log(m))$, which is achieved by integrating the sorting procedure into the pivoting procedure after an initial heapsort construction, allowing for the dynamic computation of $\bar{k}_1$.
This approach makes use of the max-heap data structure (\citealp[\eg see][]{williams1964heap}), which can be constructed with a cost of $O(m)$, and offers substantial speed advantages over a full sort especially when $k\ll m$.

\subsection{Differentiability of the projection map \texorpdfstring{$\mathopfont{proj}_{{\cal B}_k}$}{ProjBk}}
In this section, we discuss the differentiability  of the projection map $\mathopfont{proj}_{{\cal B}_k}$, which is needed in the construction of the generalized Jacobian of this projection. Our derivations are adapted from the work \cite{wu2014moreau}, which focused on the projection onto the vector $k$-norm ball.

Based on the definition of the Fr{\'e}chet differentiability, the projection map $\mathopfont{proj}_{{\cal B}_k}$ is differentiable at a point $y\in\R^m$ when its one-sided directional derivative $\proj'_{\mathcal{B}_k}{ y;\,\cdot\,}$ is linear in the second argument.
The latter directional derivative can be calculated by projecting the direction onto the critical cone of $\mathcal{B}_k$ at $y$ \cite[Theorem 4.1.1]{facchinei2007finite}.
Our subsequent derivations are based on these observations.

For convenience, we will consider a sorted point $y\in\R^m$ with the  permutation $\pi$ such that $y=y_\pi=\sorth{y}$.
In general, one may first obtain the permutation $\pi$ and then apply the following derivations. 
In the sorted case, the critical cone of ${\cal B}_k$ at $y$ is given by
\begin{align*}
{\cal C}_{{\cal B}_k}(y)&\coloneqq\tangentcone[]{\mathcal{B}_k}{\bar{y}}\cap(\bar{y}-y)^\perp \quad  \mbox{  with $\bar{y}=\proj_{\mathcal{B}_k}{y}$},
\end{align*}
where  
$\tangentcone[]{\mathcal{B}_k}{\bar{y}}$ is the tangent cone of $\mathcal{B}_k$ at $\bar y$ and $v^\perp$ represents the orthogonal space of a vector $v$.
The KKT conditions \eqref{eq:kkt_maxksum_projection_sort_2} of the sorted projection problem \eqref{eq:maxksum_projection_sorted} gives
\begin{align*}
    (\bar{y}-y)^\perp
    &%
    =\{d\in\R^m:\ind_{\bar\alpha}^\top d_{\bar\alpha}+\bar\mu_{\bar\beta}^\top d_{\bar\beta}=0\},\quad \bar\mu_{\bar\beta}\in[0,1]^{\abs{\bar\beta}},\;\ind^\top \bar\mu_{\bar\beta}=k-k_0,
\end{align*}
where $\{\bar\alpha,\bar\beta,\bar\gamma\}$ are the partition of the indices associated with $\bar y$ in \eqref{eq:index}.
This holds since $\bar{y} = y - \bar\lambda\bar\mu$ for $\bar{\mu}\in\partial\maxsum_k{\bar{y}}$ and $\bar\lambda>0$.
By direct computation, the tangent cone at any  $\bar{y}\in {\cal B}_k$ is given by 
\begin{align*}
    \tangentcone{\mathcal{B}_k}{\bar{y}}  = 
    \begin{cases}
        \R^m,\quad&\text{if $\maxsum_k{\bar{y}} < 0$},\\
        \bigl\{d\in\R^m : \ind_{\bar\alpha}^\top d_{\bar\alpha} + \maxsum_{k-\bar{k}_0}[\big]{d_{\bar\beta}}\leq0\bigr\},\quad&\text{if $\maxsum_{k}{\bar{y}} = 0$}.
    \end{cases}
\end{align*}
Therefore, when $\maxsum_k{y}>0$ so that $\maxsum_k{\bar{y}}=0$, the critical cone can be expressed as
\begin{align}
\label{eq:criticalcone_analytic}
{\cal C}_{{\cal B}_k}(y)&=\bigl\{d\in\R^m:\ind_{\bar\alpha}^\top d_{\bar\alpha} + \maxsum_{k-\bar{k}_0}[\big]{d_{\bar\beta}}\leq0,\;\ind_{\bar\alpha}^\top d_{\bar\alpha}+\bar\mu_{\bar\beta}^\top d_{\bar\beta}=0 \bigr\}\notag\\*
     &=\bigl\{d\in\R^m:\maxsum_{k-\bar{k}_0}[\big]{d_{\bar\beta}}\leq \bar\mu_{\bar\beta}^\top d_{\bar\beta},\;\ind_{\bar\alpha}^\top d_{\bar\alpha}+\bar\mu_{\bar\beta}^\top d_{\bar\beta}=0 \bigr\}.
\end{align}
When $\maxsum_k{y}=0$ so that $y=\bar{y}$, we have $y-\bar{y} = \{0\}$ and ${\cal C}_{{\cal B}_k}(y)$ coincides with the tangent cone $\tangentcone[]{\mathcal{B}_k}{\bar{y}}$. %
The next result states the conditions under which the metric projector is differentiable, \ie when the critical cone ${\cal C}_{{\cal B}_k}(y)$ is a linear subspace.
It is readily obtained from conditions (i), (iii), and (iv) of \citet[Theorem 4.2]{wu2014moreau}.

\begin{lemma}\label{prop:differentiability}
    Let $1\leq k\leq m$, $y=\sorth y\in\R^m$ be a given sorted vector, and $\bar{y}$ be its projection with associated index-pair $(\bar{k}_0, \bar{k}_1)$ in \eqref{eq:order_structure} and index-sets $\{\bar{\alpha}, \bar{\beta}, \bar{\gamma}\}$ in \eqref{eq:index}.
    Let $(\bar y, \bar\mu, \bar{\lambda}, \bar{\theta})$ satisfy the optimality conditions in \eqref{eq:kkt_maxksum_projection_sort_2}.
    Then  $\proj_{\mathcal{B}_k}{\,\cdot\,}$ is differentiable at $y$ if and only if $y$ satisfies one of the following three properties:
    \begin{enumerate}
        \item $\maxsum_k{y}<0$;
        \item $\maxsum_k{y}>0$, $k<\bar{k}_1$, and ${y}_{\bar{k}_0+1}-\bar\lambda < \bar\theta < {y}_{\bar{k}_1}$;
        \item $\maxsum_k{y}>0$ and $k=\bar{k}_1$.
    \end{enumerate}
\end{lemma}

Next we use \cref{prop:differentiability} to compute $J_{\mathcal{B}_k}$, the Jacobian of the map $\mathopfont{proj}_{{\cal B}_k}$ at differentiable points corresponding to the three cases in the above lemma.
For simplicity, we assume that $y=\sorth{y}$ is sorted. For any positive integer $p$, we write $I_p$ as the $p\times p$ identity matrix. The proof of the following lemma, which is adapted from \citet[][eq.~(68) in Section 4.3.4]{wu2014moreau}, is given in \cref{apx:proof}.

\begin{lemma}
\label{prop:explicit_gen_jac}
    Let $1\leq k\leq m$, $y=\sorth y\in\R^m$ be a given sorted vector and $\bar{y}$ be its projection with associated index-pair $(\bar{k}_0, \bar{k}_1)$ in \eqref{eq:order_structure} and index-sets $\{\bar{\alpha}, \bar{\beta}, \bar{\gamma}\}$ in \eqref{eq:index}.  
 
{\rm (a)} If $\maxsum_k{y}<0$, then $\proj'_{\mathcal{B}_k}{y;d}=d$ for any $d\in \mathbb{R}^m$ and $J_{\mathcal{B}_k}\in \mathbb{R}^{m \times m}$ is the identity matrix.

{\rm (b)} If $\maxsum_k{y}>0$, $k<\bar{k}_1$ and  ${y}_{\bar{k}_0+1}-\bar\lambda < \bar\theta < {y}_{\bar{k}_1}$.
		Then
		\begin{align}
			\label{eq:generalized_jacobian}
			J_{\mathcal{B}_k} =
			\begin{bmatrix}
				I_{|\bar\alpha| + |\bar{\beta}|} \, - B^\top (B B^\top)^{-1} B\\%_{\pi^{-1},\pi^{-1}}
				& I_{|\bar\gamma|}
			\end{bmatrix},
		\end{align}
where 
		\begin{align}
			\label{eq:criticalcone_projection_case2}
		B \coloneqq
			\begin{bmatrix}
		\ind_{\bar\alpha}^\top & \bar\mu_{\bar\beta}^\top\\
				0_{(\abs{\bar\beta}-1),\abs{\bar\alpha}} & C
			\end{bmatrix} \in \R^{(\abs{\bar\alpha} + \abs{\bar\beta}-1)\times (\abs{\bar\alpha}+\abs{\bar\beta})},\quad
			C \coloneqq
			\begin{bmatrix}
			1 & -1\\
			&\ddots & \ddots\\
			&& 1 & -1
			\end{bmatrix}\in\R^{(\abs{\bar\beta}-1)\times \abs{\bar\beta}}
		\end{align}
  and the matrix $ B B^\top$ is invertible.
       
{\rm (c)} If $\maxsum_k{y}>0$ and $k=\bar{k}_1$, then
        \begin{align*}
            J_{\mathcal{B}_k} = \begin{bmatrix}
                I_{|\bar\alpha| + |\bar{\beta}|}\;  - \displaystyle\frac{ c  c^\top}{\| c\|_2^2} & \\ & I_{|\bar\gamma|}
            \end{bmatrix}
           \quad \mbox{with ${c} \coloneqq \ind_{(\bar\alpha,\bar\beta)}$}.
        \end{align*}
\end{lemma}
Finally, since the locally Lipschitz continuous map $\projection_{\mathcal{B}_k}$ is differentiable almost everywhere, any point $y\in \mathbb{R}^m$ can be viewed as the limit of a sequence of differentiable points.
Hence,  one can further derive the generalized Jacobian of $\projection_{\mathcal{B}_k}$ at any point based on \cref{prop:explicit_gen_jac}.

\section{The algorithmic framework}
\label{sec:methodology}

In this section, we introduce our algorithmic framework to solve  the superquantile-based problem \eqref{eq:source_problem}. Since one can always move the superquantile term   $\superquantile_{\tau_0}[\big]{G^{\,0}(x;\omega^{[m]})}$ in the objective function to the constraint set using the epigraphic transformation, we assume for simplicity that $G^{\,0}(x;\omega^{[m]}) = \{0\}$ in \eqref{eq:source_problem} so that the objective function only consists of the convex smooth term $f(x)$. By introducing auxiliary variables $\{y^{\ell}\}_{\ell=1}^L\subseteq \mathbb{R}^m$ and $z\in \mathbb{R}^n$, one can reformulate problem \eqref{eq:source_problem} into 
\begin{equation}
    \label{eq:source_problem_smooth}
    \begin{array}{rl}
         \displaystyle \min_{x,y,z} & f(x) \\
         \text{s.t.} & y^\ell \geq \left(g^{\,\ell}(x;\omega^{1}), \ldots, g^{\,\ell}(x;\omega^{m}) \right)^\top,\;y^\ell\in {\cal B}_{k_\ell}, \quad \ell\in\{1,\ldots,L\},\\[0.05in]
         & z = x, \quad z\in X,
    \end{array}
\end{equation}
where $k_\ell = (1-\tau_\ell)m$
and ${\cal B}_k$ is the non-positive top-$k$-sum set defined in \eqref{eq:defn set b}.
For exposition, we consider $X=[p,q]\subseteq\R^n$ to be a box with $p\leq q\in\{\R\cup\pm\infty\}^n$.
Our overall algorithm has double loops, where the outer loop uses the ALM to deal with the constraints, and the inner loop adopts the semismooth Newton method to solve the ALM subproblems.

\subsection{An augmented Lagrangian framework to solve problem \texorpdfstring{\eqref{eq:source_problem_smooth}}{Src}}
\label{sec:p_alm_framework}

For notational simplicity, we write $y = ((y^1)^\top, \ldots, (y^L)^\top)^\top \in \mathbb{R}^{mL}$, 
\[
G(x) \coloneqq \left(\underbrace{g^{1}(x;\omega^1), \ldots, g^1(x;\omega^m)}_{G^1(x;\omega^{[m]})}, \ldots, \underbrace{g^{L}(x;\omega^1),\ldots, g^{L}(x;\omega^m)}_{G^L(x;\omega^{[m]})} \right)^\top: \mathbb{R}^n \to \mathbb{R}^{mL},
\]
and $\overline{\mathcal{B}} \coloneqq \mathcal{B}_{k_1}\times\cdots\times\mathcal{B}_{k_L}\subseteq \mathbb{R}^{mL}$.  Let  $\lambda (\geq 0)\in\R^{mL}$ be the multiplier of the inequality constraint $y\geq G(x)$ and 
 $\mu\in\R^n$ be the multiplier of the equality constraint $z=x$. 
The augmented Lagrangian function associated with problem \eqref{eq:source_problem_smooth} for a parameter $\sigma>0$ takes the form of
\begin{align*}
    L_\sigma(x,y,z;\lambda,\mu)
	&\coloneqq f(x) + \frac\sigma2\norm[\big]{ \max\bigl(G(x)+\lambda/\sigma - y,0\bigr)}_2^2 + \frac\sigma2\norm[\big]{x+\mu/\sigma - z}_2^2+\texttt{const}
\end{align*}
for $x\in \mathbb{R}^n$, $y\in \overline{\mathcal{B}}$ and $z\in X$, 
where \texttt{const} is a constant depending on the multipliers.
To make the subproblems easier to solve,  at the $\nu$-th iteration,  we may also add a proximal term $\displaystyle\frac{1}{2\sigma}\norm{x-x^\nu}_M^2\coloneqq\frac{1}{2\sigma}(x-x^\nu)^\top M(x-x^\nu)$ for some  positive semidefinite matrix $M$ to regularize the augmented Lagrangian function. 
The overall algorithm is summarized in \cref{alg:palm}.

\begin{algorithm}
	\caption{The ALM framework for problem \eqref{eq:source_problem_smooth}}\label{alg:palm}
	\begin{algorithmic}[1]
		\Initialize{Choose $x^0=z^0\in\R^n$, $y^0,\lambda^0\in\R^{mL}$, and $\mu\in\R^n$.
    Choose $\sigma_0>0$  and a positive semidefinite matrix $M$.
		Set iteration $\nu=0$, and repeat the following steps until a proper stopping criterion is satisfied.}
		\State Primal update:
            \begin{align}\label{eq:palm_subproblem}
            (x^{\nu+1}, y^{\nu+1}, z^{\nu+1}) \approx \argmin_{x,\,y\in \overline{{\cal B}},\,z\in X}  L_{\sigma_\nu}(x,y,z;\lambda^\nu,\mu^\nu) + \frac{1}{2\sigma_\nu}\norm{x-x^\nu}_{M}^2.
		\end{align}

		\State Dual update:
            \begin{align*}
                \lambda^{\nu+1}&\coloneqq \max\Bigl(\lambda^\nu + \sigma_\nu \bigl(G(x^{\nu+1}-y^{\nu+1})\bigr),0\Bigr)\\
                \mu^{\nu+1}&\coloneqq \mu^\nu + \sigma_\nu (x^{\nu+1} - z^{\nu+1}).
		\end{align*}

            \State Check termination criteria and update $\sigma_\nu$.
    
	\end{algorithmic}
\end{algorithm}

Let $J G(x)\in \mathbb{R}^{(mL) \times n}$ denote the Jacobian matrix of $G$ at $x$. 
The KKT optimality condition of problem \eqref{eq:source_problem_smooth} is  %
\begin{equation}
\left\{
\begin{aligned}
    \label{eq:source_problem_smooth_kkt}
    0 &= \nabla f(x) + J G(x)^\top \lambda + \mu\\
    \lambda &\in {\cal N}_{\overline{\mathcal{B}}}(y), \quad 
    \mu \in   {\cal N}_X(z)\\
    0&\leq \lambda \perp y-G(x) \geq 0, \quad 
  x-z = 0.
\end{aligned}
\right.
\end{equation}
The above condition is necessary and sufficient to characterize a solution to the convex problem \eqref{eq:source_problem_smooth} when the Slater condition holds.

It is known that under mild conditions,  any accumulation point of the sequence $\{(x^\nu, y^\nu, z^\nu, \lambda^\nu, \mu^\nu)\}$ generated by the ALM  converges to a tuple satisfying the KKT condition \eqref{eq:source_problem_smooth_kkt} \citep[\eg see][]{rockafellar1976augmented}.
Conditions ensuring the local superlinear convergence rate of ALM methods are  also well studied in the literature; see the early seminal works by \citet{rockafellar1976augmented,rockafellar1976monotone} and the later refinement by \cite{cui2019superlinear}. 
For example,  consider the following constraint non-degeneracy condition at a tuple $(\bar{x}, \bar{y}, \bar{z})$ due to \citet{robinson1984local,robinson1987local}:
\begin{align*}
    0&\in\mathopfont{int}\bigl\{F(\bar{x},\bar{y},\bar{z}) + JF(\bar{x},\bar{y},\bar{z})\left(\R^n \times \R^{mL}\times \R^n\right) - \R^{mL}_+\times0_n\times\overline{\mathcal{B}}\times X\bigr\},
\end{align*}
where     $F(x,y,z)\coloneqq \bigl(y-G(x),z-x,y,z \bigr)$,  $\mathopfont{int}$ denotes the interior of a set and $\mathbb{R}^{mL}_+$ denotes the $(mL)$-dimensional non-negative orthant.
See \cref{sec:cq_non-degeneracy} for a discussion of the non-degeneracy assumption based on a simplified form of problem \eqref{eq:source_problem_smooth}.
The convergence result in \cite{rockafellar1976augmented} implies that if this condition is satisfied at an optimal solution $(x^*, y^*, z^*)$, the sequence of  dual variables $\{(\lambda^\nu, \mu^\nu)\}$ generated by the ALM converges at an asymptotically superlinear rate to the unique multiplier  when subproblems are solved to a sufficient accuracy. A proper termination condition to ensure the rate is detailed in the numerical experiments (\cref{sec:experiments}).

\subsection{A semismooth Newton method for  subproblem \texorpdfstring{\eqref{eq:palm_subproblem}}{}}
\label{sec:subproblem_alm}

The  effectiveness of the ALM for solving problem \eqref{eq:source_problem} hinges on the efficient computation of the subproblem \eqref{eq:palm_subproblem}. 
In this section, we present a semismooth Newton method to solve this subproblem. Despite the seemingly  high computational cost typically associated with such procedures, the actual cost for the superquantile-constrained problem is in fact quite low. This efficiency is largely attributed to the unique structural properties brought by the projection onto the non-positive top-$k$-sum set $\overline{\cal B}$.

For exposition, we begin by considering the  augmented Lagrangian subproblem associated with the (possibly nonconvex) equality constraint $y = G(x)$ before returning to  the convex inequality case.
For fixed ${\lambda^\nu}$ and $\mu^\nu$, consider the subproblem
\begin{align}\label{eq:alm_subproblem_2}
 (x^{\nu+1},y^{\nu+1},z^{\nu+1}) \in \Argmin_{{x,\,y\in\overline{\mathcal{B}},\,z\in X}} \Big\{ f(x)
    &+ \frac{\sigma_\nu}2 \norm[\big]{y-\bigl(G(x)+\sigma_\nu^{-1}{\lambda^\nu}\bigr)}_2^2\notag\\*
	&+ \frac{\sigma_\nu}2 \norm[\big]{z-(x+\sigma_\nu^{-1}\mu^\nu)}_2^2
        + {\frac{1}{2\sigma_\nu}\norm{x-x^\nu}_M^2 \Big\}}
\end{align}

for some positive semidefinite matrix $M$ which may help to regularize the subproblem.
Observe that the unique solution in terms of $(y,z)$ satisfies
\begin{align}
        \label{eq:yz_projection}
		y^{\nu+1} = \proj_{\overline{\mathcal{B}}}[\big]{G(x^{\nu+1})+\sigma_\nu^{-1}{\lambda^\nu}},\quad
		z^{\nu+1} = \proj_{X}[\big]{x^{\nu+1}+\sigma_\nu^{-1}{\mu^\nu}}.
\end{align}

Substituting the above solutions into \eqref{eq:alm_subproblem_2} for $y$ and $z$ gives the problem in terms of $x$ alone:
\begin{align}
    \label{eq:varphi}
    x^{\nu+1} \in \Argmin_x \Big\{\varphi_{\sigma_\nu}(x)\coloneqq f(x) 
 &+ \frac{\sigma_\nu}{2} \norm[\big]{G(x)+\sigma_\nu^{-1}{\lambda^\nu} - \proj_{\overline{\mathcal{B}}}[\big]{G(x)+\sigma_\nu^{-1}{\lambda^\nu}}}_2^2\notag\\*
	&+\frac{\sigma_\nu}{2} \norm[\big]{x+\sigma_\nu^{-1}{\mu^\nu} - \proj_{X}[\big]{x+\sigma_\nu^{-1}{\mu^\nu}}}_2^2
        + {\displaystyle\frac{1}{2\sigma_\nu}} \norm{x-x^\nu}_M^2 \Big\}.
\end{align}
Therefore,  the function $\varphi_{\sigma_\nu}$ may serve as the subproblem objective, which has a Lipschitz continuous gradient that takes the following form:
\begin{align}
	\label{eq:varphi_gradient}
	\nabla\varphi_{\sigma_\nu}(x) = \nabla f(x) 
 &+ {\sigma_\nu} \, J G(x)^\top  \big[ G(x) + {\sigma}_\nu^{-1}\lambda^\nu - \proj_{\overline{\mathcal{B}}}{G(x) + {\sigma}_\nu^{-1}\lambda^\nu} \big]\notag\\*[0.1in]
	&+ {\sigma_\nu} \left[\,x + {\sigma}_\nu^{-1}\mu^\nu - \proj_{X}{x+\sigma_\nu^{-1}\mu^\nu}\,\right]+ \sigma_\nu^{-1} M(x-x^\nu).
\end{align}

Since $\overline{\mathcal{B}}$ and $X$ are polyhedral sets,  the projections onto them  are strongly semismooth \citep{shapiro1988sensitivity}. Therefore, one can apply the the semismooth Newton method to solve the nonsmooth equation $\nabla \varphi_{\sigma_\nu}(x) = 0$ in \eqref{eq:varphi_gradient} and use the line search for globalization.

The following lemma, whose proof is deferred to \cref{apx:proof}, shows that we may perform the same procedure to solve the augmented Lagrangian subproblem associated with \emph{inequality} constraints.
\begin{lemma}\label{lem:ineqalm}
    Let $g\in\R^m$ be given.
    Let ${\mathcal{B}}_k$ be the non-positive top-$k$-sum set defined in \eqref{eq:defn set b}.
    Consider the set
    \begin{equation}
    \label{eq:ineqalm_proj}
    \begin{array}{rl}
        \displaystyle \mathcal{P}\coloneqq\Bigl\{y^*:y^*\in\underset{y}{\argmin}\; \bigl\{\|\max(\sorth  g-y,0)\|_2^2:y\in {\mathcal{B}_k}\bigr\}\Bigr\}
    \end{array}
    \end{equation}
    and the element
    \begin{equation}
    \label{eq:eqalm_proj}
    \begin{array}{rl}
        \displaystyle \bar y=\underset{y}{\Argmin}\; \bigl\{\|\sorth  g -y\|_2^2:y\in {\cal B}_k\bigr\}.
    \end{array}
    \end{equation}
    Then it holds that $\bar y\in\mathcal{P}$.
    In fact, it holds that $\tilde{y}\in\mathcal{P}$ for any $\tilde{y}$ such that $(\tilde{y})_{1:\bar{k}_1}=(\bar{y})_{1:\bar{k}_1}$ with $(\tilde{y})_{\bar{k}_1+1:m}\leq(\bar{y})_{\bar{k}_1+1:m}$, where $\bar{k}_1$ associated with $\bar y$ is defined in \eqref{eq:order_structure}.
\end{lemma}
Therefore, given $\lambda^\nu$ and $\mu^\nu$, instead of \eqref{eq:varphi} for the equality constrained case, we may consider the following objective function 
\begin{align*}
    \widetilde\varphi_{\sigma_\nu}(x) \coloneqq f(x) &+ \frac{\sigma_\nu}{2}\norm[\Big]{\max\bigl(G(x) + \sigma_\nu^{-1}\lambda^\nu - \proj_{\overline{\mathcal{B}}}{G(x) + \sigma_\nu^{-1}\lambda^\nu},0 \bigr)}_2^2\\*
    \phantom{\varphi_{\sigma_\nu}^\tau(x) \coloneqq f(x)} &+ \frac{\sigma_\nu}{2}\norm[\big]{(x + \sigma_\nu^{-1}\mu^\nu) - \proj_{X}{x + \sigma_\nu^{-1}\mu^\nu} }_2^2 + \displaystyle{\frac 1{2\sigma_\nu}}\norm{x-x^\nu}_M^2.
\end{align*}
Observe that
if $x^+\in \Argmin_x \widetilde\varphi_{\sigma_\nu}(x)$, we may take
\begin{align*}
    y^{+} \coloneqq \proj_{\overline{\mathcal{B}}}{G(x^+) + \sigma_\nu^{-1}\lambda^\nu},\quad
    z^{+} \coloneqq \proj_{X}{x^+ + \sigma_\nu^{-1}\mu^\nu}
\end{align*}
to obtain a solution to the ALM subproblem
\begin{align*}
    (x^+,y^+,z^+) \in \argmin_{x,\,y\in\overline{\mathcal{B}},\,z\in X} \; f(x) &+ \frac{\sigma_\nu}{2}\norm[\big]{\max(G(x) + \sigma_\nu^{-1}\lambda^\nu - y,0)}_2^2 \\*
    &+ \frac{\sigma_\nu}{2}\norm[\big]{x+\sigma_\nu^{-1}\mu^\nu-z}_2^2 + \frac 1{2\sigma_\nu}\norm{x-x^\nu}_M^2.
\end{align*}
Though there may exist multiple solutions $y^+$, the value of the subproblem objective coincides with $\widetilde\varphi_{\sigma_\nu}(x^+)$ for any optimal $(y^+,z^+)$.
Therefore, the ALM framework can readily handle inequality constraints.
In the inequality constrained case, the first- and second-order information of the function $\varphi$ can be computed analogous to \eqref{eq:varphi_gradient} and \eqref{eq:varphi_hessian} by replacing $J G(x)$ with $I_{G^+}\cdot J G(x)$, where $I_{G^+}$ is the diagonal matrix with $(I_{G^+})_{ii}=1$ when $G_i(x)-y_i>0$ and zero otherwise.
In addition, the semismooth Newton method summarized in \cref{alg:ssnewton} can be used to obtain $x^+$, and as a result, $y^+$ and $z^+$.

Below we outline the semismooth Newton algorithm with line search to solve  the nonsmooth equation \eqref{eq:alm_subproblem_2}, which is the optimality condition of the $\nu$-th subproblem of the ALM. Throughout the statement of this algorithm, we maintain 
$\nu$ as the fixed outer iterate of the ALM, and write the semismooth Newton iteration as 
$t\geq 0$. 
\begin{algorithm}
	\caption{A semismooth Newton method for the $\nu$-th subproblem \eqref{eq:alm_subproblem_2}}\label{alg:ssnewton}
	\begin{algorithmic}[1]
		\Initialize{Set $t=0$, choose an initial point $x^{\nu,0}$, and choose backtracking line search parameters $\delta,c\in(0,1)$. Perform the following steps until a convergence criterion is reached.}
		\State Semismooth Newton direction: using direct factorization or an indirect method, obtain an approximate solution $d^{\,t}$ to the linear system
		\begin{align}\label{eq:ssnewton_linearsystem}
			V d + \nabla\varphi_{\sigma_\nu}(x^{\nu,t})=0,
		\end{align}
		where $V\in \partial \left({\nabla}\varphi_{\sigma_\nu}(x^{\nu,t})\right)$. 
		\State Line search (backtracking): Set $\rho_t=\rho$, where $\rho$ is the smallest nonnegative integer such that
		\begin{align*}
			\varphi_{\sigma_\nu}(x^{\nu,t} + \delta^\rho d^{\,t}) \leq \varphi_{\sigma_\nu}(x^{\nu,t}) + c\,  \delta^\rho \, \nabla\varphi_{\sigma_\nu}(x^{\nu,t})^\top d^{\,t}.
		\end{align*}

		\State Update solution:
		\begin{align*}
			x^{\nu,t+1} = x^{\nu,t} + \delta^{\rho_t} d^t.
		\end{align*}
	\end{algorithmic}
\end{algorithm}

In fact, it is not easy to obtain an element from $\partial \left({\nabla}\varphi_{\sigma_\nu}(x)\right)$ in the equation \eqref{eq:ssnewton_linearsystem} even if $G(x)$ is an affine function, see, for example, the discussion in \cite{dolgopolik2024note}. One remedy to overcome this difficulty is to construct an element $V$ in the following way:
\begin{equation}
\begin{aligned}    
\label{eq:varphi_hessian}
\left.
\begin{array}{l}
     V\coloneqq   \nabla^2 f(x) + \displaystyle\sum_{i=1}^{mL}\nabla^2 G_i(x) \cdot\Bigl(G(x) + \sigma_\nu^{-1}\lambda^\nu - \proj_{\overline{\mathcal{B}}}{G(x) + \sigma_\nu^{-1}\lambda^\nu}\Bigr)_{\!i}\\
        \phantom{\hat{\nabla}^2\varphi_{\sigma_\nu}(x) = 
        \nabla^2 f(x)} + \sigma_\nu \, J G(x)^\top \cdot (I - J_{\overline{\mathcal{B}}})\cdot JG(x)  + \sigma_\nu(I-J_X)+{\sigma_\nu}^{-1}M
\end{array}
\right\},\\*[0.15in]
\mbox{with}\quad J_{\overline{\mathcal{B}}} \in \partial \proj_{\overline{\mathcal{B}}}{G(x)-\sigma_\nu^{-1}\lambda},
\qquad
 J_X\in \partial \proj_{X}{x-\sigma_\nu^{-1}\mu}.
\end{aligned}
\end{equation}
An explicit form for $J_{\overline{\mathcal{B}}}$ is given in \cref{prop:explicit_gen_jac}, and the matrix $J_X$ is a diagonal matrix with the $i$-th diagonal element given by
\begin{align*}
	[J_X]_{ii} =
	\begin{cases}
		1,\quad&\text{if $x_i-\mu_i/\sigma_\nu\in[p_i,q_i]$};\\
		0,\quad&\text{otherwise.}
	\end{cases}
\end{align*}
Although such a matrix $V$ may not belong to $\partial \left({\nabla}\varphi_{\sigma}(x)\right)$, following the proof of
\cite{li2018highly}, it can be shown that $Vd \in \widehat{\partial}(\nabla \varphi_{\sigma_\nu}(x))(d)$ for some upper-semicontinuous and convex-valued set-valued map $\widehat{\partial}(\nabla \varphi_{\sigma_\nu}(\,\cdot\,))$, which is a superset of $\partial (\nabla \varphi_{\sigma_\nu}(x))$ at any $x$. Furthermore, the gradient function $\nabla \varphi_{\sigma_\nu}$ is semismooth with respect to  $\widehat{\partial}(\nabla \varphi_{\sigma_\nu}(\,\cdot\,))$. Therefore, one can safely use the matrix $V$ in \eqref{eq:varphi_hessian} in the semismooth Newton equation \eqref{eq:ssnewton_linearsystem}. The details of how to efficiently solve this linear equation are deferred to the next section.

Before ending this section, we  highlight a few extensions of problem \eqref{eq:source_problem} that can be similarly handled by our proposed semismooth-Newton based ALM.
Consider a problem involving only a single superquantile objective.
Let $\omega^j\coloneqq\{(A_{j,:},b_j)\}_{j=1}^m$, where in the context of the supervised learning problem,  we interpret $b_j\in\R$ as the response and $A_{j,:}^\top\in\R^{n}$ as the covariates of the $j$-th observation. The risk-averse regression problem takes the form of 
\[
\displaystyle\operatornamewithlimits{minimize}_{x\in \mathbb{R}^n} \quad \mathopfont{CVaR}_\tau\left(\psi(A_{j,:}^\top x- b_j)\bigr\}_{j=1}^m\right), 
\]
where  $\psi:\R\to\R$ measures the individual loss for each covariate-response pair, such as $\psi(t) = t^2$ for the regression problem. This problem is known as the performance at the top or hard example mining for (imbalanced) classification in machine learning \citep{lyu2020average, Shrivastava2016, Bengio2009,yang2022algorithmic}.
Our proposed framework  can directly handle convex continuously differentiable loss functions.
In some cases, it is preferable to consider nonsmooth loss functions such as the absolute loss $\psi(t) = \abs{t}$ or the hinge loss $\psi(t) = \max(1-t,0)$.  It is then convenient to consider a redefined risk set associated with a composite (empirical) superquantile operator $\widetilde{B}_{k(\tau), r} \coloneqq  \{y\in \mathbb{R}^m: \sum_{i=1}^{k(\tau)} \psi(y_i) \leq r\}$ for some $r\geq0$,
where $k(\tau) = (1-\tau)m$.
The reason that we consider the shifted set with $r\geq 0$ is because the newly defined operator $\mathopfont{CVaR}_\tau$ may not be translation equivariant. 
The projections onto these modified sets and their Jacobians need to be computed accordingly, but they retain similar sparsity structures to those detailed below in \cref{sec:sparsity}.
In particular, for the absolute value loss, the composite superquantile operator corresponds to the vector-$k$-norm, where an explicit form for the derivative of the associated projection operator is given in \cite{wu2014moreau}. We omit the detail here for brevity. 

\section{Implicit scenario reduction via the semismooth Newton equation}
\label{sec:scenario_reduction}

In this section, we explore how the inherent nonsmoothness of the superquantile operator enables implicit scenario reduction within the augmented Lagrangian subproblems, resulting in significant computational efficiencies in solving the linear equations of the semismooth Newton method.
For simplicity, we focus our discussion throughout this section to  a single superquantile constraint (\ie $L=1$) and omit the superscript $\ell$. Each $g(x;\omega^{j})\coloneqq A_{j,:}x+b_j$ is an affine function with $\omega^j\coloneqq\{(A_{j,:}\, , \, b_j)\}$ so that the constraint $\superquantile_{\tau}{G(x; \omega^{[m]})} \leq 0$ can be written as $\mathopfont{T}_k(Ax+b) \leq 0$ for
$k = (1-\tau)m$,
or equivalently, $Ax+b \in {\cal B}_k$.
Extending the discussion to multiple superquantile constraints with $L>1$ is straightforward, using the appropriate concatenations of $A^\ell$, $b^\ell$, and the Cartesian product  of the sets $\{\mathcal{B}_{k_{\ell}}\}$.

\subsection{Structured sparsity of the linear system \texorpdfstring{\eqref{eq:ssnewton_linearsystem}}{}}
\label{sec:sparsity}

Since the function $G$ is assumed to be affine in $x$, 
we consider the objective function of the ALM subproblem to be $\varphi_\sigma$ as in \eqref{eq:varphi} to simplify the discussion.
The matrix $V$ in \eqref{eq:varphi_hessian} then takes the form of
\begin{equation}\label{jacobian}
 V=   \nabla^2 f(x) + \sigma_\nu \, A^\top  (I - J_{{\mathcal{B}_k}})A  + \sigma_\nu(I-J_X)+{\sigma_\nu}^{-1}M.
        \end{equation}
After computing the gradient $\nabla\varphi_\sigma$ at an iterate $\bar{x}$ (which requires the sorting of $\{A_{j,:}^\top \bar{x} + b_j\}_{j=1}^m$), one can efficiently compute the generalized Jacobian of $\nabla\varphi_\sigma$ based on
 the already sorted values of $\{A_{j,:}^\top \bar{x} + b_j\}_{j=1}^m$, as specified in \cref{prop:explicit_gen_jac}. Below we will show that the extra cost for this computation of the generalized Jacobian is actually very low even if $m$ is large.

Let $\bar\alpha,\bar\beta,\bar\gamma$ be the associated index-sets of $\bar{y}\coloneqq\proj_{\mathcal{B}_k}{A\bar{x}+b-\sigma\bar{\lambda}}$ defined in \eqref{eq:order_structure} 
for $k = (1-\tau)m$.
When $\bar{x}$ and $\bar{\lambda}$ are such that the projection onto ${\cal B}_k$ is differentiable at $A\bar{x}+b-\sigma\bar{\lambda}$, we use the Jacobian formula given in \cref{prop:explicit_gen_jac}.
Otherwise, a sequence of differentiable points converging to $A\bar{x}+b-\sigma\bar{\lambda}$ can be constructed, and we use the limit of the gradients at this sequence based on whether $\bar k_1=k$ (Case 3) or not (Case 2) from \cref{prop:explicit_gen_jac}.

In Case 1 of \cref{prop:explicit_gen_jac} (when $\maxsum_k{A\bar{x}+b-\sigma^{-1}\bar{\lambda}}<0$), then the matrix in linear system \eqref{eq:ssnewton_linearsystem} has no contribution from the metric projection and can be constructed easily.

In Case 2, we have
$I-J_{\mathcal{B}_k}=\begin{bmatrix} Q & \\ & 0_{\bar\gamma} \end{bmatrix}$ where $Q = \begin{bmatrix} Q_{11} & Q_{12} \\ Q_{21} & Q_{22} \end{bmatrix}\in\R^{\bar k_1\times\bar k_1}$
is a symmetric projection matrix that is computable in closed form as a function of the index-sets.
In particular, for $\bar\rho \coloneqq k^2 - 2 k \bar k_0 + \bar k_0 \bar k_1$, direct computation shows that
    \begin{alignat*}{2}
        Q_{11} &= {\bar\rho}^{-1}({\bar k_1-\bar k_0})\ind_{\bar\alpha}\ind_{\bar\alpha}^\top, &\qquad
        Q_{12} &= {\bar\rho}^{-1}({k-\bar k_0}) \ind_{\bar\alpha}\ind_{\bar\beta}^\top,
        \\*
        Q_{21} &= Q_{12}^\top, &\qquad
        Q_{22} &= I - {\bar\rho}^{-1}{\bar k_0} \ind_{\bar\beta}\ind_{\bar\beta}^\top.
    \end{alignat*}
To ease notation, let $\bar\kappa\coloneqq(\bar\alpha,\bar\beta)$ with $\abs{\bar\kappa}=\bar{k}_1$.
Since $Q$ is a symmetric projection matrix, we have
\begin{align}
    \label{eq:TtildeTtilde}
    \sigma A^\top(I-J_{\mathcal{B}_k})A &= \sigma\cdot(A_{\bar\kappa,:})^\top Q A_{\bar\kappa,:} = \sigma\cdot(A_{\bar\kappa,:})^\top Q^\top Q A_{\bar\kappa,:} = \sigma T^\top T,\quad T \coloneqq QA_{\bar\kappa,:}.
\end{align}
Observe that $T\in\R^{\abs{\kappa}\times n}$ can be computed as follows:
    \begin{alignat*}{3}
        T_{i,:} &= {\bar\rho}^{-1}({\bar k_1-\bar k_0}) \ind_{\pi(\bar\alpha)}^\top A_{\pi(\bar\alpha),:} + {\bar\rho}^{-1}({k-\bar k_0}) \ind_{\pi(\bar\beta)}^\top A_{\pi(\bar\beta),:}\quad &i&\in\bar\alpha,\\
        T_{i,:} &= A_{\pi(i),:} + {\bar\rho}^{-1}({k-\bar k_0}) \ind_{\pi(\bar\alpha)}^\top A_{\pi(\bar\alpha),:} - {\bar\rho}^{-1}{\bar k_0} \ind_{\pi(\bar\beta)}^\top A_{\pi(\bar\beta),:}\quad &i&\in\bar\beta,
    \end{alignat*}
where $\pi(\mathcal{I})$ denotes the $\mathcal{I}$ components of $\pi$.
Since the rows of $T$ corresponding to $\bar\alpha$ are identical, we may construct $T^\top T = \widetilde{T}^\top\widetilde{T}$ for some $\widetilde{T}$ in a reduced space where we define, for any $i\in\bar{\alpha}$,
\begin{align}
    \widetilde T \coloneqq 
			\begin{bmatrix} \sqrt{\abs{\bar\alpha}} T_{i,:}\\T_{\bar\beta,:} \end{bmatrix} \in \R^{(\abs{\bar\beta}+1)\times n}.
\end{align}

\begin{figure}[h]
\begin{subfigure}{0.75\linewidth}
    \begin{tikzpicture}[baseline=(current bounding box.south),decoration=brace]
      \node (a) at (0,0) {};
      \matrix (m) [table,minimum height=1cm] {
        A_{\bar\alpha,:}^\top\\
        A_{\bar\beta,:}^\top\\
        A_{\bar\gamma_1,:}^\top\\
        \cdot\\
        \cdot\\
        \cdot\\
        A_{\bar\gamma_{\abs{\bar\gamma}}}\\
      };
      \node at (m-1-1.north east) {$\qquad\quad{\top}$};
    \end{tikzpicture}
    \begin{tikzpicture}[baseline=(current bounding box.south),decoration=brace]
    \node (a) at (0,0) {};
    \matrix (m) [table,minimum height=1cm, minimum width=1cm] {
    {}  \& {}  \& {}  \& {}   \&  {} \& {}     \& {}\\%1
    {}  \& {}  \& {}  \& {}   \&  {} \& {}     \& {}\\%2
    {}  \& {}  \& {} \& {}   \&  {} \& {}     \& {}\\%3
    {}  \& {}  \& {}  \& {{}}   \&  {} \& {}     \& {}\\%4
    {}  \& {}  \& {}  \& {}   \&  {{}} \& {}     \& {}\\%5
    {}  \& {}  \& {}  \& {}   \&  {} \& {{}}     \& {}\\%6
    {}  \& {}  \& {}  \& {}   \&  {} \& {}     \& {}\\%7
    };
    \node [anchor=north west] (big0) at ($(m-3-3.south east)-(-5mm,8mm)$) {\resizebox{0.1\linewidth}{!}{ 0}};
    
    \draw[red, opacity=0.5, ultra thick, fill] (m-1-1.north west) rectangle ($(m-1-1.south east)-(-0.01,0.01)$); %
    \draw[red, ultra thick] (m-1-1.north west) rectangle ($(m-1-1.south east)-(-0.01,0.01)$); %
    \draw[red!50!orange, opacity=0.2, fill] ($(m-1-2.north west)+(0.5mm,0mm)$) rectangle (m-1-2.south east); %
    \draw[red!50!orange, opacity=0.2, fill] ($(m-2-1.north west)+(0mm,-0.5mm)$) rectangle (m-2-1.south east); %
    \draw[orange, opacity=0.5, fill, ultra thick, dashed] ($(m-2-2.north west)+(0.4mm,-0.4mm)$) rectangle (m-2-2.south east); %
    \draw[orange, ultra thick, dashed] ($(m-2-2.north west)+(0.4mm,-0.4mm)$) rectangle (m-2-2.south east); %
    \draw[arrows = {Stealth[inset=0pt, length=0.5mm, angle'=90]-Stealth[inset=0pt, length=0.5mm, angle'=90]}, line width=3pt, orange] ($(m-2-2.north west)+(0.6mm,-0.6mm)$) -- (m-2-2.south east); %
    \draw[green!50!black, line width=1.6pt, dotted, opacity=0.25,fill] ($(m-3-3.north west)+(0.5mm,-0.5mm)$) rectangle (m-7-7.south east); %
    \draw[green!50!black, line width=1.6pt, dotted] ($(m-3-3.north west)+(0.5mm,-0.5mm)$) rectangle (m-7-7.south east); %
    \draw[decorate,decoration={calligraphic brace,amplitude=10pt},transform canvas={yshift=+0.2em},ultra thick] (m-1-1.north west) -- node[above=8pt] {\Large$\abs{\bar\alpha}+\abs{\bar\beta}$} (m-1-2.north east); %
    \draw[decorate,decoration={calligraphic brace,amplitude=10pt},transform canvas={yshift=+0.2em},ultra thick] (m-1-2.north east) -- node[above=8pt] {\Large$\abs{\bar\gamma}$} (m-1-7.north east); %
    \end{tikzpicture}
    \begin{tikzpicture}[baseline=(current bounding box.south),decoration=brace]
      \node (a) at (0,0) {};
      \matrix (m) [table,minimum height=1cm] {
        A_{\bar\alpha,:}^\top\\
        A_{\bar\beta,:}^\top\\
        A_{\bar\gamma_1,:}^\top\\
        \cdot\\
        \cdot\\
        \cdot\\
        A_{\bar\gamma_{\abs{\bar\gamma}}}\\
      };
    \end{tikzpicture}
    \caption{}
    \label{tks:second_order}
\end{subfigure}
    \;\;
\begin{subfigure}{0.2\linewidth}
  \begin{tikzpicture}[baseline=(current bounding box.south),x=0.5cm,y=2.5cm]%
  \node (a) at (0,0) {};
  \def\xzerosort{{%
    1.3-0.39,%
    1.25-0.39,%
    1.2-0.39,%
    0.45,%
    0.4,%
    0.3,%
    0.05,%
    -0.1,%
    -0.2,%
    -0.4,%
    -0.5,%
    -0.6,%
    -0.7,%
    -0.8,%
    -0.9,%
    -1.0,%
    -1.1,%
    -1.2,%
    -1.5,%
    -1.7,%
  }}
  \def\lambda{0.15}
  \def\xbarsort{{%
    1.3-0.39-\lambda,%
    1.25-0.39-\lambda,%
    1.2-0.39-\lambda,%
    0.25,%
    0.25,%
    0.25,%
    0.05,%
    -0.1,%
    -0.2,%
    -0.4,%
    -0.5,%
    -0.6,%
    -0.7,%
    -0.8,%
    -0.9,%
    -1.0,%
    -1.1,%
    -1.2,%
    -1.5,%
    -1.7,%
  }}
  
  \foreach \i in {0,1,...,19}{
    \node[draw, circle, fill=black, inner sep=0pt, minimum size=2mm, opacity=0.6] (a) at (+1,\xzerosort[\i]) {};
    \node[draw, circle, fill=blue, inner sep=0pt, minimum size=2mm, opacity=0.6] (b) at (-1,\xbarsort[\i]) {};
    \draw[->, thick, dotted] (a) -- (b);
  }
  \node[label={[label distance=0.2cm]0:$\substack{1}$}] (a) at (+1,\xzerosort[0]) {};
  \node[label={[label distance=0.2cm]0:$\substack{\bar{k}_0}$}] (a) at (+1,\xzerosort[2]-0.05) {};
  \node[label={[label distance=0.2cm]0:$\substack{\bar{k}_0+1}$}] (a) at (+1,\xzerosort[3]) {};
  \node[label={[label distance=0.2cm]0:$\substack{\bar{k}_1}$}] (a) at (+1,\xzerosort[5]-0.02) {};
  \node[label={[label distance=0.2cm]0:$\substack{\bar{k}_1+1}$}] (a) at (+1,\xzerosort[6]) {};
  \node[label={[label distance=0.2cm]0:$\substack{m}$}] (a) at (+1,\xzerosort[19]) {};

  \node[label={[label distance=0.2cm]180:$\substack{1}$}] (a) at (-1,\xbarsort[0]+0.05) {};
  \node[label={[label distance=0.2cm]180:$\substack{\bar{k}_0}$}] (a) at (-1,\xbarsort[2]) {};
  \node[label={[label distance=0.2cm]180:$\substack{k}$}] (a) at (-1,\xbarsort[3]) {};
  \node[label={[label distance=0.2cm]180:$\substack{\bar{k}_1+1}$}] (a) at (-1,\xbarsort[6]) {};
  \node[label={[label distance=0.2cm]180:$\substack{m}$}] (a) at (-1,\xbarsort[19]) {};
  
  \draw[-, ultra thick, red] (1.48,0.97) rectangle (-1.5,0.59);
  \draw[-, ultra thick, orange, dashed] (1.5,0.56) rectangle (-1.5,0.19);
  \draw[-, ultra thick, green!50!black, dotted] (1.5,0.16) rectangle (-1.5,-1.8);
  \node (a) at (-1.5,-1.8) {};
  \matrix [above left] at (-0.5,1.0) {
      \node [shape=rectangle, draw=red, line width=1, label={[label distance=1.6mm]above:$\bar\alpha$}] {};\\
  };
  \matrix [above left] at (0.75,1.0) {
      \node [shape=rectangle, draw=orange, line width=1, label=above:$\bar\beta$, dashed] (rectbeta) {};\\
  };

  \matrix [above left] at (2.0,1.0) {
    \node [shape=rectangle, draw=green!50!black, line width=1, label={[label distance=1mm]above:$\bar\gamma$}, draw=none] (rectgamma) {};\\
  };
  \coordinate (topmid) at ($(rectgamma.north west)!0.5!(rectgamma.north east)$);
  \coordinate (bottommid) at ($(rectgamma.south west)!0.5!(rectgamma.south east)+(0mm,0.5mm)$);
  \coordinate (leftmid) at ($(rectgamma.north west)!0.5!(rectgamma.south west)+(0mm,0.25mm)$);
  \coordinate (rightmid) at ($(rectgamma.north east)!0.5!(rectgamma.south east)+(0mm,0.25mm)$);
  \coordinate (bottomright) at ($(rectgamma.south east)+(0mm,0.5mm)$);
  \coordinate (bottomleft) at ($(rectgamma.south west)+(0mm,0.5mm)$);
  \foreach \corner in {rectgamma.north west, topmid, rectgamma.north east, bottomright, bottommid, bottomleft, leftmid, rightmid}{
    \fill [green!50!black] ($(\corner) + (-0.5mm,-0.5mm)$) rectangle ++(0.5mm,0.5mm); %
  }

  \matrix [above left] at (3.25,1.025) {
      \node [draw, circle, fill=black, inner sep=0pt, minimum size=2mm, opacity=0.6, label={[label distance=1mm]above:$\;\;\sorth{y}^0$}] {};\\
  };
  \matrix [above left] at (-1.75,1.025) {
      \node [draw, circle, fill=blue, inner sep=0pt, minimum size=2mm, opacity=0.6, label={[label distance=1mm]above:$\bar{y}$}] {};\\
  };
  \end{tikzpicture}
  \caption{}
    \label{tks:projection}
\end{subfigure}
  \caption{\textbf{(a)} Schematic of the structured sparsity (with affine $G(x) \coloneqq Ax + b$ and a single superquantile constraint) for constructing the generalized Jacobian at a sorted point $\sorth{y}^0\in\R^m$ based on \textbf{(b)} sorted projection $\bar{y}\in\mathcal{B}_k$ with index sets $\bar\alpha,\bar\beta,\bar\gamma$. Only the rows $A_{j,:}$ for  $j\in \bar\alpha\cup \bar\beta$ are relevant in the semismooth Newton equation. See \eqref{eq:varphi_hessian} for the case of the nonlinear function $G$.}
  \label{fig:tks}
\end{figure}
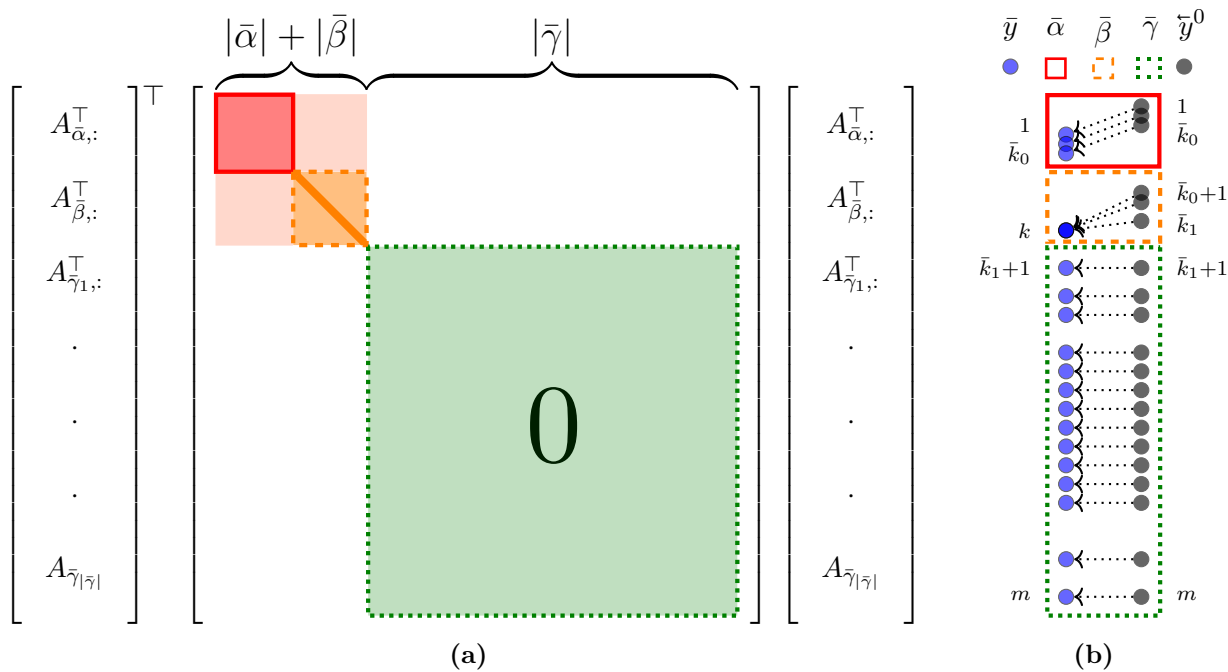

We further illustrate the special structured pattern of the generalized Jacobian for this case in \cref{fig:tks}. This depiction highlights how the Newton matrix incorporates only a subset of the scenarios, resulting from the form of the non-positive top-$k$-sum projection  in \cref{tks:projection}.
\cref{tks:second_order}  demonstrates that the projection matrix onto ${\cal B}_k$ is composed of four components:
(i) the $(\bar\alpha,\bar\alpha)$ block, which averages the scenarios within the index set $\bar{\alpha}$;
(ii) the $(\bar\alpha,\bar\beta)$ block and its transpose, both of which average the scenarios between $\bar\alpha$ and $\bar\beta$;
(iii) the $(\bar\beta,\bar\beta)$ block, which adjusts the $\bar\beta$ scenarios; and
(iv) the $(\bar\gamma,\bar\gamma)$ block, which excludes scenarios in $\bar\gamma$.

In Case 3 of \cref{prop:explicit_gen_jac}, we have $Q = \begin{bmatrix} (\abs{\bar\alpha}+\abs{\bar\beta})^{-1}\ind_{\bar\kappa}\ind_{\bar\kappa}^\top &\\& I_{|\bar\gamma|}\end{bmatrix}$. Hence, $A^\top (I-J_{\mathcal{B}_k}) A = A^\top Q A = \widetilde{T}\widetilde{T}^\top$ where $\widetilde{T} = \ind_{\bar\kappa}^\top A_{\bar\kappa,:} / \sqrt{{\abs{\bar\alpha}+\abs{\bar\beta}}}$,
\ie the scaled average of the $\bar\alpha$ and $\bar\beta$ rows of $A$.

Note that the Sherman-Morrison-Woodbury identity can be used to reduce the computational cost of solving the linear system \eqref{eq:ssnewton_linearsystem} as follows:
\begin{align}
    \label{eq:smw}
	V^{-1} &= (D + \widetilde{T}^\top(\sigma I)\widetilde{T})^{-1}
	= {D}^{-1} - {D}^{-1}\widetilde{T}^\top(\sigma^{-1} I + \tilde{T}{D}^{-1}\widetilde{T}^\top)^{-1}\widetilde{T}{D}^{-1},
\end{align}
where $D\coloneqq \nabla^2 f(\bar{x}) + \sigma^{-1} M + \sigma(I-J_X)$ is assumed to be invertible.

  {\small  
\begin{table}[t]
	\centering
	\begin{tabular}{ll}
	Terms &  Computational cost\\\toprule\\[-0.05in]
         Objective function value \eqref{eq:varphi}  & $O(m\log(m)+mn)$\\[0.1in]
         \midrule \\[-0.05in]

     Gradient \eqref{eq:varphi_gradient}  & $O(n\cdot\bar{k}_1)$ \\[0.1in]
     & ({\sl \small after computing the objective function \eqref{eq:varphi}})\\[0.1in]
     \midrule \\[-0.05in]

    Generalized Jacobian matrix \eqref{jacobian}  & $O\bigl(n\bar{k}_1+\min\bigl\{n\abs{\bar\beta}^2 + n^3,
            n\abs{\bar\beta}(n+\abs{\bar\beta}) + \abs{\bar\beta}^3\bigr\}\bigr)$\\[0.1in]

         and its inversion  &   ({\sl \small after computing the gradient \eqref{eq:varphi_gradient}}) \\[0.05in]
            
\bottomrule
	\end{tabular}
        \caption{\small 
        Computational cost for ALM subproblem \eqref{eq:palm_subproblem} when the objective is linear or diagonal convex quadratic. The regularization matrix $M$ is assumed to be diagonal. 
        The cost for the generalized Jacobian is according to the most expensive situation  in \cref{prop:explicit_gen_jac} (Case 2).
        }%
	\label{table:computational_cost}
\end{table}
}

\cref{table:computational_cost} summarizes the extra computational cost of obtaining first-vs-zeroth-order and second-vs-first-order information in our method.
After computing the value $\varphi_\sigma(\bar{x})$, the incremental cost of the gradient computation is  due to a matrix-vector product of the form $A^\top v$, where ${v} = \bar{y}-\proj_{\mathcal{B}_k}{\bar{y}}$ is a sparse vector with $\bar{k}_1$ nonzero elements, as can be seen from \cref{fig:tks}.
After evaluating $\nabla\varphi$, the extra cost to compute its generalized Jacobian is contributed by two steps: (i) the computation of the matrix for  the semismooth Newton equation, which involves summing across $A_{j,:}$ for $j\in\bar\alpha\cup\bar\beta$ and can be performed more efficiently than a sparse matrix-vector product of the same theoretical cost in the gradient computation; and (ii) solving the linear equation, which can be negligible in the large-$m$ case.
The first line in \cref{table:computational_cost} indicates that the function evaluation depends explicitly on the number of scenarios, $m$, via the sorting, linear transformation, and projection procedures, whereas the cost of computing the gradient, generalized Jacobian, and solving the linear equations is not directly dependent on $m$.

When $\abs{\bar\beta}\ll m$ at the current iterate $(\bar{x}, \bar{\lambda})$, \cref{table:computational_cost} shows that obtaining the generalized Jacobian of $\nabla \varphi_\sigma(\bar{x})$ and inverting the matrix incurs little extra cost after knowing its function value and the gradient, indicating the scalability of \cref{alg:palm} with respect to the number of scenarios.
In particular, due to the low rank structure of $T$ and additional simplifications arising when $f$ is linear or separable quadratic, the Newton system in \eqref{eq:ssnewton_linearsystem} may be computed in $n$-dimensions or $\abs{\bar\beta}+1$-dimensions due to the Sherman-Morrison-Woodbury identity, where a scalar parameter $\theta>0$ may be used to tune this trade-off.
Numerical experiments on synthetic data support the conclusion that in certain instances, the incremental cost of obtaining the Newton direction may be smaller than the incremental cost of evaluating the gradient; see \cref{tab:detail_synthetic} for detailed timing results.
Note that $J G^\top \, Q \, J G$ still may be computed with the reduced  dimension $\abs{\bar\beta}$ when $G$ is nonlinear.

\subsection{Implicit vs explicit scenario reduction}
\label{sec:implicit_scenario_reduction}

As discussed in the last section, our semismooth-Newton-based ALM method can be viewed as an implicit scenario reduction approach. 
We will now compare this approach to an explicit scenario reduction strategy used in the outer-approximation (OA) scheme  in \citet[Chapter 6.C]{royset2021primer}. This comparison highlights the potential advantages of our approach.

The OA method leverages the intuition that many $u_{\ell j}=0$ at the solution to the reformulation \eqref{eq:nlp} when $\tau_{\ell}$ is close to $1$, and in its simplest form, proceeds by maintaining active sets ${\cal A}_{\ell}^\nu$ of candidate effective scenarios with ${\cal A}_{\ell}^\nu = \{1,\ldots,m\}\setminus\{j : u_{\ell j}=0\}$ with the property that ${\cal A}_\ell^\nu \subseteq {\cal A}_\ell^{\nu+1}$ for iteration $\nu$.
The set for the next iteration is obtained by ${\cal A}_\ell^{\nu+1} = {\cal A}_\ell^{\nu} \cup J_\ell^\nu$ where $J_\ell^\nu$ contains the indices of a fraction of the
$(1-\tau_\ell)m$-largest
elements of the values $\{u^{\nu}_{\ell j} - g^\ell(x^\nu; \omega^{j})+ t_\ell^\nu\}_{j=1}^m$. 
The method terminates in iteration $\nu+1$ when ${\cal A}_\ell^{\nu} = {\cal A}_\ell^{\nu+1}$.
Achieving optimal performance with the OA method necessitates fine-tuning the update strategy for ${\cal A}_\ell^\nu$ and often involves solving subproblems $\nu$ to progressively tighter tolerances. In fact, the problem of adding an appropriate amount of additional scenarios is challenging even when the candidate solution is near the true solution set, as there is a trade-off between number of subproblems solved and the difficulty of solving each subproblem.
In addition, the set of effective scenarios  at the solution may be larger than
$(1-\tau_\ell)m$,
but this number is not known \emph{a priori}.

In contrast,  our ALM method implicitly maintains an approximation to the optimal set of effective scenarios through the partition set $\bar\alpha\cup\bar\beta$, and it uses an inherent update rule for ${\cal A}_\ell^\nu$ that is guided by the  generalized Jacobian of $\mathopfont{proj}_{\overline{\cal B}}$.
Consequently, the ALM method effectively prunes and approximates this set in each iteration, continuously refining the pool of scenarios considered.

\subsection{Comparison to the existing semismooth-Newton based ALM  that leverages the structured sparsity}

There are a few works in the existing literature using seemingly similar semismooth-Newton based ALM to solve problems with structured sparsity, such as   \citet{li2018fused,li2018highly, li2021fast, wu2022convex}.
All of these works are based on the observation that the solution of the problem is sparse. 
Among them, \cite{wu2022convex} is most closely related to our superquantile problem \eqref{eq:source_problem},
but a major distinction lies in how the algorithm leverage the sparsity. The optimization model in the latter paper explicitly added an $\ell_1$ regularizer to the objective function so that the solution is sparse. 
In the case when $G^{\,0}(x;\omega^{[m]})\coloneqq \{\abs{A_{i,:}^\top x+b_i}\}_{i=1}^m$ as in \citet{wu2022convex}, the ALM subproblems for the dual formulation requires the linear solution of a system in the form of $APA^\top + Q$, where $P$ is a sparse diagonal  matrix based on the generalized Jacobian of the proximal mapping of the $\ell_1$ regularizer, and $Q$ corresponds to the generalized Jacobian of the proximal operator for the superquantile.
However, when the sparse regularizer is not present, the system instead requires the solution of a linear system involving a matrix of the form $AA^\top + Q$, which can be significantly more costly.

In contrast, the corresponding linear solve in our formulation takes the form $A^\top Q A + \nabla^2 f(x)$.
When the objective $f$ is linear, we have $\nabla^2 f(x)=0$.
In this way, our framework can employ careful numerical linear algebra to reduce the computational cost of the linear solve without the use of a sparse regularizer in the objective.
The matrix $Q$ serves as a form of \emph{implicit} scenario reduction.
This is made explicit by the form of $Q$, which subsets certain rows  of the matrix $A$ and uses variational analysis to adjust the active $(1-\tau_0)$-tail scenarios rather than an ad hoc, manually-tuned procedure.

Since our technique scales in the number of scenarios, we expect our method to be complementary to the decision-space scalability developed in \cite{wu2022convex}.

\section{Experiments}
\label{sec:experiments}

To study the performance of our proposed framework, we conduct two simulation studies: the first uses synthetic data (\cref{sec:experiments:synthetic}), and the second  demonstrates its practical application through quantile regression on a realistic dataset containing more than 30 million records (\cref{sec:experiments:qr}).
The implementation in Julia is available at \url{https://github.com/jacob-roth/superquantile-opt}.

We compare the performance of our semismooth-Newton-based ALM implemented in Julia v1.9.3 \citep{julia} against several other approaches: (i) The interior point methods for solving the nonlinear reformulation \eqref{eq:nlp} by Gurobi v11.0 \citep{gurobi}; (ii) The alternating direction method of multipliers implemented by OSQP v0.6.3 \citep{osqp}; (iii) the Portfolio Safeguard v3.2  \citep{psg}; and (iv)
an implementation of the outer approximation algorithm based on  \citet[Chapter 6.C]{royset2021primer} that calls Gurobi as a subproblem solver.
For simplicity, we refer to these methods as (a) \underline{GRB}, (b) \underline{OSQP}, (c) \underline{PSG}, and (d) \underline{G-OA} respectively in the later discussion.

\subsection{Implementation details}
\label{sec:implementation_detail}
Before showing the numerical results, we first provide some implementation details of our approach.
At a candidate primal-dual tuple $(x,y,z,\lambda,\mu)$ for problem \eqref{eq:source_problem_smooth} and a  given tolerance $\epsilon$, we stop the algorithm if $\eta\coloneqq\max\{\eta_{\text{p}},\eta_{\text{d}},\eta_{\text{r}}\} \leq \epsilon$,
where
\begin{subequations}
\label{eq:stopping_condition}
\begin{align*}
    \eta_{\mathrm{p}} &\coloneqq \max\Bigl\{\norm[\big]{
        \max( Ax + b - 
        y,0)
    }\, /(1+\norm{b}),\;\; \norm[\big]{y - \proj_{\overline{\mathcal{B}}}{y}} \, / (1+\norm{y}),\\*[0.05in]
    &\phantom{{}=\max\Bigl\{} \norm{x - z} \, / (1 + \norm{z}),\;\;\norm[\big]{z - \proj_{X}{z}} \, / (1+\norm{z})
    \Bigr\},\\[0.05in]
   \eta_{\mathrm{d}} &\coloneqq \norm{ A^\top\lambda +  \mu + \nabla f(x)} \, /(1+\norm{\nabla f(x)}),\\[0.05in]
    \eta_{\mathrm{r}} &\coloneqq {|\text{obj}_{\text{p}} - \text{obj}_{\text{d}}|}/{(1 + \abs{\text{obj}_{\text{p}}})},\\
   & \qquad  \mbox{with} \; \text{obj}_{\text{p}} \coloneqq f(x) \;\,\mbox{and}\;\, 
    \text{obj}_{\text{d}} \coloneqq -f^*(-A^\top\lambda- \mu) -\delta^*_{\overline{\mathcal{B}}}(-\lambda) - \delta^*_X(-\mu) + b^\top\lambda.
\end{align*}
\end{subequations}
The above conditions ensure the KKT conditions  \eqref{eq:source_problem_smooth_kkt} hold when $\epsilon=0$.
For the ALM subproblems \eqref{eq:palm_subproblem}, we use the following termination criteria:
\begin{align}
    \label{eq:stop_a}
    \norm[\big]{\nabla \varphi_{\sigma_\nu}(x^{\nu+1})}_2 &\leq \frac{\epsilon_\nu }{\sigma_\nu},\quad \epsilon_\nu \geq 0,\quad \sum_{\nu=0}^\infty \epsilon_\nu<+\infty\\*
    \norm[\big]{\nabla \varphi_{\sigma_\nu}(x^{\nu+1})}_2 &\leq \frac{\delta_\nu}{\sigma_\nu}  \norm[\big]{(x^{\nu+1},y^{\nu+1},z^{\nu+1}) - (x^{\nu},y^{\nu},z^{\nu})}_{M}, \; 0\leq\delta_\nu<1,\quad\sum_{\nu=0}^\infty\delta_\nu<+\infty,
\end{align}
where $M$ is the positive semidefinite matrix chosen at iteration $\nu$ in \cref{alg:palm} and $\varphi_{\sigma_\nu}$ is defined  in \eqref{eq:varphi}.
These conditions are given in \cite{li2020lp,wu2022convex} and are based on the guidelines developed by \cite{rockafellar1976augmented} for the proximal point algorithm.

Note that if the subproblem is solved exactly with $M\equiv0$, then the resulting iterate is dual feasible, \ie satisfies $\nabla f(x) + A^\top \lambda+\mu = 0$.
This holds because the ALM updates give
\begin{align*}
    0 &= \nabla f(x^{\nu+1}) + \sigma_\nu A^\top \big(Ax^{\nu+1}+b+\sigma_\nu^{-1}\lambda^\nu - y^{\nu+1}\big) + \sigma_\nu(x^{\nu+1} + \sigma_\nu^{-1}\mu^\nu - z^{\nu+1}),\\
    \lambda^{\nu+1} &= \lambda^{\nu} - \sigma_\nu(y^{\nu+1} - Ax^{\nu+1}-b),\\
	\mu^{\nu+1} &= \mu^{\nu} - \sigma_\nu(z^{\nu+1} - x^{\nu+1}),
\end{align*}
which in turn gives
\begin{align*}
    \nabla f(x^{\nu+1}) &= -\sigma_\nu A^\top \big(Ax^{\nu+1}+b+\sigma_\nu^{-1}\lambda^\nu - y^{\nu+1}\big) - \sigma_\nu(x^{\nu+1}+\sigma_\nu^{-1}\mu^\nu-z^{\nu+1})= -A^\top\lambda^{\nu+1} - \mu^{\nu+1}.
\end{align*}

Next, an important aspect of the numerical implementation involves mitigating the computational burden associated with sorting operations, as a sorting step is executed in each line search iteration.
Although the optimal index $\bar{k}_1$ cannot be determined \emph{a priori}, an existing permutation may already be optimal, based on the observation that the projection onto the top-$k$-sum constraint set only requires the  elements  in the set $\bar\alpha\cup\bar\beta$ to be in sorted order.
At the $\nu$-th iteration of the ALM, {within the line search procedure}, we find that the candidate solution's current $\bar{k}_1^{(\nu)}$ augmented by a small buffer  provides an effective estimate for the next iterate's index $\tilde{k}_1^{(\nu+1)}$.
We leverage this observation to perform a partial quick-sort of the top $\tilde{k}_1^{(\nu+1)}$ elements and thereafter compute the projection based on this partial permutation.
A projection $\tilde{z}$ computed from the partial permutation of an initial point $z^0\in\R^n$ can be certified as correct by verifying that  $z^0_i > \tilde{z}_{\tilde{k}_1}$ for $i\in \{\tilde{k}_1,\tilde{k}_1+1,\ldots,n\}$.
Effectively, this checks whether or not the projection $\tilde{z}$ is sorted and obeys the KKT conditions $\tilde{z}_{\tilde{k}_1+1} < \tilde{z}_{\tilde{k}_1}$ and $\tilde{z}_{i}\geq \tilde{z}_{i+1}$ for $i\in \{\tilde{k}_1,\tilde{k}_1+1,\ldots,n\}$.
We repeatedly solve the projection problem associated with progressively larger partial-sorts until the  above conditions are satisfied.
In practice, we find that re-sorts are required infrequently (approximately 5\% of iterations necessitate this adjustment) and can be efficiently warm-started.

Through experimentation, we find that the OA method performs relatively well by targeting $\approx$10 total subproblems and solving subproblem $\nu$ to progressively tighter tolerances when using Gurobi's barrier solver.
Gurobi has no warm-start feature for the barrier method, but we find that performance did not improve using more frequent warm-starts with primal/dual simplex methods.
In tuning the $L=1$, $k=1\%m$ case, we set $J^\nu$ in iteration $\nu$ to be the largest $3k$ elements of $(Ax+b)-\ind t-u$ in an effort to perform an aggressive scan of the risk space (where $u_j=0$ for $j\notin \mathcal{A}^\nu$).
Note that the size of $\mathcal{A}^{\nu+1}$ may increase modestly relative to $3k$ if many active elements in iteration $\nu$ remain active in $\nu+1$, which is often the case in later iterations.

\subsection{Results on synthetic data}
\label{sec:experiments:synthetic}
In this section, we generate synthetic data for the superquantile-constrained problem  \eqref{eq:source_problem}.
Problems instances are randomly generated according to the following procedure, which ensures that a solution exists.
The box constraint is $X\coloneqq\{x\in\R^n : -\ind\leq x\leq \ind\}$.
For each $\ell\in\{1,\ldots,L\}$, constraint data $A^{\ell}_{ij}\sim \texttt{Normal}(0,10^2)$ with unit column infinity-norm (\ie $A^{\ell}_{:,j}$ is normalized by $\norm{A^{\ell}_{:,j}}_\infty$), and $\tilde{b}^\ell\sim \texttt{Normal}(0,1)$. To ensure that the problem feasible, we define $b^\ell$ to be a shift of $\tilde{b}^\ell$ via the following procedure: for $p\in\{1,\ldots,\ceil{\log(N)}\}$, generate feasible points $\tilde{x}^p\in X$ uniformly by $\tilde{x}^p \coloneqq -\ind + 2u$ for $u\sim \texttt{Uniform}(0,1)^n$.
Subsequently, for each $\ell\in\{1,\ldots,L\}$, define $\tilde{y}^{\ell,p}\coloneqq A^\ell \tilde{x}^p + \tilde{b}^\ell$ and set $b^\ell\coloneqq\tilde{b}^\ell - \min_p\bigl(\maxsum_{k_\ell}{\tilde{y}^{\ell,p}}\bigr)$ to generate difficult instances.
We consider two types of objective functions: (i) diagonal (stongly) convex quadratic function
$f(x)=\tfrac12 x^\top Cx + c^\top x$ where $C_{ii}\sim \abs{\texttt{Normal}(0,1)}$ and $c_i\sim \texttt{Normal}(0,1)$; (ii) linear function $f(x) = c^\top x$, where $c$ is generated in the same way as (i). 
The ALM method is initialized at the feasible point $x^0=0$, $\mu=0$, and $\lambda=0$, and we find that the method performs similarly for initial points with elements drawn uniformly within the bounds of $X$.

The experiments are conducted on a Dell workstation running a Windows environment, equipped with 8 physical processors and 32GB of RAM. We compare our ALM method with GRB, G-OA, and PSG for computing highly accurate solutions  (where $\epsilon=10^{-8}$), and with OSQP for computing solutions with lower accuracy (where $\epsilon=10^{-3}$).
While both ALM and Gurobi barrier methods are able to use multiple threads (via BLAS/LAPACK routines and parallelized matrix factorization, respectively), establishing a fair comparison is somewhat ambiguous since Gurobi's procedure is not well documented.
Since OSQP's default solver is single-threaded, we limit GRB, G-OA, and ALM to one execution thread %
in our experiments.
PSG's multithreading capabilities are not documented, and we use the default settings.
In addition, we turn off Gurobi's presolve option and OSQP's polishing procedure, as our findings indicate that these features do not enhance performance.

In scenarios where $m\gg n$, we find that Gurobi's barrier method is more efficient than either primal simplex or dual simplex; when $m\approx n$, the barrier method is comparable to the dual simplex for less stringent risk tolerance ($1-k(\tau)/m$) and may be slower for more stringent risk tolerance, but both are more efficient than primal simplex; when $m\ll n$, the dual simplex is most efficient and may significantly outperform the barrier method, as well as ALM.
Since the focus of the comparison is on instances where $m\gg n$, we report solve times using the minimum of Gurobi's barrier and dual-simplex methods in a grid of instances with $m\geq n$ when reporting GRB solve statistics.
We set the feasibility and optimality tolerance to $10^{-8}$.
To obtain OSQP solutions with $\epsilon=10^{-3}$ precision, we set the relative and absolute tolerances to $10^{-5}$ and $10^{-3}$, respectively.
In the PSG suite, we set the precision to $10^{-8}$.
For PSG, we use the \texttt{HELI} solver, which uses Gurobi as a subproblem solver, for both linear and quadratic problems.
The \texttt{HELI} solver managed to solve instances to the specified tolerance more reliably than the PSG suggested solvers,\footnote{In fact, the quadratic linear (quadratic) solvers suggested in the PSG \href{http://www.aorda.com/html/PSG_Help_HTML/index.html?optimization_solvers_in_matlab.htm}{documentation} failed to find solutions with the desired tolerance in 75\% (96\%) of the instances tested below.} \texttt{CAR} and \texttt{VAN}, as well as their Gurobi-based variants \texttt{CARGRB} and \texttt{VANGRB}.

Finally, problems are generated for a variety of parameters: $m\in\{2^{14}/L,2^{15}/L,\ldots,2^{20}/L\}$, $n\in\{2^7,2^8,\ldots,2^{13}\}$, $L\in\{1,10\}$, and $k=\{0.01,0.1\}\times m$.
However, some combinations of $m$, $n$, and $L$ exceed the 32GB RAM capacity,  so we exclude them in the experiments.
In addition, PSG's Gurobi-based solvers are constrained by  the number of decision variables they can handle.
Consequently,  we only select instances that are solvable by PSG.

\cref{fig:synthetic_heatmap}, which splits over the next couple of pages, shows the relative solve time of method $i$ versus ALM with $i\in\{$a = GRB, b = G-OA, c = PSG, d = OSQP$\}$ for a range of problems across parameters $m$, $n$, $L$, and $k$ for both linear and diagonal quadratic objectives.

The first row of \cref{fig:heatmap_linear:grb} shows that across linear and quadratic objectives, ALM was never slower than GRB for instances with high risk aversion ($k=1\%m$).
We highlight that ALM's relative performance generally improves both as $m$ increases within each subplot and as $L$ increases (between consecutive subplots).
In the largest-$m$ instances, ALM demonstrates improvement over GRB of at least 20x and up to 70x when both trends are maximized.
Similar trends are observed in the second row of \cref{fig:heatmap_linear:grb}, corresponding to lower risk aversion instances (larger $k$).
In these instances, ALM's relative performance degrades due to the increased cost of constructing $\widetilde{T}$ associated with larger $k$ but maintains an advantage of at least 5x and up to 30x on large-$m$ problems.
For relatively small-$m$ and large-$n$ problems, GRB's dual simplex can be very efficient, beating ALM by up to 4x on problems with a linear objective, but losing its advantage for problems with quadratic objective.

In \cref{fig:heatmap_linear:grb_oa}, G-OA obtained solutions satisfying primal and dual feasibility in each instance, but objective gap tolerances were violated by a factor of $\approx 10$, contributing to the vertical shading.
The top row of \cref{fig:heatmap_linear:grb_oa} shows that ALM retains a solve-time advantage of about 5-10x over G-OA in the large-$m$ cases for the high risk-aversion instances.
The bottom row shows that ALM retains a more significant advantage in lower risk-aversion instances.
In contrast to GRB, this difference is due to the fact that G-OA must scan through a larger set of effective scenarios, degrading its performance as the size of the OA subproblem grows.
This issue could perhaps be overcome by using adaptive scenario pruning strategies, which indicates the importance of adaptive scenario selection in ALM for problems with a larger effective scenario tail.

\captionsetup[subfigure]{skip=-5pt plus 2pt minus 2pt}
\setlength{\abovecaptionskip}{5pt plus 2pt minus 2pt}
\edef\myTrim{\trimlr{} \trimbottom{} \trimlr{} \trimtop}

\begin{figure}[ht]
    \centering
    \setcounter{subfigure}{0}
    \setcounter{subsubfigure}{0}
    \begin{subfigure}[b]{\textwidth}
        \centering
        \setcounter{subfigure}{0}
        \setcounter{subsubfigure}{0}
        
        \begin{minipage}[b]{0.49\textwidth}
            \centering
            \caption*{GRB: linear objective} %
            \begin{subsubfigure}[b]{0.49\linewidth}
                {\resizebox{\linewidth}{!}{\includegraphics[trim=1mm 5mm 1mm 15mm,clip]{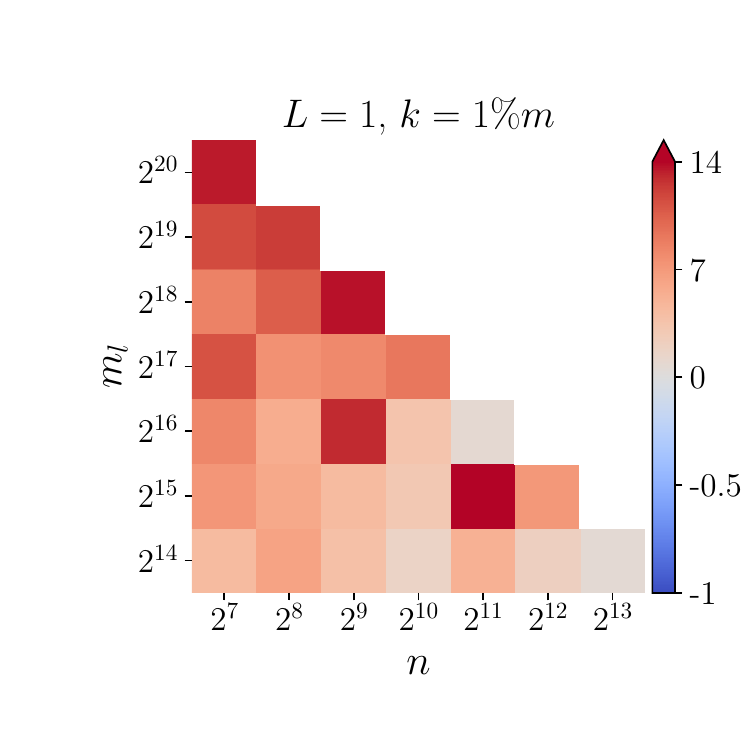}}}
                \caption{\!\!\!\!\!\!\!\!\!\!}
                \label{fig:heatmap_linear:grb:a}
            \end{subsubfigure}
            \begin{subsubfigure}[b]{0.49\linewidth}
                {\resizebox{\linewidth}{!}{\includegraphics[trim=1mm 5mm 1mm 15mm,clip]{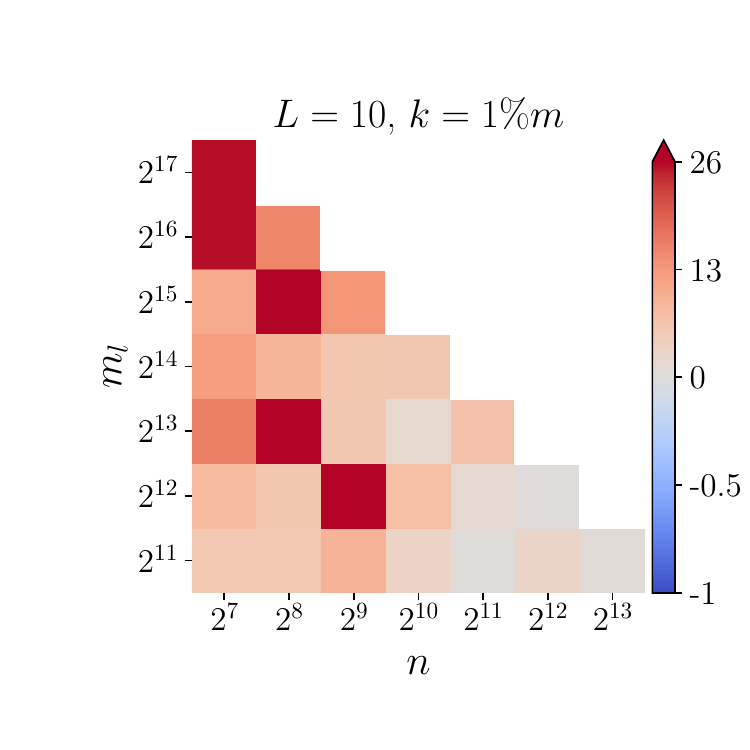}}}
                \caption{\!\!\!\!\!\!\!\!\!\!}
                \label{fig:heatmap_linear:grb:b}
            \end{subsubfigure}
            \begin{subsubfigure}[b]{0.49\linewidth}
                {\resizebox{\linewidth}{!}{\includegraphics[trim=1mm 5mm 1mm 15mm,clip]{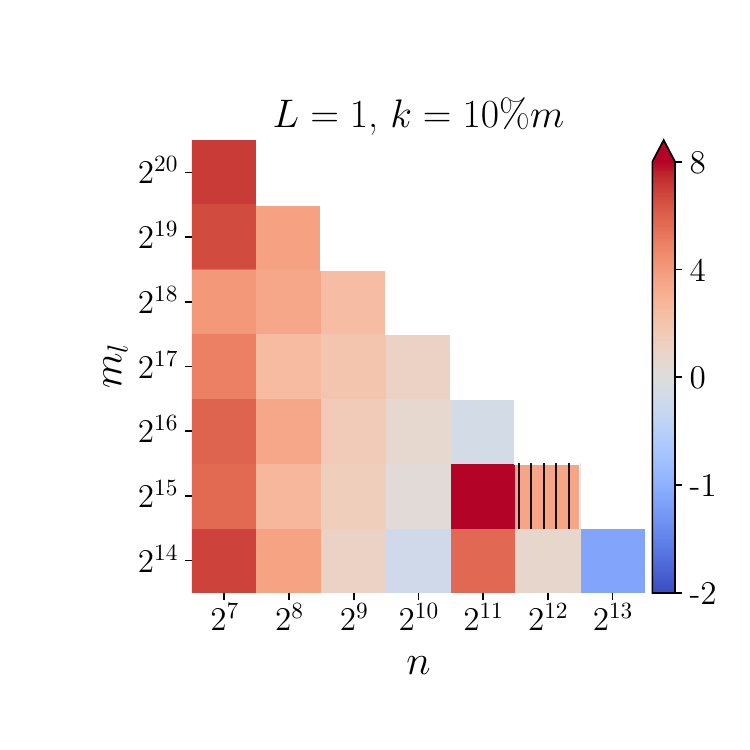}}}
                \caption{\!\!\!\!\!\!\!\!\!\!}
                \label{fig:heatmap_linear:grb:c}
            \end{subsubfigure}
            \begin{subsubfigure}[b]{0.49\linewidth}
                {\resizebox{\linewidth}{!}{\includegraphics[trim=1mm 5mm 1mm 15mm,clip]{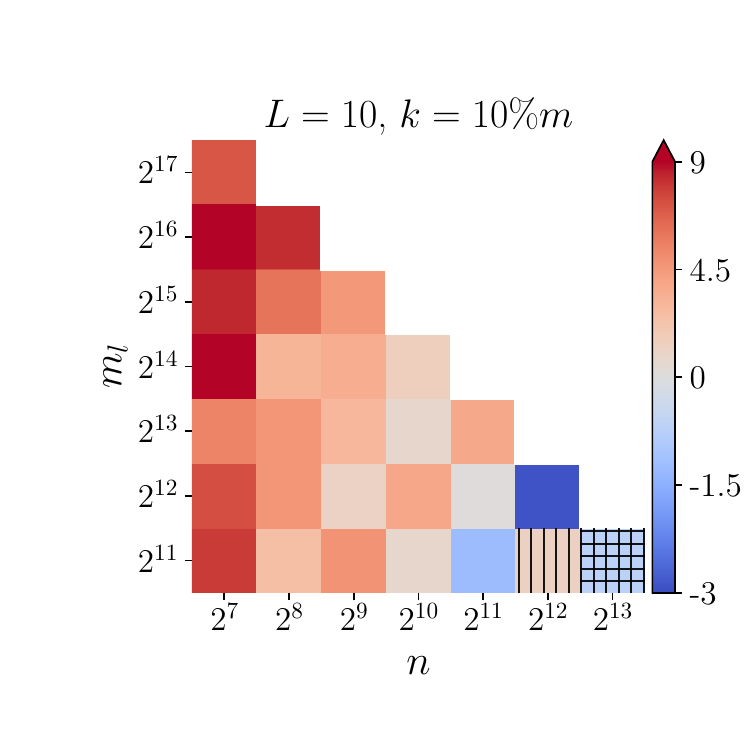}}}
                \caption{\!\!\!\!\!\!\!\!\!\!}
                \label{fig:heatmap_linear:grb:d}
            \end{subsubfigure}
            \captionsetup{labelformat=empty}
            \caption{}
            \label{fig:heatmap_linear:grb}
        \end{minipage}
        \hfill
        \setcounter{subfigure}{0}
        \begin{minipage}[b]{0.49\textwidth}
            \centering
            \caption*{GRB: quadratic objective} %
            \begin{subsubfigure}[b]{0.49\linewidth}
                {\resizebox{\linewidth}{!}{\includegraphics[trim=1mm 5mm 1mm 15mm,clip]{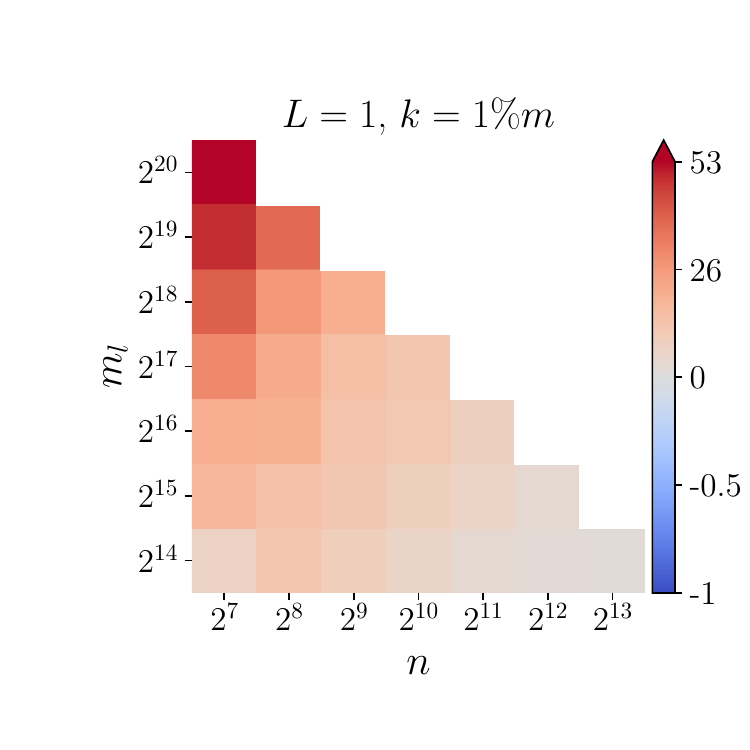}}}
                \caption{\!\!\!\!\!\!\!\!\!\!}
                \label{fig:heatmap_quadratic:grb:a}
            \end{subsubfigure}
            \begin{subsubfigure}[b]{0.49\linewidth}
                {\resizebox{\linewidth}{!}{\includegraphics[trim=1mm 5mm 1mm 15mm,clip]{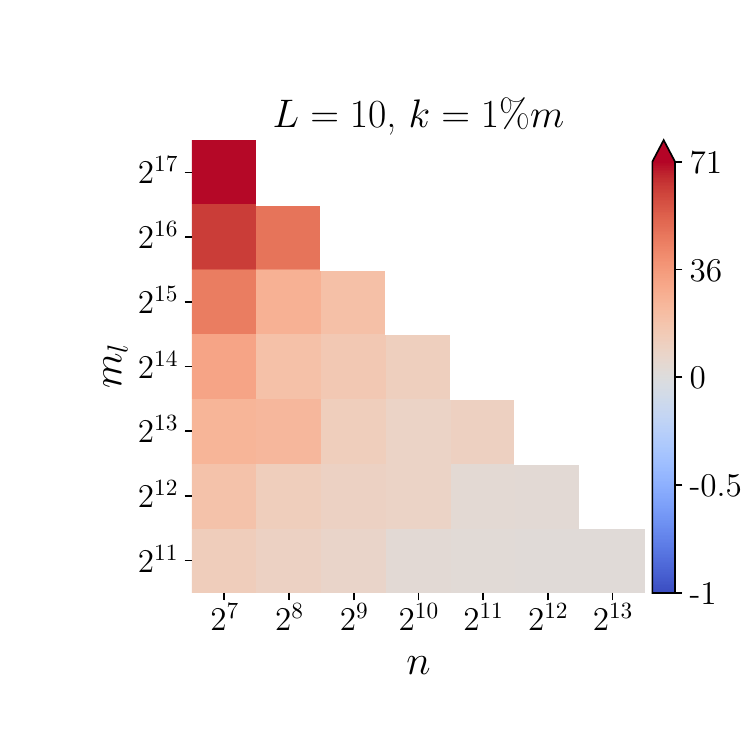}}}
                \caption{\!\!\!\!\!\!\!\!\!\!}
                \label{fig:heatmap_quadratic:grb:b}
            \end{subsubfigure}
            \begin{subsubfigure}[b]{0.49\linewidth}
                {\resizebox{\linewidth}{!}{\includegraphics[trim=1mm 5mm 1mm 15mm,clip]{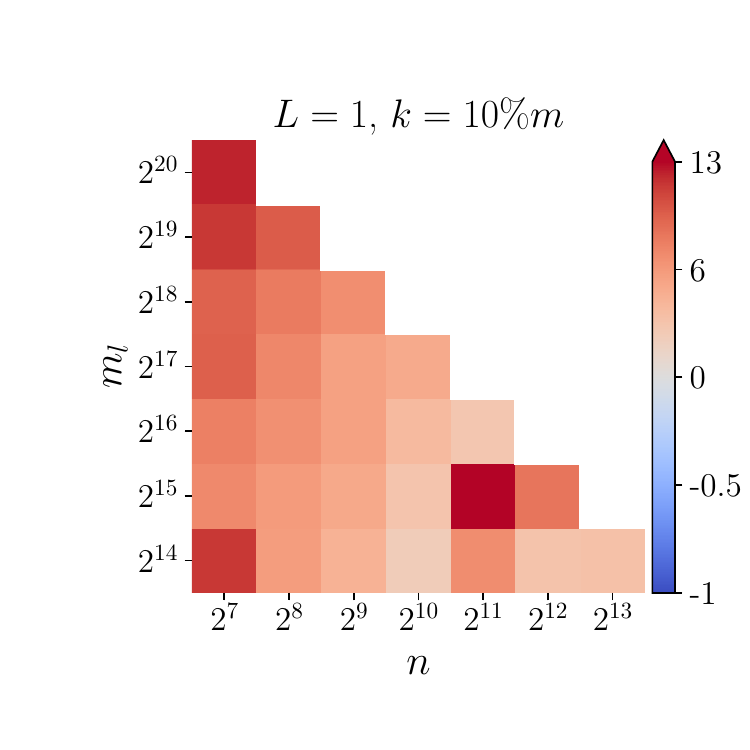}}}
                \caption{\!\!\!\!\!\!\!\!\!\!}
                \label{fig:heatmap_quadratic:grb:c}
            \end{subsubfigure}
            \begin{subsubfigure}[b]{0.49\linewidth}
                {\resizebox{\linewidth}{!}{\includegraphics[trim=1mm 5mm 1mm 15mm,clip]{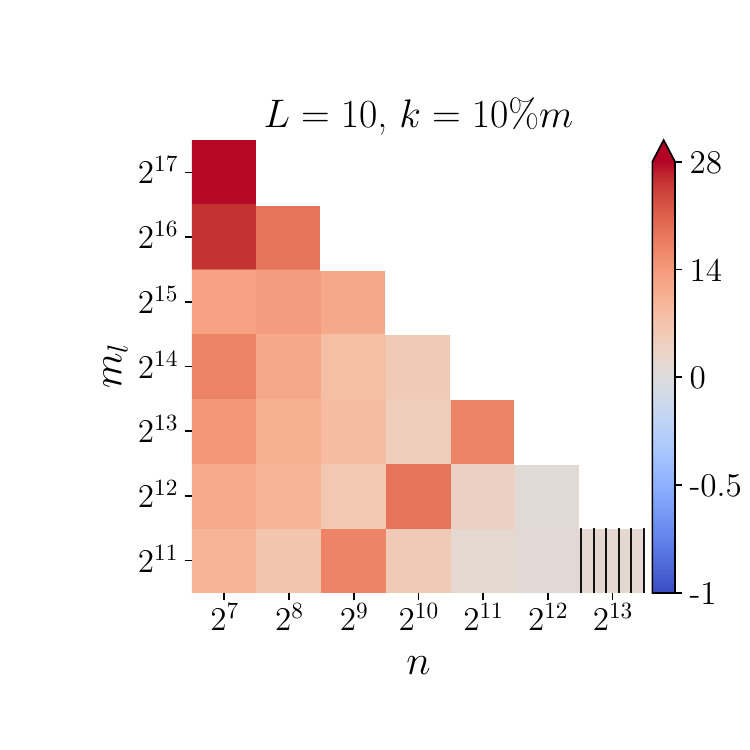}}}
                \caption{\!\!\!\!\!\!\!\!\!\!}
                \label{fig:heatmap_quadratic:grb:d}
            \end{subsubfigure}
            \captionsetup{labelformat=empty}
            \caption{}
            \label{fig:heatmap_quadratic:grb}
        \end{minipage}
        \captionsetup{labelformat=empty}
        \caption{}
        \label{fig:heatmap:grb}
    \end{subfigure}
\end{figure}
\begin{figure}[ht]
    \ContinuedFloat
    \centering

    \vspace{-6mm}
    \begin{subfigure}[b]{\textwidth}
        \centering
        \setcounter{subfigure}{1}
        \setcounter{subsubfigure}{0}
        
        \begin{minipage}[b]{0.49\textwidth}
            \centering
            \caption*{G-OA: linear objective} %
            \begin{subsubfigure}[b]{0.49\linewidth}
                {\resizebox{\linewidth}{!}{\includegraphics[trim=1mm 5mm 1mm 15mm,clip]{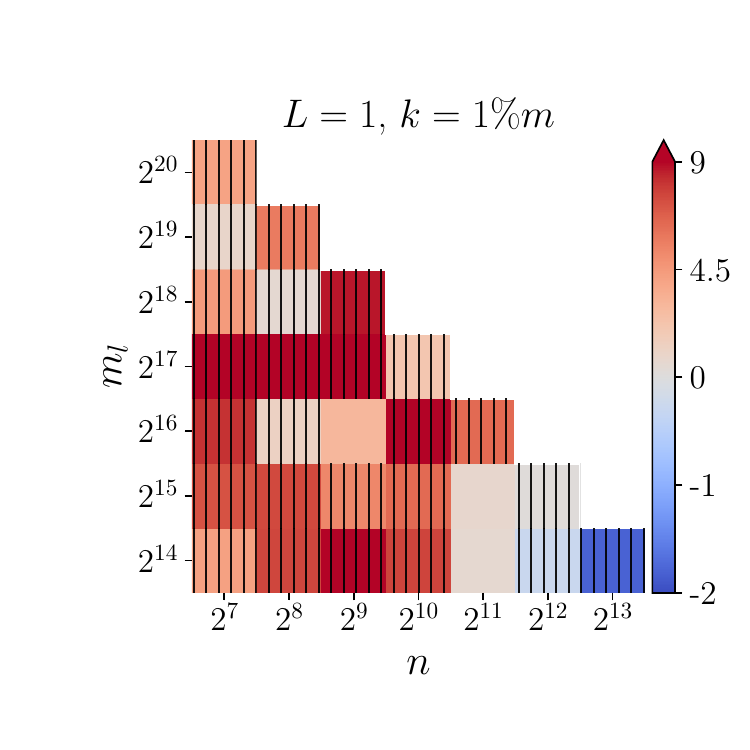}}}
                \caption{\!\!\!\!\!\!\!\!\!\!}
                \label{fig:heatmap_linear:grb_oa:a}
            \end{subsubfigure}
            \begin{subsubfigure}[b]{0.49\linewidth}
                {\resizebox{\linewidth}{!}{\includegraphics[trim=1mm 5mm 1mm 15mm,clip]{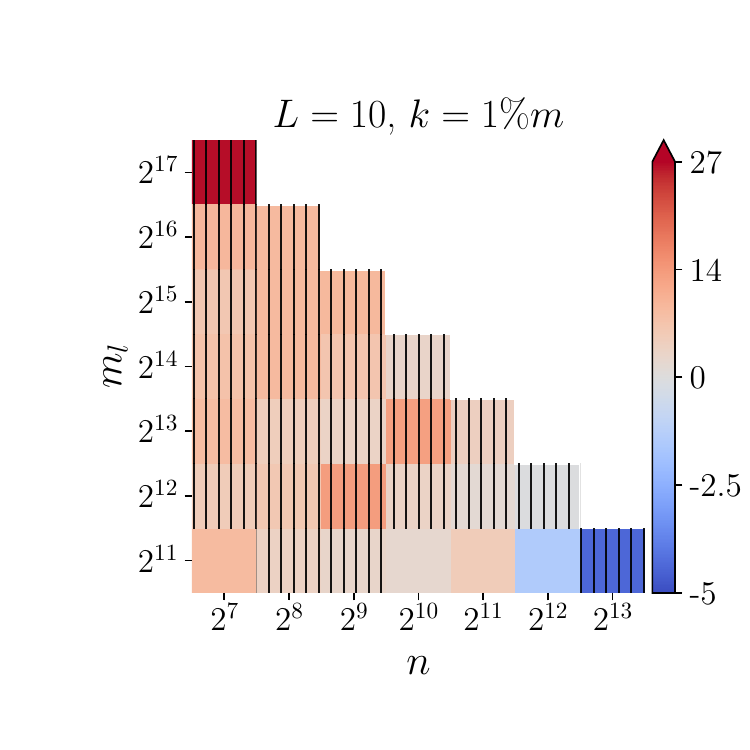}}}
                \caption{\!\!\!\!\!\!\!\!\!\!}
                \label{fig:heatmap_linear:grb_oa:b}
            \end{subsubfigure}
            \begin{subsubfigure}[b]{0.49\linewidth}
                {\resizebox{\linewidth}{!}{\includegraphics[trim=1mm 5mm 1mm 15mm,clip]{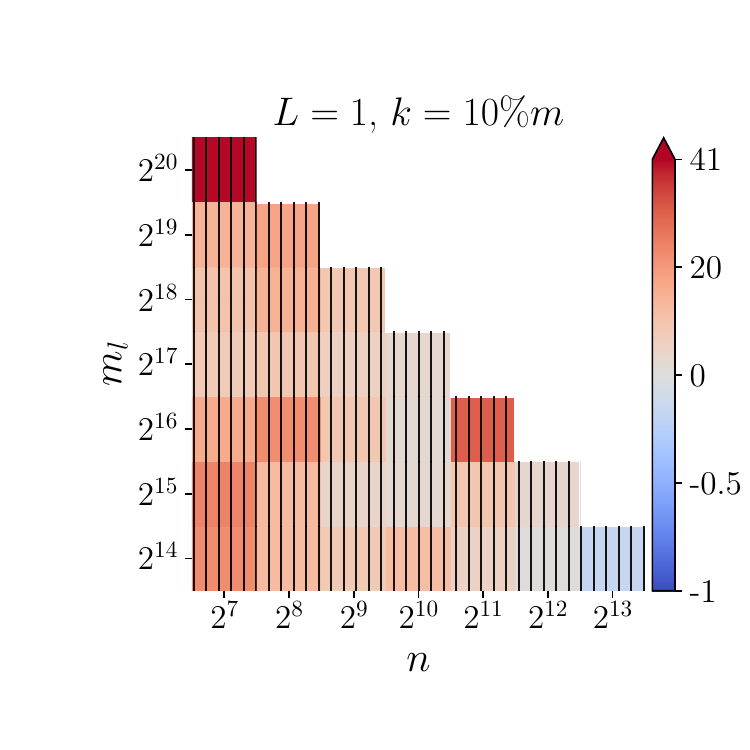}}}
                \caption{\!\!\!\!\!\!\!\!\!\!}
                \label{fig:heatmap_linear:grb_oa:c}
            \end{subsubfigure}
            \begin{subsubfigure}[b]{0.49\linewidth}
                {\resizebox{\linewidth}{!}{\includegraphics[trim=1mm 5mm 1mm 15mm,clip]{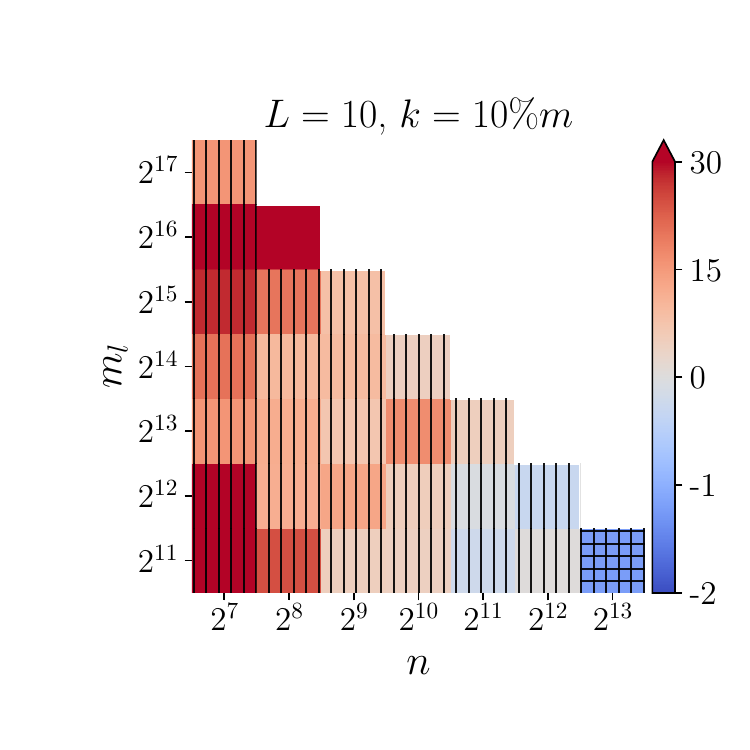}}}
                \caption{\!\!\!\!\!\!\!\!\!\!}
                \label{fig:heatmap_linear:grb_oa:d}
            \end{subsubfigure}
            \captionsetup{labelformat=empty}
            \caption{}
            \label{fig:heatmap_linear:grb_oa}
        \end{minipage}
        \hfill
        \setcounter{subfigure}{1}
        \begin{minipage}[b]{0.49\textwidth}
            \centering
            \caption*{G-OA: quadratic objective} %
            \begin{subsubfigure}[b]{0.49\linewidth}
                {\resizebox{\linewidth}{!}{\includegraphics[trim=1mm 5mm 1mm 15mm,clip]{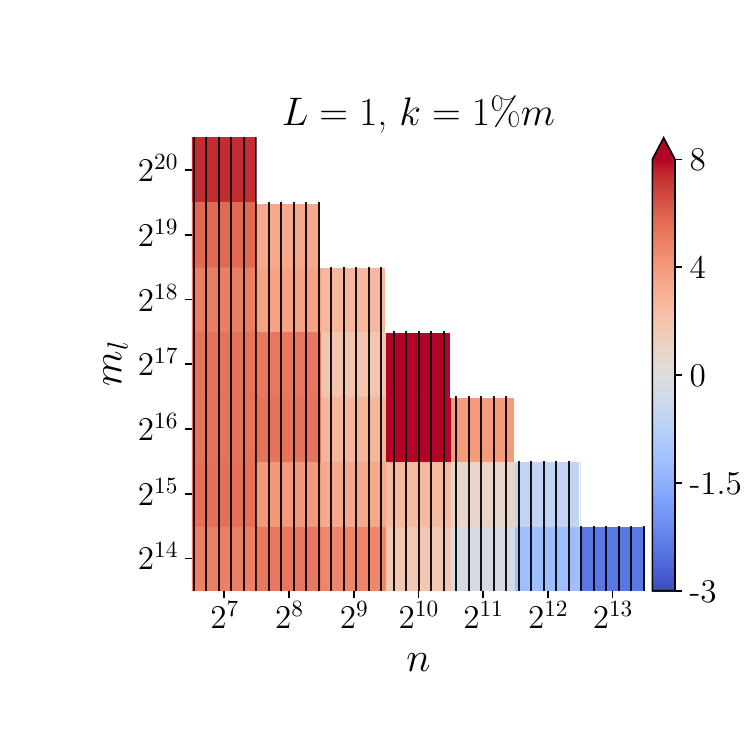}}}
                \caption{\!\!\!\!\!\!\!\!\!\!}
                \label{fig:heatmap_quadratic:grb_oa:a}
            \end{subsubfigure}
            \begin{subsubfigure}[b]{0.49\linewidth}
                {\resizebox{\linewidth}{!}{\includegraphics[trim=1mm 5mm 1mm 15mm,clip]{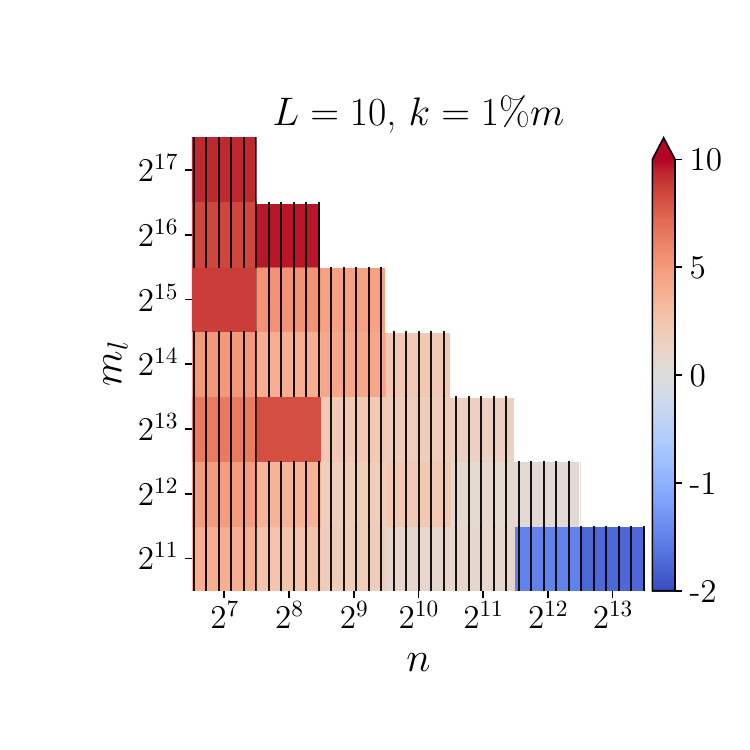}}}
                \caption{\!\!\!\!\!\!\!\!\!\!}
                \label{fig:heatmap_quadratic:grb_oa:b}
            \end{subsubfigure}
            \begin{subsubfigure}[b]{0.49\linewidth}
                {\resizebox{\linewidth}{!}{\includegraphics[trim=1mm 5mm 1mm 15mm,clip]{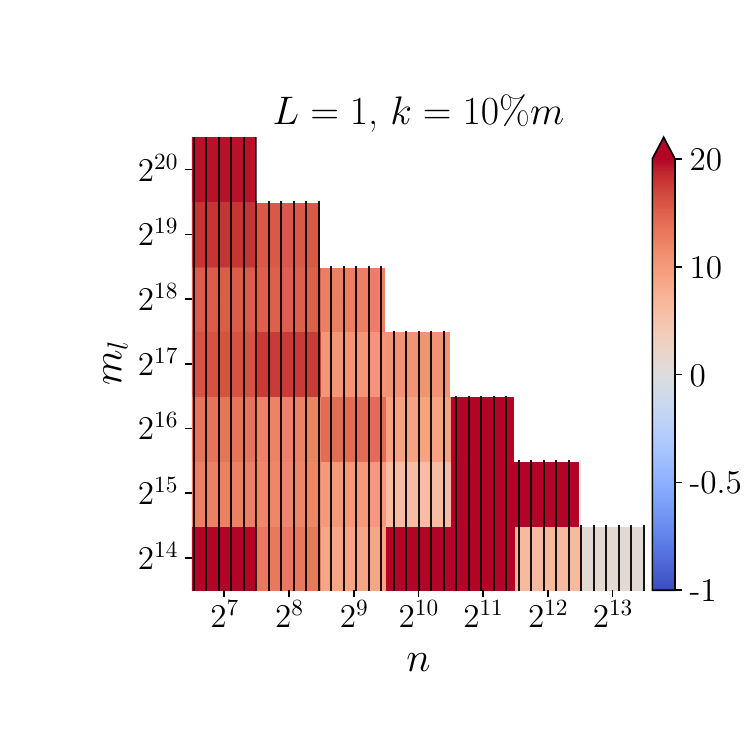}}}
                \caption{\!\!\!\!\!\!\!\!\!\!}
                \label{fig:heatmap_quadratic:grb_oa:c}
            \end{subsubfigure}
            \begin{subsubfigure}[b]{0.49\linewidth}
                {\resizebox{\linewidth}{!}{\includegraphics[trim=1mm 5mm 1mm 15mm,clip]{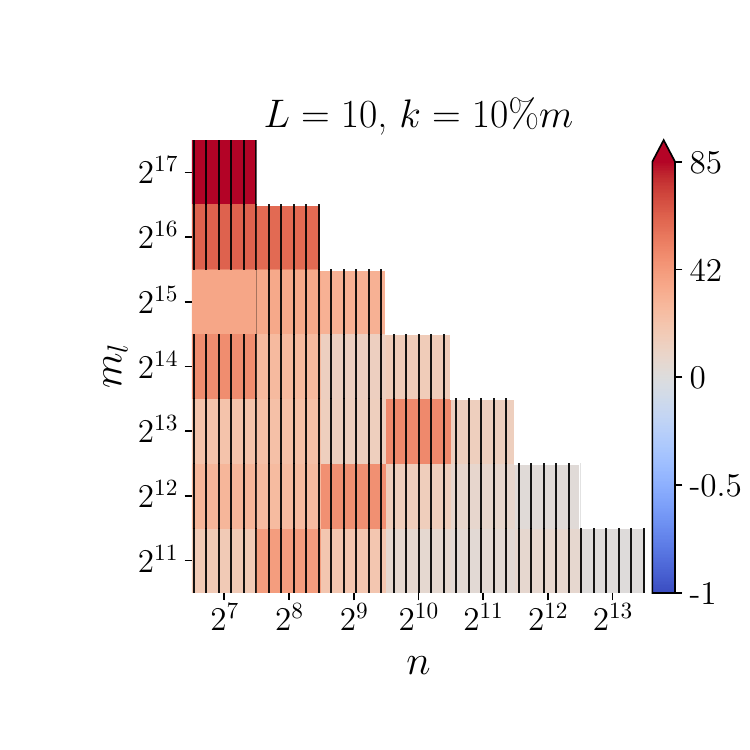}}}
                \caption{\!\!\!\!\!\!\!\!\!\!}
                \label{fig:heatmap_quadratic:grb_oa:d}
            \end{subsubfigure}
            \captionsetup{labelformat=empty}
            \caption{}
            \label{fig:heatmap_quadratic:grb_oa}
        \end{minipage}
        \captionsetup{labelformat=empty}
        \caption{}
        \label{fig:heatmap:grb_oa}
    \end{subfigure}
\end{figure}
\begin{figure}[!ht]
    \ContinuedFloat
    \centering
    \vspace{-6mm}
    \begin{subfigure}[b]{\textwidth}
        \centering
        \setcounter{subfigure}{2}
        \setcounter{subsubfigure}{0}
        
        \begin{minipage}[b]{0.49\textwidth}
            \centering
            \caption*{PSG: linear objective} %
            \begin{subsubfigure}[b]{0.49\linewidth}
                {\resizebox{\linewidth}{!}{\includegraphics[trim=1mm 5mm 1mm 15mm,clip]{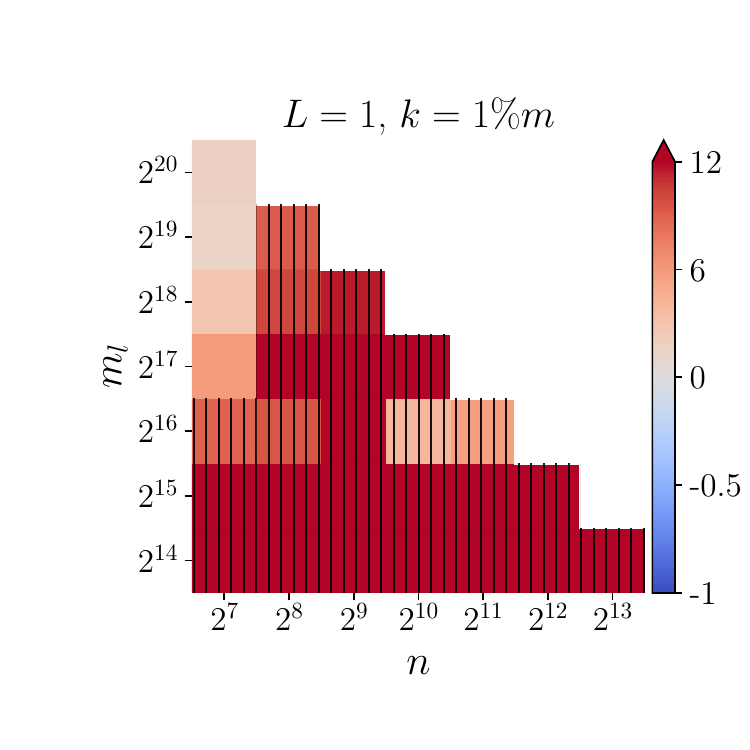}}}
                \caption{\!\!\!\!\!\!\!\!\!\!}
                \label{fig:heatmap_linear:psg:a}
            \end{subsubfigure}
            \begin{subsubfigure}[b]{0.49\linewidth}
                {\resizebox{\linewidth}{!}{\includegraphics[trim=1mm 5mm 1mm 15mm,clip]{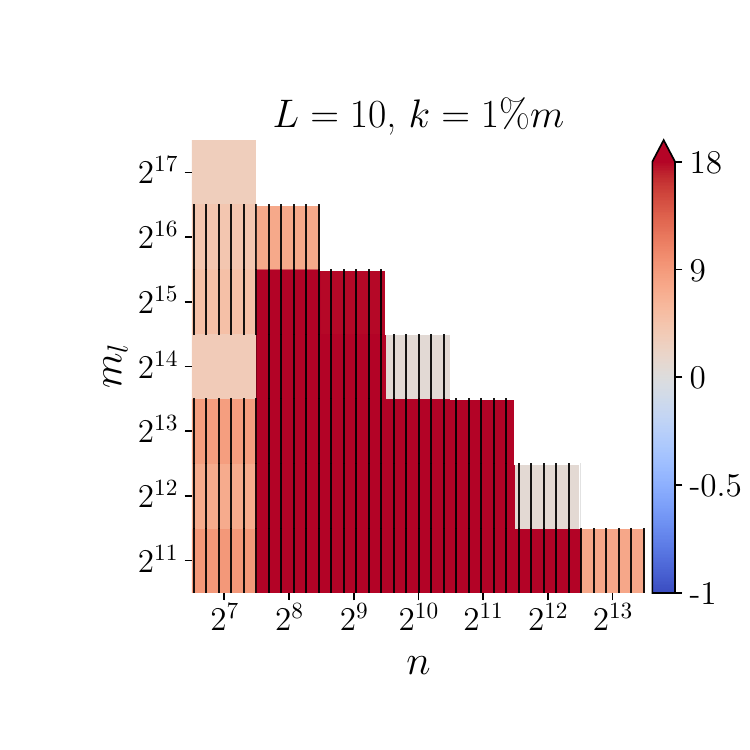}}}
                \caption{\!\!\!\!\!\!\!\!\!\!}
                \label{fig:heatmap_linear:psg:b}
            \end{subsubfigure}
            \begin{subsubfigure}[b]{0.49\linewidth}
                {\resizebox{\linewidth}{!}{\includegraphics[trim=1mm 5mm 1mm 15mm,clip]{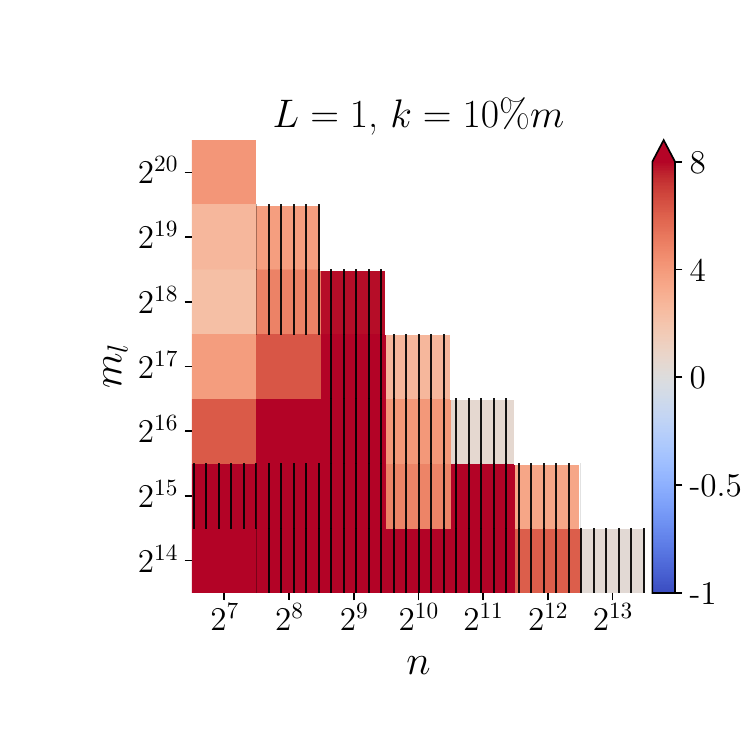}}}
                \caption{\!\!\!\!\!\!\!\!\!\!}
                \label{fig:heatmap_linear:psg:c}
            \end{subsubfigure}
            \begin{subsubfigure}[b]{0.49\linewidth}
                {\resizebox{\linewidth}{!}{\includegraphics[trim=1mm 5mm 1mm 15mm,clip]{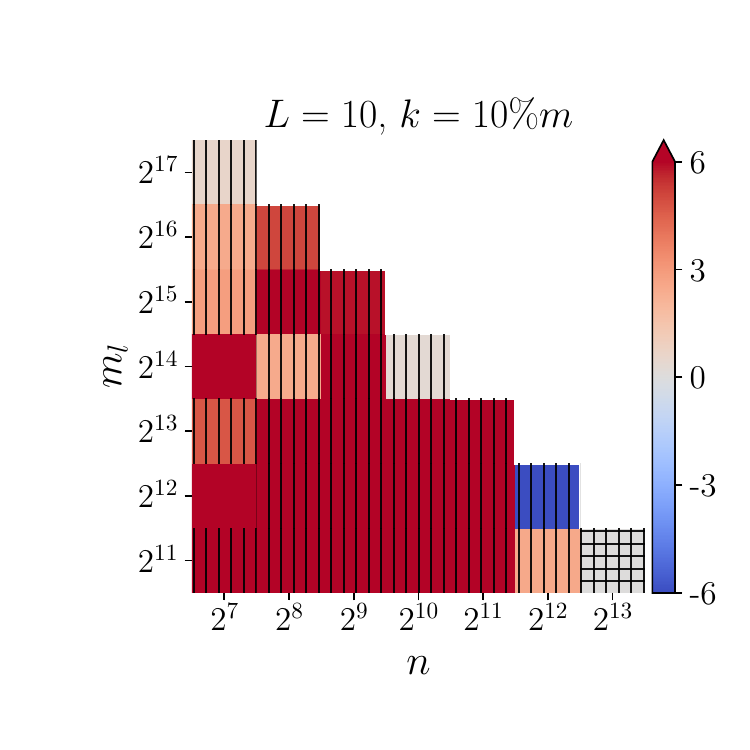}}}
                \caption{\!\!\!\!\!\!\!\!\!\!}
                \label{fig:heatmap_linear:psg:d}
            \end{subsubfigure}
            \captionsetup{labelformat=empty}
            \caption{}
            \label{fig:heatmap_linear:psg}
        \end{minipage}
        \hfill
        \setcounter{subfigure}{2}
        \begin{minipage}[b]{0.49\textwidth}
            \centering
            \caption*{PSG: quadratic objective} %
            \begin{subsubfigure}[b]{0.49\linewidth}
                {\resizebox{\linewidth}{!}{\includegraphics[trim=1mm 5mm 1mm 15mm,clip]{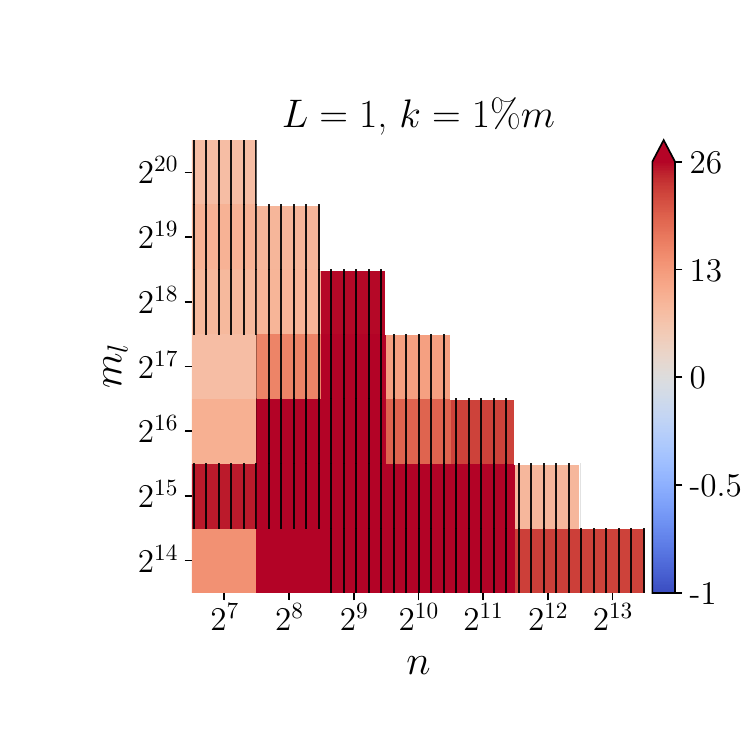}}}
                \caption{\!\!\!\!\!\!\!\!\!\!}
                \label{fig:heatmap_quadratic:psg:a}
            \end{subsubfigure}
            \begin{subsubfigure}[b]{0.49\linewidth}
                {\resizebox{\linewidth}{!}{\includegraphics[trim=1mm 5mm 1mm 15mm,clip]{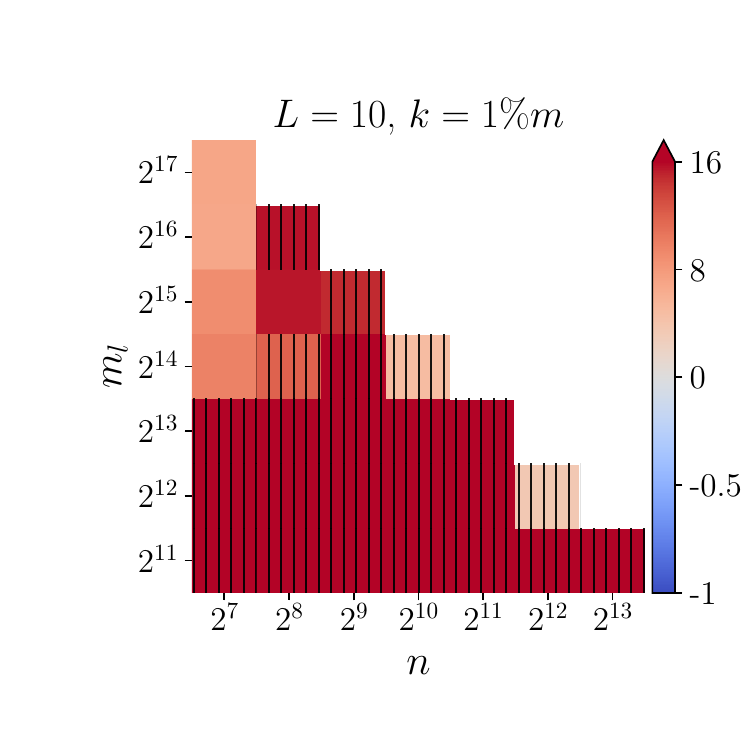}}}
                \caption{\!\!\!\!\!\!\!\!\!\!}
                \label{fig:heatmap_quadratic:psg:b}
            \end{subsubfigure}
            \begin{subsubfigure}[b]{0.49\linewidth}
                {\resizebox{\linewidth}{!}{\includegraphics[trim=1mm 5mm 1mm 15mm,clip]{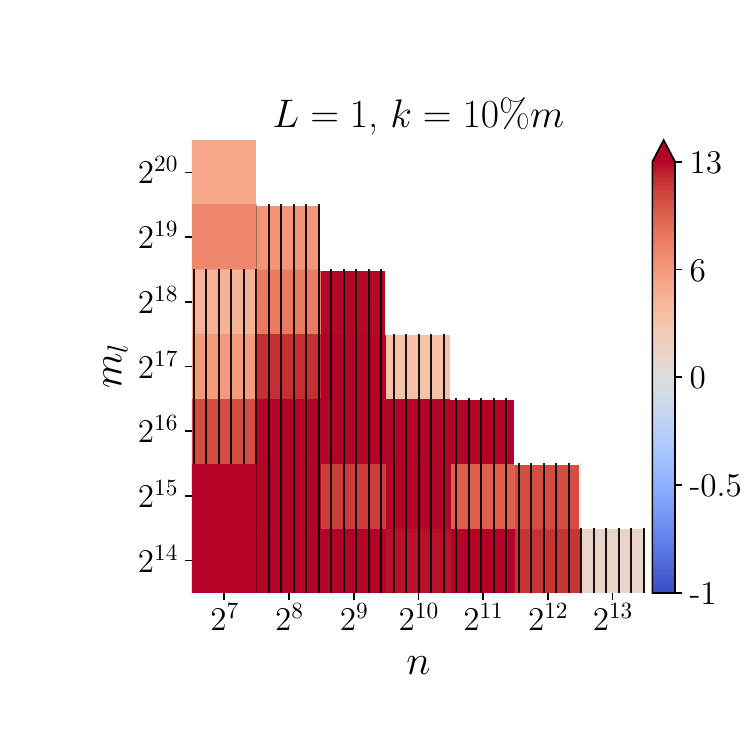}}}
                \caption{\!\!\!\!\!\!\!\!\!\!}
                \label{fig:heatmap_quadratic:psg:c}
            \end{subsubfigure}
            \begin{subsubfigure}[b]{0.49\linewidth}
                {\resizebox{\linewidth}{!}{\includegraphics[trim=1mm 5mm 1mm 15mm,clip]{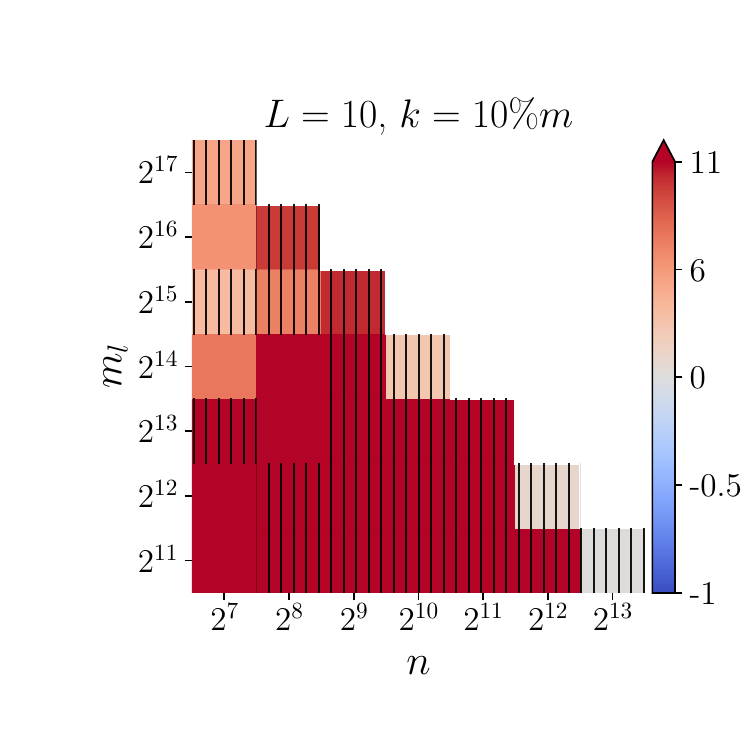}}}
                \caption{\!\!\!\!\!\!\!\!\!\!}
                \label{fig:heatmap_quadratic:psg:d}
            \end{subsubfigure}
            \captionsetup{labelformat=empty}
            \caption{}
            \label{fig:heatmap_quadratic:psg}
        \end{minipage}
        \captionsetup{labelformat=empty}
        \caption{}
        \label{fig:heatmap:psg}
    \end{subfigure}
    \vskip\baselineskip %
\end{figure}
\begin{figure}[ht]
    \ContinuedFloat
    \centering

    \vspace{-6mm}
    \begin{subfigure}[b]{\textwidth}
        \centering
        \setcounter{subfigure}{3}
        \setcounter{subsubfigure}{0}
        
        \begin{minipage}[b]{0.49\textwidth}
            \centering
            \caption*{OSQP: linear objective} %
            \begin{subsubfigure}[b]{0.49\linewidth}
                {\resizebox{\linewidth}{!}{\includegraphics[trim=1mm 5mm 1mm 15mm,clip]{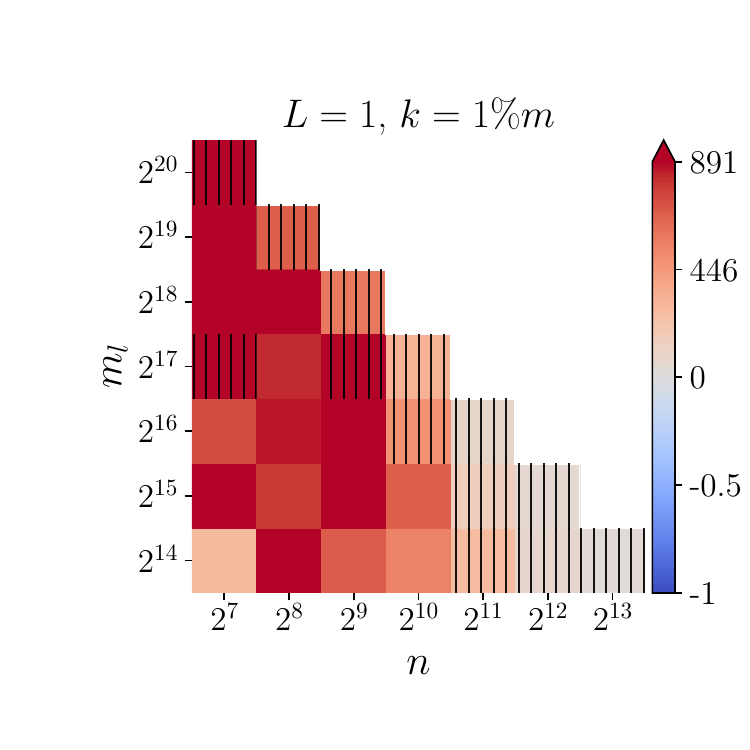}}}
                \caption{\!\!\!\!\!\!\!\!\!\!}
                \label{fig:heatmap_linear:osqp:a}
            \end{subsubfigure}
            \begin{subsubfigure}[b]{0.49\linewidth}
                {\resizebox{\linewidth}{!}{\includegraphics[trim=1mm 5mm 1mm 15mm,clip]{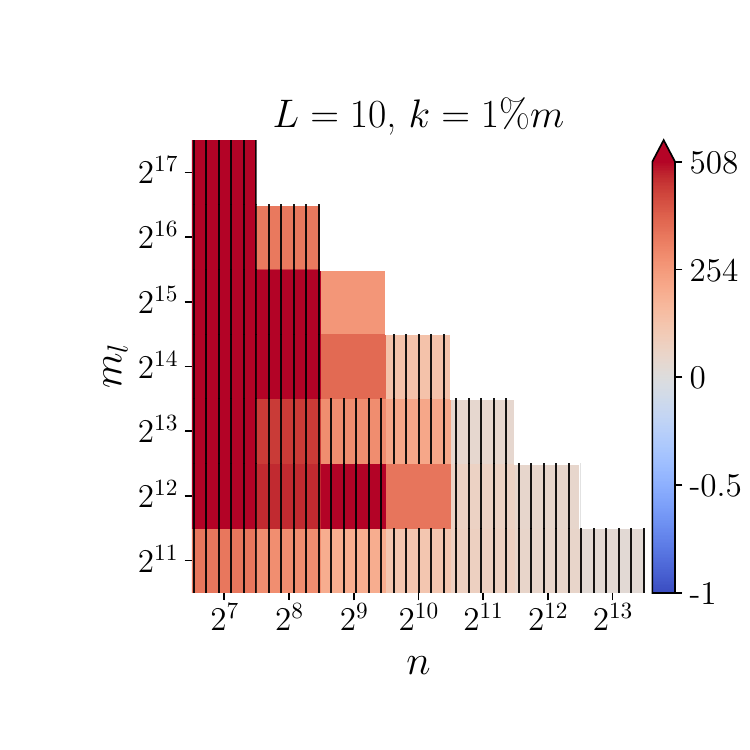}}}
                \caption{\!\!\!\!\!\!\!\!\!\!}
                \label{fig:heatmap_linear:osqp:b}
            \end{subsubfigure}
            \begin{subsubfigure}[b]{0.49\linewidth}
                {\resizebox{\linewidth}{!}{\includegraphics[trim=1mm 5mm 1mm 15mm,clip]{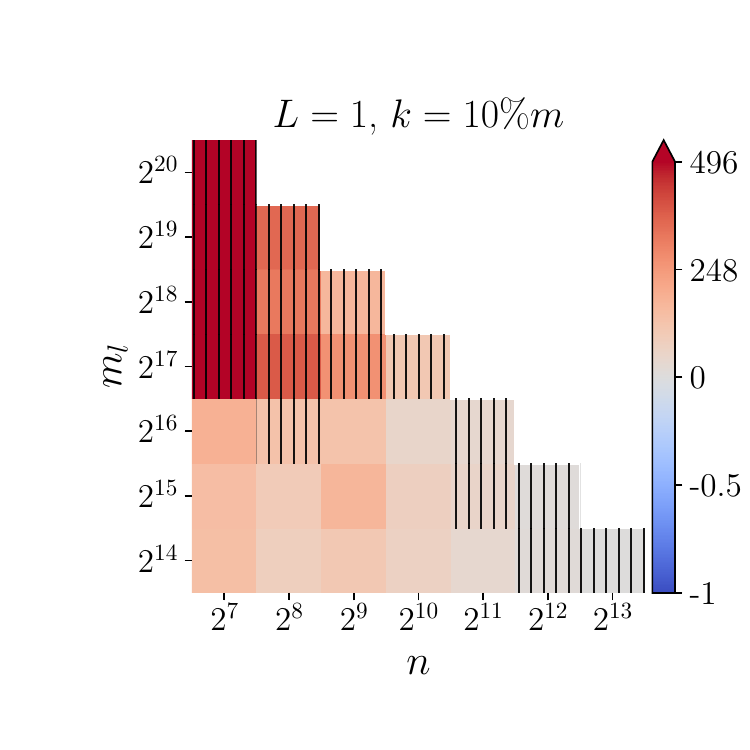}}}
                \caption{\!\!\!\!\!\!\!\!\!\!}
                \label{fig:heatmap_linear:osqp:c}
            \end{subsubfigure}
            \begin{subsubfigure}[b]{0.49\linewidth}
                {\resizebox{\linewidth}{!}{\includegraphics[trim=1mm 5mm 1mm 15mm,clip]{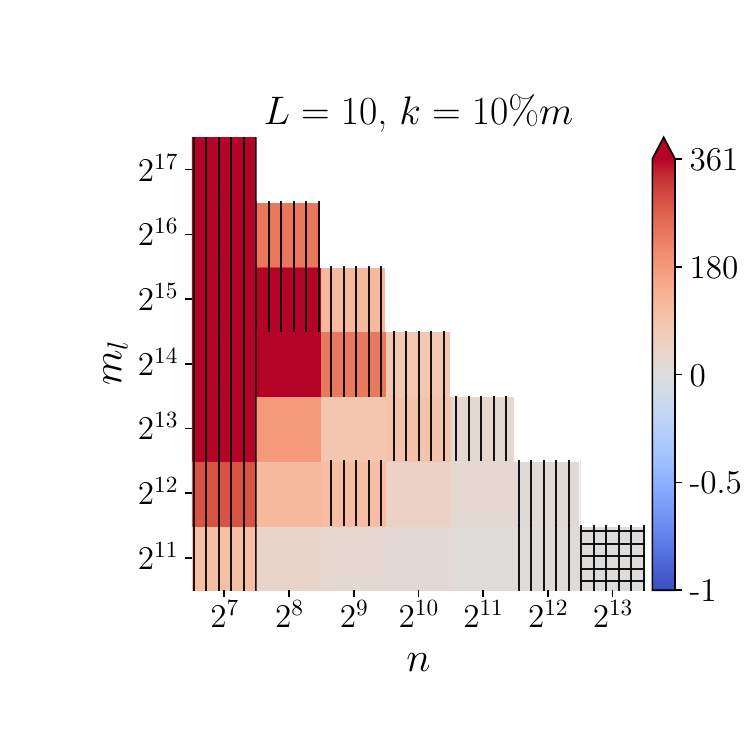}}}
                \caption{\!\!\!\!\!\!\!\!\!\!}
                \label{fig:heatmap_linear:osqp:d}
            \end{subsubfigure}
            \captionsetup{labelformat=empty}
            \caption{}
            \label{fig:heatmap_linear:osqp}
        \end{minipage}
        \hfill
        \setcounter{subfigure}{3}
        \begin{minipage}[b]{0.49\textwidth}
            \centering
            \caption*{OSQP: quadratic objective} %
            \begin{subsubfigure}[b]{0.49\linewidth}
                {\resizebox{\linewidth}{!}{\includegraphics[trim=1mm 5mm 1mm 15mm,clip]{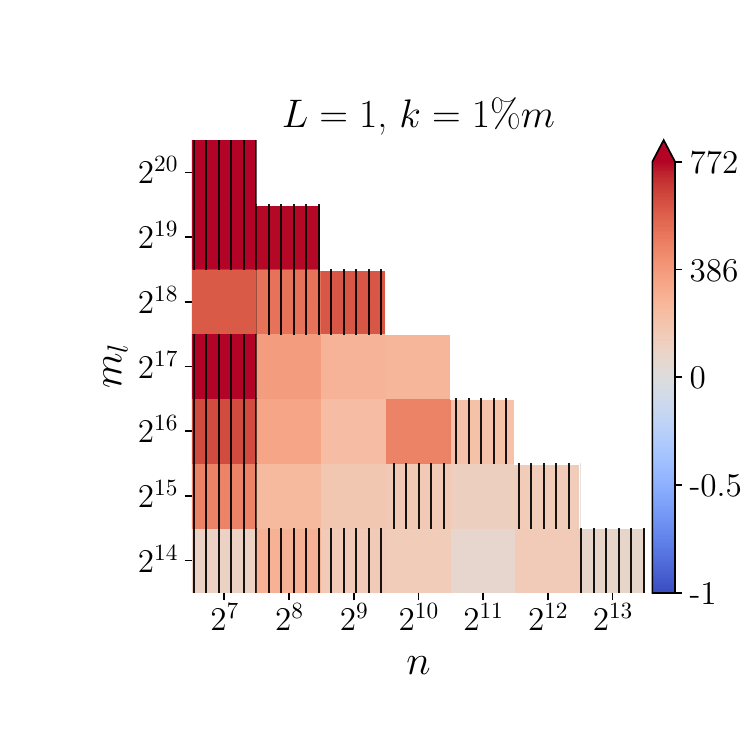}}}
                \caption{\!\!\!\!\!\!\!\!\!\!}
                \label{fig:heatmap_quadratic:osqp:a}
            \end{subsubfigure}
            \begin{subsubfigure}[b]{0.49\linewidth}
                {\resizebox{\linewidth}{!}{\includegraphics[trim=1mm 5mm 1mm 15mm,clip]{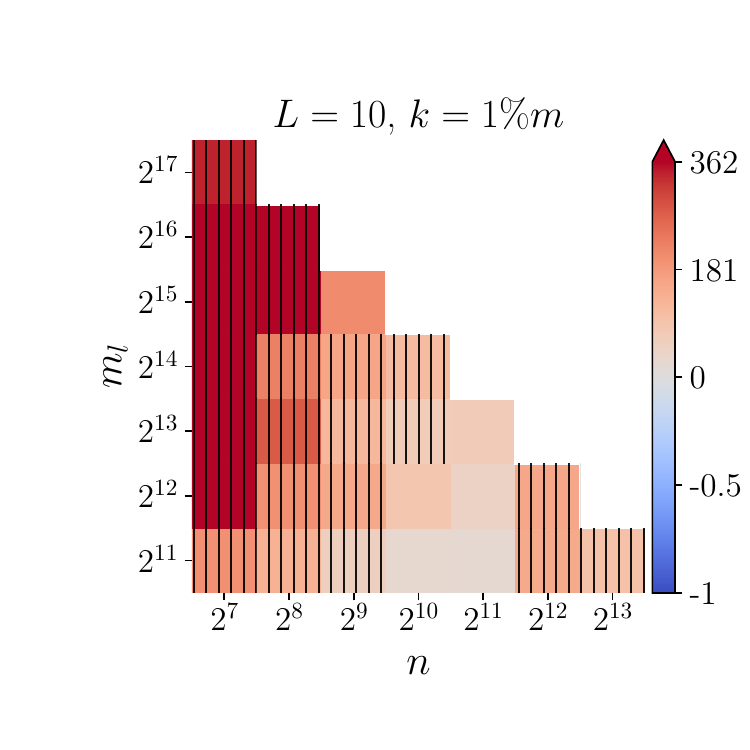}}}
                \caption{\!\!\!\!\!\!\!\!\!\!}
                \label{fig:heatmap_quadratic:osqp:b}
            \end{subsubfigure}
            \begin{subsubfigure}[b]{0.49\linewidth}
                {\resizebox{\linewidth}{!}{\includegraphics[trim=1mm 5mm 1mm 15mm,clip]{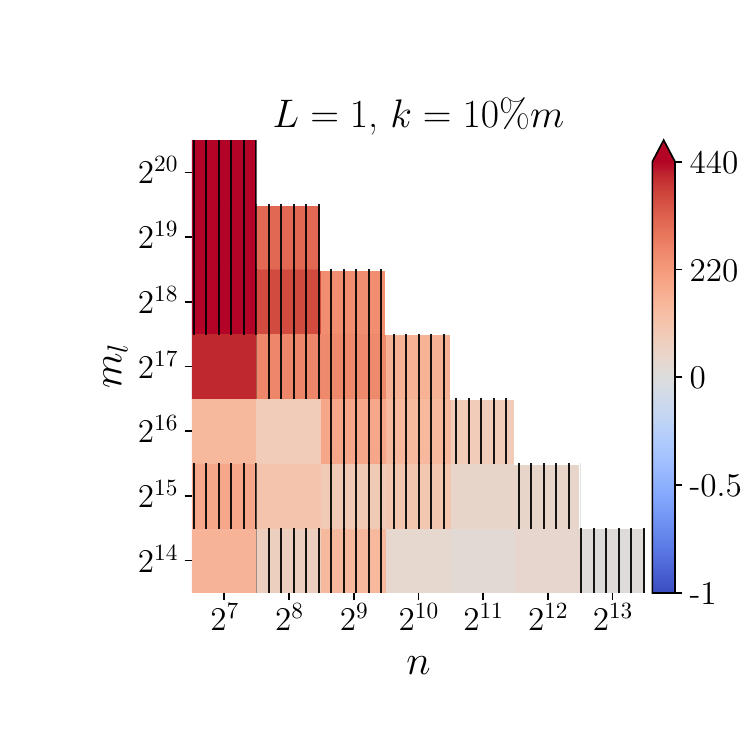}}}
                \caption{\!\!\!\!\!\!\!\!\!\!}
                \label{fig:heatmap_quadratic:osqp:c}
            \end{subsubfigure}
            \begin{subsubfigure}[b]{0.49\linewidth}
                {\resizebox{\linewidth}{!}{\includegraphics[trim=1mm 5mm 1mm 15mm,clip]{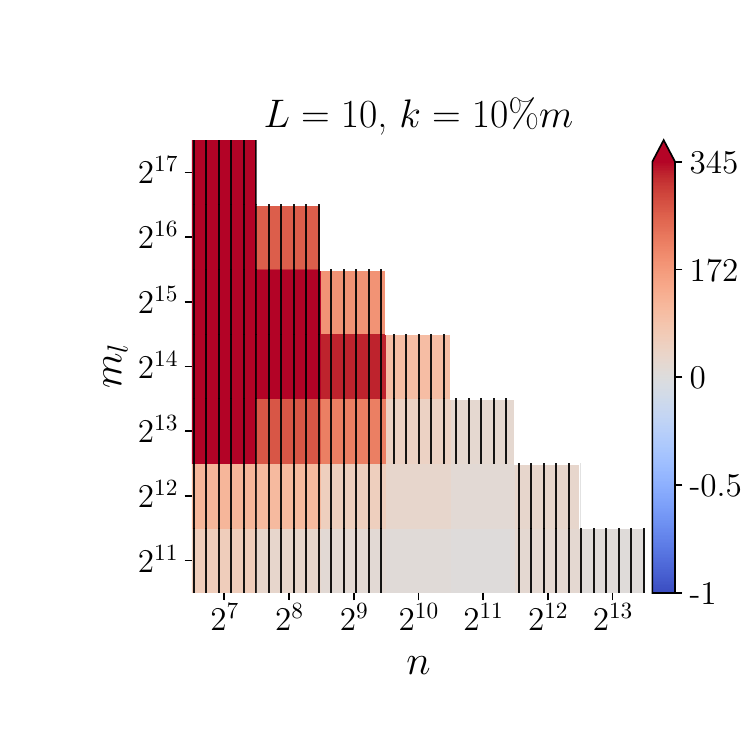}}}
                \caption{\!\!\!\!\!\!\!\!\!\!}
                \label{fig:heatmap_quadratic:osqp:d}
            \end{subsubfigure}
            \captionsetup{labelformat=empty}
            \caption{}
            \label{fig:heatmap_quadratic:osqp}
        \end{minipage}
        \captionsetup{labelformat=empty}
        \caption{}
        \label{fig:heatmap:osqp}
    \end{subfigure}
    \caption{Heatmap of relative solve time (method $i$ $-$ ALM)/ALM, for $i\in\{$a = GRB, b = G-OA, c = PSG, d = OSQP$\}$.
    Subfigures (i)--(iv) ((v)--(vii)) correspond to instances with linear (separable quadratic) objective.
    Red shading (positive) indicates that ALM was faster than method $i$; blue shading (negative) indicates that method $i$ was faster than ALM.
    Shading intensity represents relative solve time.
    Vertical (horizontal) hatching indicates that method $i$ (ALM) failed to obtain a satisfactory solution within 1hr.
    For convenience: $2^{20}=1.0\times10^6$, $2^{17}=1.3\times10^5$, and $2^{13}=0.8\times10^4$.}
    \label{fig:synthetic_heatmap}
\end{figure}

\cref{fig:heatmap_linear:psg} shows that in some instances, PSG can obtain solutions to the large-$m$ cases very fast, almost matching ALM on problems with linear objective (left two columns).
However, PSG failed to obtain solutions with satisfactory primal feasibility in many instances.
Finally, \cref{fig:heatmap_linear:osqp} shows that ALM ($\epsilon=10^{-3})$ was at least 50x faster, and sometimes more than 750x faster, than OSQP ($\epsilon=10^{-3}$) in large-$m$ instances.

\cref{tab:synthetic} provides additional detail on the comparison of ALM against other methods across the three instances with the largest $m$.
It shows that PSG can outperform GRB and GRB-OA in certain instances, but PSG can also fail and return solutions with significant primal feasibility violation (\eg the third subproblem in each group).
In addition, the table indicates that G-OA finds feasible solutions that are nearly optimal.

Detailed timing statistics for ALM on the ten largest-$m$ instances (top left region of subplots in \cref{fig:synthetic_heatmap}) are summarized in \cref{tab:detail_synthetic}.
Notably, the incremental cost of obtaining second-order information (after computing $\nabla\varphi_\sigma$) can be 3-7x \emph{cheaper} than the incremental cost of obtaining first-order information (after computing $\varphi_\sigma$) when $m\gg n$ and $k\ll m$.
However, this cost advantage lessens as the ratio  $m/n$ decreases.

In addition, comparison of the three timing columns indicates several key insights. 
First, solutions with the accuracy around $\eta=10^{-3}$  can be obtained in approximately the same amount of time for both linear and quadratic objectives.
Second, ALM is highly effective at refining a  solution from moderate precision ($\eta=10^{-6}$) to  high precision ($\eta=10^{-8}$).
This stands in contrast to barrier methods, which  often encounter difficulties in the later stages of solving large-scale instances due to poor numerical conditioning. The efficiency of ALM can be attributed to the small cardinality of ${\bar\beta}$ near a solution, as noted in the last column, which enhances performance particularly when the set of effective variables remains consistent, thereby minimizing the need for frequent sorting. The penultimate column highlights that the matrix $\widetilde{T}^\top \widetilde{T}$ is constructed in a space considerably smaller than $m$, and the Newton system can be solved in such a small dimensional space.

\begin{sidewaystable}
\begin{subtable}{\linewidth}
    \hspace*{-6mm}%
    \begin{adjustbox}{max width=1.05\textwidth}
    \begin{tabular}{llllllllllllllllllllllllllll}
\toprule 
\multicolumn{1}{c}{Problem} &  & \multicolumn{4}{c}{Relative Pobj} &  & \multicolumn{5}{c}{Pinfeas} &  & \multicolumn{4}{c}{Dinfeas} &  & \multicolumn{4}{c}{Relgap} &  & \multicolumn{5}{c}{Walltime [s]} \\
\cmidrule{1-1} \cmidrule{3-6} \cmidrule{8-12} \cmidrule{14-17} \cmidrule{19-22} \cmidrule{24-28} ($m,n$) &  & GRB & G-OA & PSG & OSQP &  & ALM & GRB & G-OA & PSG & OSQP &  & ALM & GRB & G-OA & OSQP &  & ALM & GRB & G-OA & OSQP &  & ALM & GRB & G-OA & PSG & OSQP \\
\midrule 
\midrule 
\multicolumn{28}{c}{$L={1},\;k={1}\%m$} \\
\midrule 
$(2^{20}, 2^{7})$ &  & -2.6e-11 & 3.6e-8 & 5.3e-9 & $-$ &  & 3.3e-9 & 5.2e-12 & 0.0e0 & 7.0e-9 & $-$ &  & 3.3e-11 & 2.8e-14 & 1.3e-14 & $-$ &  & 7.3e-11 & -1.9e-12 & 2.2e-6 & $-$ &  & 2.6e1 & 3.8e2 & 1.3e2 & 6.9e1 & 3.6e3 \\
$(2^{17}, 2^{7})$ &  & 4.1e-10 & -2.5e-6 & 4.0e-8 & -7.8e-4 &  & 2.3e-9 & 0.0e0 & 7.2e-6 & 7.4e-9 & 2.3e-3 &  & 8.2e-10 & 1.1e-14 & 1.4e-14 & 1.9e-4 &  & 1.8e-10 & 1.8e-9 & 5.4e-7 & -7.5e-4 &  & 1.6e0 & 1.9e1 & 2.8e1 & 1.1e1 & 8.8e2 \\
$(2^{17}, 2^{10})$ &  & 8.6e-10 & -2.7e-6 & 2.3e-5 & -4.0e-4 &  & 9.4e-9 & 0.0e0 & 3.5e-5 & 4.4e-8 & 5.4e-3 &  & 1.1e-10 & 5.4e-11 & 1.3e-10 & 2.3e-4 &  & -4.3e-10 & 7.5e-10 & -1.6e-6 & -3.9e-4 &  & 2.3e1 & 2.3e2 & 7.0e1 & 6.9e2 & 3.6e3 \\
\midrule 
\multicolumn{28}{c}{$L={1},\;k={10}\%m$} \\
\midrule 
$(2^{20}, 2^{7})$ &  & 6.0e-9 & -4.6e-6 & 1.1e-9 & $-$ &  & 7.5e-9 & 0.0e0 & 8.0e-6 & 7.8e-9 & $-$ &  & 2.1e-10 & 2.9e-14 & 1.4e-14 & $-$ &  & -1.3e-9 & 5.2e-9 & -3.8e-6 & $-$ &  & 2.7e1 & 2.2e2 & 1.1e3 & 1.4e2 & 3.6e3 \\
$(2^{17}, 2^{7})$ &  & 6.4e-9 & -6.7e-5 & 1.4e-8 & -3.2e-3 &  & 6.7e-11 & 0.0e0 & 1.3e-4 & 5.0e-9 & 6.0e-3 &  & 7.5e-11 & 7.4e-15 & 4.5e-15 & 7.5e-5 &  & 1.5e-12 & 6.3e-12 & -6.4e-5 & -3.1e-3 &  & 2.9e0 & 1.7e1 & 2.5e1 & 1.4e1 & 5.5e2 \\
$(2^{17}, 2^{10})$ &  & 7.7e-10 & -9.3e-8 & -2.8e-1 & 1.0e-2 &  & 9.6e-9 & 3.3e-12 & 7.8e-7 & 2.8e0 & 4.6e-3 &  & 7.0e-11 & 2.0e-12 & 3.2e-11 & 1.9e-4 &  & 2.8e-11 & 1.4e-9 & 6.7e-8 & 1.0e-2 &  & 5.6e1 & 1.1e2 & 2.3e2 & 2.0e2 & 3.6e3 \\
\midrule 
\multicolumn{28}{c}{$L={10},\;k={1}\%m$} \\
\midrule 
$(2^{17}, 2^{7})$ &  & 5.5e-11 & -4.6e-8 & 9.4e-9 & $-$ &  & 2.0e-9 & 1.4e-11 & 1.4e-7 & 8.6e-10 & $-$ &  & 3.0e-10 & 1.0e-12 & 1.8e-13 & $-$ &  & 3.9e-10 & 5.3e-10 & -4.4e-8 & $-$ &  & 2.5e1 & 6.7e2 & 7.0e2 & 1.0e2 & 3.6e3 \\
$(2^{14}, 2^{7})$ &  & 3.8e-9 & -2.3e-6 & 1.3e-7 & -2.5e-4 &  & 4.3e-7 & 0.0e0 & 8.6e-6 & 1.7e-14 & 2.3e-3 &  & 8.7e-12 & 2.2e-12 & 1.8e-13 & 2.8e-4 &  & -1.2e-9 & 2.2e-9 & -2.3e-6 & -2.8e-4 &  & 4.6e0 & 6.1e1 & 3.7e1 & 2.0e1 & 9.9e2 \\
$(2^{14}, 2^{10})$ &  & -7.8e-11 & 3.3e-8 & 8.0e-3 & -5.8e-5 &  & 4.5e-9 & 0.0e0 & 2.0e-7 & 8.7e-1 & 3.0e-3 &  & 1.7e-11 & 1.4e-12 & 4.0e-12 & 6.9e-4 &  & 1.6e-10 & 1.7e-10 & 5.3e-8 & -5.0e-5 &  & 4.2e1 & 2.8e2 & 1.6e2 & 7.6e1 & 3.6e3 \\
\midrule 
\multicolumn{28}{c}{$L={10},\;k={10}\%m$} \\
\midrule 
$(2^{17}, 2^{7})$ &  & 5.0e-10 & -1.5e-7 & 1.1e-3 & -9.1e-3 &  & 6.1e-9 & 1.5e-9 & 3.4e-7 & 6.9e-9 & 2.2e-2 &  & 3.7e-11 & 7.9e-14 & 2.7e-14 & 5.8e-4 &  & -3.8e-10 & 2.7e-9 & -1.2e-7 & -3.0e-3 &  & 6.0e1 & 4.9e2 & 1.0e3 & 9.5e1 & 3.6e3 \\
$(2^{14}, 2^{7})$ &  & 1.6e-8 & -4.4e-8 & 3.2e-7 & -1.2e-3 &  & 5.8e-10 & 0.0e0 & 1.6e-7 & 9.2e-9 & 3.8e-3 &  & 4.2e-9 & 7.7e-13 & 1.8e-13 & 2.7e-4 &  & 4.3e-11 & 1.2e-10 & 3.8e-8 & -1.1e-3 &  & 4.8e0 & 4.9e1 & 1.0e2 & 3.4e1 & 1.9e3 \\
$(2^{14}, 2^{10})$ &  & 4.1e-9 & 3.9e-7 & -1.1e-1 & -5.6e-4 &  & 8.7e-9 & 3.6e-12 & 0.0e0 & 2.0e0 & 1.8e-2 &  & 8.0e-11 & 7.0e-13 & 9.7e-11 & 4.4e-4 &  & -7.8e-10 & 4.0e-9 & 9.5e-7 & -5.1e-4 &  & 8.2e1 & 2.0e2 & 4.6e2 & 1.1e2 & 3.6e3 \\
\bottomrule 
\end{tabular}
    \end{adjustbox}
    \caption{Linear objective.}
    \label{tab:comparison_synthetic_linear}
\end{subtable}
\newline
\vspace*{0.2cm}
\newline
\begin{subtable}{\linewidth}
    \hspace*{-6mm}
    \begin{adjustbox}{max width=1.05\textwidth}
    \begin{tabular}{llllllllllllllllllllllllllll}
\toprule 
\multicolumn{1}{c}{Problem} &  & \multicolumn{4}{c}{Relative Pobj} &  & \multicolumn{5}{c}{Pinfeas} &  & \multicolumn{4}{c}{Dinfeas} &  & \multicolumn{4}{c}{Relgap} &  & \multicolumn{5}{c}{Walltime [s]} \\
\cmidrule{1-1} \cmidrule{3-6} \cmidrule{8-12} \cmidrule{14-17} \cmidrule{19-22} \cmidrule{24-28} ($m,n$) &  & GRB & G-OA & PSG & OSQP &  & ALM & GRB & G-OA & PSG & OSQP &  & ALM & GRB & G-OA & OSQP &  & ALM & GRB & G-OA & OSQP &  & ALM & GRB & G-OA & PSG & OSQP \\
\midrule 
\midrule 
\multicolumn{28}{c}{$L={1},\;k={1}\%m$} \\
\midrule 
$(2^{20}, 2^{7})$ &  & 3.2e-10 & 2.0e-7 & 7.3e-11 & $-$ &  & 1.4e-9 & 0.0e0 & 0.0e0 & 1.4e-8 & $-$ &  & 1.2e-11 & 4.9e-14 & 8.3e-15 & $-$ &  & -1.0e-10 & 4.2e-10 & 3.9e-7 & $-$ &  & 1.1e1 & 5.9e2 & 9.2e1 & 9.5e1 & 3.6e3 \\
$(2^{17}, 2^{7})$ &  & 3.7e-9 & 4.7e-7 & 3.0e-9 & -5.3e-4 &  & 9.7e-9 & 0.0e0 & 0.0e0 & 6.1e-14 & 2.0e-3 &  & 5.5e-10 & 8.3e-15 & 1.0e-14 & 1.5e-5 &  & -1.2e-9 & 4.9e-9 & 9.2e-7 & -5.2e-4 &  & 9.5e-1 & 3.0e1 & 6.1e0 & 8.1e0 & 4.4e2 \\
$(2^{17}, 2^{10})$ &  & 9.9e-10 & 4.2e-8 & -5.0e-2 & -5.3e-5 &  & 1.1e-9 & 0.0e0 & 0.0e0 & 4.5e0 & 7.1e-4 &  & 4.1e-11 & 5.1e-13 & 6.3e-10 & 7.2e-6 &  & -7.5e-10 & 4.6e-10 & 8.3e-8 & -5.1e-5 &  & 1.2e1 & 1.5e2 & 2.6e2 & 1.5e2 & 2.4e3 \\
\midrule 
\multicolumn{28}{c}{$L={1},\;k={10}\%m$} \\
\midrule 
$(2^{20}, 2^{7})$ &  & 3.2e-9 & 1.0e-5 & 1.7e-9 & $-$ &  & 7.7e-9 & 0.0e0 & 0.0e0 & 8.2e-9 & $-$ &  & 2.6e-11 & 2.8e-14 & 1.5e-14 & $-$ &  & -1.9e-9 & 2.4e-9 & 2.0e-5 & $-$ &  & 1.8e1 & 2.4e2 & 3.7e2 & 1.2e2 & 3.6e3 \\
$(2^{17}, 2^{7})$ &  & 4.7e-9 & 8.0e-7 & 4.4e-10 & -3.2e-4 &  & 2.1e-9 & 0.0e0 & 0.0e0 & 1.2e-8 & 7.4e-4 &  & 1.6e-11 & 8.9e-15 & 8.6e-15 & 4.2e-5 &  & -4.3e-10 & 8.3e-9 & 1.6e-6 & -3.3e-4 &  & 2.0e0 & 2.1e1 & 3.3e1 & 1.4e1 & 4.0e2 \\
$(2^{17}, 2^{10})$ &  & 5.5e-9 & 3.5e-7 & -7.4e-2 & -1.7e-4 &  & 6.8e-9 & 0.0e0 & 0.0e0 & 3.4e0 & 1.7e-3 &  & 9.7e-10 & 3.7e-14 & 2.1e-14 & 1.1e-5 &  & -2.6e-9 & 5.7e-9 & 6.8e-7 & -1.7e-4 &  & 2.7e1 & 1.7e2 & 3.1e2 & 1.2e2 & 3.2e3 \\
\midrule 
\multicolumn{28}{c}{$L={10},\;k={1}\%m$} \\
\midrule 
$(2^{17}, 2^{7})$ &  & 4.5e-10 & 9.1e-8 & 3.2e-10 & -1.1e-2 &  & 3.3e-9 & 0.0e0 & 0.0e0 & 8.2e-14 & 5.6e-2 &  & 5.6e-11 & 8.2e-13 & 1.4e-14 & 1.3e-3 &  & -3.0e-10 & 1.3e-9 & 5.9e-7 & -2.1e-3 &  & 1.5e1 & 1.1e3 & 1.5e2 & 1.2e2 & 3.6e3 \\
$(2^{14}, 2^{7})$ &  & 9.0e-10 & 2.2e-6 & 2.1e-9 & -1.7e-3 &  & 4.6e-9 & 0.0e0 & 0.0e0 & 1.3e-14 & 1.4e-2 &  & 1.5e-10 & 8.5e-13 & 1.5e-14 & 2.8e-4 &  & -2.5e-10 & 2.2e-9 & 7.0e-6 & -1.5e-3 &  & 1.6e0 & 5.2e1 & 9.7e0 & 1.7e1 & 5.4e2 \\
$(2^{14}, 2^{10})$ &  & 3.2e-9 & 1.9e-8 & -5.2e-2 & -9.6e-5 &  & 6.1e-9 & 0.0e0 & 0.0e0 & 6.5e0 & 9.2e-3 &  & 1.2e-9 & 1.7e-13 & 2.3e-14 & 1.3e-4 &  & -3.1e-9 & 4.7e-11 & 3.6e-8 & -9.4e-5 &  & 2.3e1 & 2.8e2 & 7.1e1 & 1.3e2 & 2.0e3 \\
\midrule 
\multicolumn{28}{c}{$L={10},\;k={10}\%m$} \\
\midrule 
$(2^{17}, 2^{7})$ &  & 1.5e-9 & 1.6e-8 & 2.1e-3 & -4.7e-3 &  & 6.3e-9 & 0.0e0 & 0.0e0 & 1.1e-8 & 1.7e-2 &  & 2.3e-11 & 5.2e-14 & 4.3e-14 & 1.0e-4 &  & -1.3e-9 & 6.9e-10 & 7.2e-8 & -6.5e-4 &  & 2.0e1 & 5.6e2 & 1.7e3 & 1.1e2 & 3.6e3 \\
$(2^{14}, 2^{7})$ &  & 4.6e-10 & 4.0e-8 & 1.5e-10 & -7.2e-4 &  & 1.2e-9 & 0.0e0 & 0.0e0 & 2.2e-9 & 2.3e-3 &  & 1.4e-10 & 1.7e-14 & 8.3e-15 & 3.0e-4 &  & -1.4e-10 & 4.4e-9 & 3.9e-7 & -6.7e-4 &  & 2.2e0 & 4.0e1 & 1.1e2 & 1.8e1 & 1.4e3 \\
$(2^{14}, 2^{10})$ &  & 2.6e-9 & 1.3e-8 & -6.4e-2 & -6.3e-4 &  & 7.9e-10 & 0.0e0 & 0.0e0 & 3.6e0 & 2.0e-2 &  & 9.4e-11 & 1.6e-14 & 4.6e-15 & 1.3e-4 &  & -2.7e-10 & 5.3e-9 & 2.8e-8 & -5.6e-4 &  & 4.9e1 & 3.3e2 & 8.2e2 & 1.7e2 & 3.6e3 \\
\bottomrule 
\end{tabular}
    \end{adjustbox}
    \caption{Diagonal quadratic objective.}
    \label{tab:comparison_synthetic_quadratic}
\end{subtable}
\caption{Performance of different methods for synthetic data. {\sl Relative Pobj} $\coloneqq (\textrm{Pobj} - \textrm{Pobj}_\textrm{ALM}) / \abs{\textrm{Pobj}_\textrm{ALM}}$, and a positive value indicates that the method achieved a higher objective value  than ALM. Tolerance for OSQP is set at $\epsilon=10^{-3}$, while for all others are set  to $10^{-8}$. Note:
PSG does not provide metrics for dual infeasibility or relative gap.}
\label{tab:synthetic}
\end{sidewaystable}

\begin{sidewaystable}
\begin{subtable}{\linewidth}
    \hspace*{-5mm}%
        \begin{adjustbox}{max width=1.05\textwidth}
        \begin{tabular}{ccccccccccccccccccccccccccccccc}
\toprule 
\multicolumn{1}{c}{Problem} &  & \multicolumn{2}{c}{KKT residual} &  & \multicolumn{2}{c}{$\text{Time}_{\text{a}}: \epsilon=10^{-3}$} &  & \multicolumn{2}{c}{$\text{Time}_{\text{b}}: \epsilon=10^{-6}$} &  & \multicolumn{2}{c}{$\text{Time}_{\text{c}}: \epsilon=10^{-8}$} &  & \multicolumn{2}{c}{\% $\text{Time}_{\text{c}}$: Sort} &  & \multicolumn{2}{c}{\% $\text{Time}_{\text{c}}$: $V^{-1}v$} &  & \multicolumn{2}{c}{\% $\text{Time}_{\text{c}}$: ${\nabla}\varphi$} &  & \multicolumn{2}{c}{$\text{ALM }\mid\text{ SSN iter}$} &  & \multicolumn{2}{c}{$\text{avg }|\bar\beta|$} &  & \multicolumn{2}{c}{$\text{final }|\bar\beta|$} \\
\cmidrule{3-4} \cmidrule{6-7} \cmidrule{9-10} \cmidrule{12-13} \cmidrule{15-16} \cmidrule{18-19} \cmidrule{21-22} \cmidrule{24-25} \cmidrule{27-28} \cmidrule{30-31} ($m,n$) &  & lin. & quad. &  & lin. & quad. &  & lin. & quad. &  & lin. & quad. &  & lin. & quad. &  & lin. & quad. &  & lin. & quad. &  & lin. & quad. &  & lin. & quad. &  & lin. & quad. \\
\midrule 
\midrule 
\multicolumn{31}{c}{$L={1},\;k={1}\%m$} \\
\midrule 
$(2^{20}, 2^{7})$ &  & 3.3e-9 & 1.4e-9 &  & 4.0e0 & 4.7e0 &  & 1.0e1 & 8.8e0 &  & 2.6e1 & 1.1e1 &  & 42 & 16 &  & 4 & 6 &  & 12 & 20 &  & 18 $\mid$ 100 & 17 $\mid$ 70$\phantom{0}$ &  & 380 & 290 &  & 110 & 80 \\
$(2^{19}, 2^{7})$ &  & 6.7e-9 & 2.0e-9 &  & 1.9e0 & 1.9e0 &  & 7.3e0 & 3.8e0 &  & 1.3e1 & 4.4e0 &  & 34 & 14 &  & 5 & 7 &  & 16 & 21 &  & 19 $\mid$ 130 & 17 $\mid$ 60$\phantom{0}$ &  & 270 & 200 &  & 120 & 80 \\
$(2^{18}, 2^{7})$ &  & 1.5e-10 & 2.6e-9 &  & 9.6e-1 & 1.1e0 &  & 2.5e0 & 1.9e0 &  & 5.5e0 & 2.1e0 &  & 39 & 14 &  & 5 & 7 &  & 14 & 21 &  & 18 $\mid$ 90$\phantom{0}$ & 16 $\mid$ 50$\phantom{0}$ &  & 210 & 150 &  & 110 & 70 \\
$(2^{17}, 2^{7})$ &  & 2.3e-9 & 9.7e-9 &  & 4.5e-1 & 4.7e-1 &  & 1.4e0 & 8.7e-1 &  & 1.6e0 & 9.5e-1 &  & 22 & 10 &  & 8 & 8 &  & 19 & 22 &  & 15 $\mid$ 70$\phantom{0}$ & 15 $\mid$ 50$\phantom{0}$ &  & 200 & 130 &  & 120 & 70 \\
$(2^{19}, 2^{8})$ &  & 7.0e-9 & 1.2e-9 &  & 5.4e0 & 4.7e0 &  & 1.2e1 & 7.6e0 &  & 1.7e1 & 8.9e0 &  & 18 & 10 &  & 8 & 9 &  & 24 & 27 &  & 18 $\mid$ 120 & 15 $\mid$ 70$\phantom{0}$ &  & 350 & 230 &  & 220 & 130 \\
$(2^{18}, 2^{8})$ &  & 6.3e-9 & 9.9e-9 &  & 2.1e0 & 2.6e0 &  & 6.8e0 & 4.1e0 &  & 9.8e0 & 4.4e0 &  & 21 & 11 &  & 9 & 10 &  & 23 & 26 &  & 19 $\mid$ 140 & 14 $\mid$ 70$\phantom{0}$ &  & 300 & 190 &  & 230 & 140 \\
$(2^{17}, 2^{8})$ &  & 4.6e-10 & 4.5e-9 &  & 1.4e0 & 1.1e0 &  & 2.8e0 & 1.7e0 &  & 4.0e0 & 1.9e0 &  & 18 & 7 &  & 11 & 11 &  & 23 & 27 &  & 18 $\mid$ 120 & 13 $\mid$ 60$\phantom{0}$ &  & 270 & 170 &  & 230 & 130 \\
$(2^{18}, 2^{9})$ &  & 4.0e-10 & 9.0e-10 &  & 6.3e0 & 6.0e0 &  & 1.5e1 & 8.4e0 &  & 1.8e1 & 9.1e0 &  & 9 & 5 &  & 16 & 13 &  & 29 & 31 &  & 20 $\mid$ 160 & 13 $\mid$ 80$\phantom{0}$ &  & 500 & 290 &  & 460 & 250 \\
$(2^{17}, 2^{9})$ &  & 9.3e-9 & 7.1e-9 &  & 4.0e0 & 3.3e0 &  & 7.4e0 & 4.7e0 &  & 8.4e0 & 5.0e0 &  & 8 & 4 &  & 21 & 17 &  & 27 & 30 &  & 17 $\mid$ 140 & 13 $\mid$ 90$\phantom{0}$ &  & 480 & 290 &  & 440 & 270 \\
$(2^{17}, 2^{10})$ &  & 9.4e-9 & 1.1e-9 &  & 1.1e1 & 9.1e0 &  & 2.1e1 & 1.1e1 &  & 2.3e1 & 1.2e1 &  & 3 & 2 &  & 38 & 24 &  & 24 & 30 &  & 16 $\mid$ 170 & 13 $\mid$ 100 &  & 900 & 500 &  & 860 & 470 \\
\midrule 
\multicolumn{31}{c}{$L={1},\;k={10}\%m$} \\
\midrule 
$(2^{20}, 2^{7})$ &  & 7.5e-9 & 7.7e-9 &  & 7.3e0 & 8.2e0 &  & 2.7e1 & 1.5e1 &  & 2.7e1 & 1.8e1 &  & 19 & 17 &  & 30 & 31 &  & 15 & 14 &  & 18 $\mid$ 110 & 14 $\mid$ 80$\phantom{0}$ &  & 450 & 400 &  & 120 & 80 \\
$(2^{19}, 2^{7})$ &  & 3.5e-9 & 3.7e-9 &  & 4.7e0 & 3.9e0 &  & 1.1e1 & 6.9e0 &  & 1.3e1 & 8.4e0 &  & 20 & 17 &  & 29 & 31 &  & 14 & 15 &  & 20 $\mid$ 100 & 15 $\mid$ 80$\phantom{0}$ &  & 310 & 270 &  & 120 & 70 \\
$(2^{18}, 2^{7})$ &  & 7.2e-9 & 2.4e-9 &  & 2.2e0 & 2.2e0 &  & 7.8e0 & 4.0e0 &  & 8.1e0 & 4.5e0 &  & 16 & 17 &  & 32 & 31 &  & 17 & 15 &  & 19 $\mid$ 160 & 15 $\mid$ 80$\phantom{0}$ &  & 240 & 190 &  & 120 & 70 \\
$(2^{17}, 2^{7})$ &  & 7.5e-11 & 2.1e-9 &  & 1.0e0 & 9.8e-1 &  & 2.9e0 & 1.7e0 &  & 2.9e0 & 2.0e0 &  & 22 & 16 &  & 27 & 30 &  & 15 & 15 &  & 18 $\mid$ 100 & 15 $\mid$ 70$\phantom{0}$ &  & 200 & 160 &  & 120 & 80 \\
$(2^{19}, 2^{8})$ &  & 1.7e-9 & 4.9e-9 &  & 1.0e1 & 1.1e1 &  & 3.0e1 & 1.7e1 &  & 3.1e1 & 1.8e1 &  & 13 & 11 &  & 37 & 37 &  & 18 & 18 &  & 19 $\mid$ 160 & 14 $\mid$ 100 &  & 370 & 280 &  & 220 & 140 \\
$(2^{18}, 2^{8})$ &  & 2.7e-9 & 8.6e-9 &  & 6.3e0 & 5.5e0 &  & 1.3e1 & 8.3e0 &  & 1.4e1 & 8.9e0 &  & 14 & 13 &  & 36 & 36 &  & 18 & 18 &  & 18 $\mid$ 140 & 14 $\mid$ 100 &  & 320 & 220 &  & 230 & 140 \\
$(2^{17}, 2^{8})$ &  & 9.9e-10 & 1.7e-9 &  & 3.0e0 & 2.5e0 &  & 7.7e0 & 3.7e0 &  & 8.0e0 & 4.1e0 &  & 16 & 7 &  & 35 & 40 &  & 17 & 19 &  & 17 $\mid$ 160 & 14 $\mid$ 100 &  & 290 & 190 &  & 230 & 140 \\
$(2^{18}, 2^{9})$ &  & 5.5e-9 & 1.9e-9 &  & 2.2e1 & 1.5e1 &  & 3.4e1 & 2.1e1 &  & 3.8e1 & 2.2e1 &  & 10 & 6 &  & 42 & 43 &  & 19 & 20 &  & 22 $\mid$ 220 & 16 $\mid$ 130 &  & 520 & 330 &  & 460 & 260 \\
$(2^{17}, 2^{9})$ &  & 6.3e-9 & 1.8e-9 &  & 1.1e1 & 8.3e0 &  & 1.7e1 & 1.1e1 &  & 1.8e1 & 1.1e1 &  & 6 & 6 &  & 46 & 44 &  & 19 & 20 &  & 21 $\mid$ 200 & 15 $\mid$ 140 &  & 500 & 280 &  & 450 & 240 \\
$(2^{17}, 2^{10})$ &  & 9.6e-9 & 6.8e-9 &  & 3.4e1 & 2.0e1 &  & 5.3e1 & 2.6e1 &  & 5.6e1 & 2.7e1 &  & 4 & 3 &  & 56 & 50 &  & 17 & 20 &  & 22 $\mid$ 280 & 15 $\mid$ 160 &  & 940 & 570 &  & 900 & 560 \\
\midrule 
\multicolumn{31}{c}{$L={10},\;k={1}\%m$} \\
\midrule 
$(2^{17}, 2^{7})$ &  & 2.0e-9 & 3.3e-9 &  & 7.1e0 & 1.1e1 &  & 1.9e1 & 1.4e1 &  & 2.5e1 & 1.5e1 &  & 27 & 14 &  & 3 & 3 &  & 18 & 24 &  & 18 $\mid$ 110 & 16 $\mid$ 90$\phantom{0}$ &  & 160 & 110 &  & 120 & 80 \\
$(2^{16}, 2^{7})$ &  & 7.3e-10 & 2.9e-9 &  & 2.9e0 & 4.2e0 &  & 1.0e1 & 5.7e0 &  & 1.2e1 & 6.3e0 &  & 27 & 15 &  & 3 & 3 &  & 17 & 22 &  & 16 $\mid$ 100 & 15 $\mid$ 70$\phantom{0}$ &  & 160 & 100 &  & 120 & 80 \\
$(2^{15}, 2^{7})$ &  & 4.0e-10 & 6.9e-10 &  & 3.3e0 & 2.7e0 &  & 6.8e0 & 3.4e0 &  & 7.8e0 & 3.6e0 &  & 20 & 15 &  & 5 & 4 &  & 22 & 23 &  & 18 $\mid$ 170 & 16 $\mid$ 80$\phantom{0}$ &  & 150 & 80 &  & 120 & 80 \\
$(2^{14}, 2^{7})$ &  & 4.3e-7 & 4.6e-9 &  & 1.4e0 & 1.0e0 &  & 3.0e0 & 1.5e0 &  & 4.6e0 & 1.6e0 &  & 11 & 13 &  & 5 & 6 &  & 17 & 23 &  & 24 $\mid$ 160 & 15 $\mid$ 70$\phantom{0}$ &  & 310 & 100 &  & 4250 & 80 \\
$(2^{16}, 2^{8})$ &  & 8.7e-11 & 3.0e-9 &  & 1.1e1 & 1.0e1 &  & 2.3e1 & 1.3e1 &  & 2.9e1 & 1.3e1 &  & 17 & 8 &  & 5 & 6 &  & 26 & 30 &  & 19 $\mid$ 190 & 14 $\mid$ 100 &  & 240 & 160 &  & 230 & 150 \\
$(2^{15}, 2^{8})$ &  & 6.8e-9 & 2.3e-9 &  & 5.9e0 & 5.0e0 &  & 1.0e1 & 6.3e0 &  & 1.2e1 & 6.9e0 &  & 12 & 8 &  & 9 & 7 &  & 27 & 30 &  & 21 $\mid$ 170 & 14 $\mid$ 100 &  & 250 & 140 &  & 220 & 140 \\
$(2^{14}, 2^{8})$ &  & 2.6e-9 & 1.2e-9 &  & 2.9e0 & 2.9e0 &  & 5.3e0 & 3.7e0 &  & 5.6e0 & 3.9e0 &  & 12 & 9 &  & 11 & 9 &  & 28 & 30 &  & 15 $\mid$ 160 & 13 $\mid$ 120 &  & 250 & 140 &  & 220 & 130 \\
$(2^{15}, 2^{9})$ &  & 1.0e-9 & 5.4e-9 &  & 1.4e1 & 1.3e1 &  & 3.1e1 & 1.7e1 &  & 3.1e1 & 1.8e1 &  & 8 & 5 &  & 14 & 11 &  & 31 & 34 &  & 20 $\mid$ 240 & 14 $\mid$ 150 &  & 470 & 260 &  & 450 & 250 \\
$(2^{14}, 2^{9})$ &  & 8.8e-9 & 2.4e-9 &  & 7.1e0 & 6.2e0 &  & 1.2e1 & 9.0e0 &  & 1.4e1 & 9.3e0 &  & 9 & 5 &  & 19 & 14 &  & 28 & 32 &  & 22 $\mid$ 210 & 15 $\mid$ 140 &  & 480 & 270 &  & 460 & 280 \\
$(2^{14}, 2^{10})$ &  & 4.5e-9 & 6.1e-9 &  & 2.7e1 & 1.8e1 &  & 3.9e1 & 2.2e1 &  & 4.2e1 & 2.3e1 &  & 3 & 3 &  & 35 & 21 &  & 26 & 33 &  & 21 $\mid$ 270 & 14 $\mid$ 180 &  & 890 & 500 &  & 860 & 490 \\
\midrule 
\multicolumn{31}{c}{$L={10},\;k={10}\%m$} \\
\midrule 
$(2^{17}, 2^{7})$ &  & 6.1e-9 & 6.3e-9 &  & 1.0e1 & 1.0e1 &  & 2.7e1 & 1.7e1 &  & 6.0e1 & 2.0e1 &  & 38 & 18 &  & 11 & 17 &  & 13 & 19 &  & 17 $\mid$ 190 & 15 $\mid$ 90$\phantom{0}$ &  & 220 & 170 &  & 120 & 80 \\
$(2^{16}, 2^{7})$ &  & 1.0e-11 & 7.9e-9 &  & 5.1e0 & 4.9e0 &  & 1.2e1 & 8.3e0 &  & 1.8e1 & 9.1e0 &  & 31 & 21 &  & 10 & 11 &  & 15 & 19 &  & 17 $\mid$ 130 & 14 $\mid$ 90$\phantom{0}$ &  & 170 & 110 &  & 120 & 80 \\
$(2^{15}, 2^{7})$ &  & 3.2e-9 & 4.8e-9 &  & 3.5e0 & 3.3e0 &  & 7.3e0 & 5.0e0 &  & 9.6e0 & 9.3e0 &  & 30 & 21 &  & 10 & 9 &  & 17 & 23 &  & 18 $\mid$ 170 & 15 $\mid$ 220 &  & 150 & 90 &  & 120 & 70 \\
$(2^{14}, 2^{7})$ &  & 4.2e-9 & 1.2e-9 &  & 1.4e0 & 1.5e0 &  & 3.9e0 & 2.1e0 &  & 4.8e0 & 2.2e0 &  & 33 & 19 &  & 6 & 8 &  & 17 & 22 &  & 16 $\mid$ 170 & 14 $\mid$ 100 &  & 130 & 80 &  & 120 & 80 \\
$(2^{16}, 2^{8})$ &  & 2.7e-9 & 6.3e-9 &  & 1.5e1 & 1.4e1 &  & 3.2e1 & 1.8e1 &  & 4.1e1 & 2.0e1 &  & 16 & 14 &  & 19 & 13 &  & 22 & 25 &  & 20 $\mid$ 230 & 15 $\mid$ 120 &  & 280 & 170 &  & 230 & 140 \\
$(2^{15}, 2^{8})$ &  & 4.8e-10 & 1.4e-9 &  & 8.7e0 & 7.5e0 &  & 1.7e1 & 1.1e1 &  & 1.9e1 & 1.2e1 &  & 19 & 15 &  & 16 & 16 &  & 23 & 25 &  & 17 $\mid$ 220 & 15 $\mid$ 150 &  & 270 & 160 &  & 230 & 140 \\
$(2^{14}, 2^{8})$ &  & 3.4e-9 & 1.5e-9 &  & 5.1e0 & 4.3e0 &  & 1.3e1 & 5.4e0 &  & 2.2e1 & 5.7e0 &  & 25 & 15 &  & 16 & 15 &  & 21 & 25 &  & 21 $\mid$ 460 & 15 $\mid$ 140 &  & 240 & 160 &  & 230 & 150 \\
$(2^{15}, 2^{9})$ &  & 9.7e-9 & 6.8e-10 &  & 3.1e1 & 2.0e1 &  & 5.0e1 & 3.4e1 &  & 5.2e1 & 3.5e1 &  & 12 & 8 &  & 27 & 27 &  & 25 & 27 &  & 17 $\mid$ 320 & 16 $\mid$ 240 &  & 480 & 260 &  & 450 & 260 \\
$(2^{14}, 2^{9})$ &  & 8.1e-10 & 4.6e-9 &  & 1.5e1 & 1.1e1 &  & 2.5e1 & 1.6e1 &  & 2.9e1 & 1.6e1 &  & 11 & 8 &  & 29 & 27 &  & 24 & 27 &  & 22 $\mid$ 350 & 15 $\mid$ 220 &  & 470 & 280 &  & 460 & 280 \\
$(2^{14}, 2^{10})$ &  & 8.7e-9 & 7.9e-10 &  & 4.7e1 & 3.6e1 &  & 7.5e1 & 4.8e1 &  & 8.2e1 & 4.9e1 &  & 5 & 4 &  & 48 & 39 &  & 21 & 25 &  & 19 $\mid$ 420 & 15 $\mid$ 310 &  & 900 & 500 &  & 900 & 530 \\
\bottomrule 
\end{tabular}
        \end{adjustbox}
\end{subtable}
\caption{Detailed ALM performance on synthetic instances with linear (lin.) and convex quadratic (quad.) objective. The column
$\% \text{T}_c$ indicates the proportion of total solve time; columns $\nabla\varphi_\sigma$ and $V^{-1} v$ report extra time, assuming that $\varphi_\sigma$ and $\nabla\varphi_\sigma$ have already been evaluated, respectively.}
\label{tab:detail_synthetic}
\end{sidewaystable}

\subsection{A quantile regression experiment on real data}\label{sec:experiments:qr}

In this section, we consider an application of the superquantile minimization problem to estimating the $\tau$-conditional quantile under a linear model \citep[Chapter 8.C]{royset2021primer}. Given input-output data $\{(a^i, b_i)\}_{i=1}^m\subseteq \mathbb{R}^n \times \mathbb{R}$ and a scalar $\tau\in (0,1)$ such that $k(\tau)={(1-\tau)\cdot m}$ is an integer, the problem takes the following form:
\begin{align}\label{eq:quantile_regression:superquantile}
\begin{array}{cl}
    \displaystyle\operatornamewithlimits{minimize}_{x,t} & \displaystyle t + x^\top \left(\frac{1}{m}\sum_{i=1}^m a^i\right)\\[0.2in]
     \text{subject to} &\maxsum_{k(\tau)}[\big]{\{b_i - x^\top {a^i} - t\}_{i=1}^m} \leq 0.
\end{array}
\end{align}
See \cref{apx:qr} for the detailed derivation of the above formulation.

Using ride-share data from the \cite{nyc2023data} and weather data from \cite{weather2023data} for October 2022 through September 2023, we construct a dataset containing $m=32,\!534,\!601$ trips with $n=19$ features to predict the tip received by the driver for rides paying with credit card.
The full dataset $\{(a^i, b_i)\}_{i=1}^{m}$ is available at \url{https://doi.org/10.5281/zenodo.11183368}.

Based on the numerical results in \cref{sec:experiments:synthetic}, we have narrowed the focus of our study to the most efficient methods: ALM, GRB, and G-OA for the superquantile formulation \eqref{eq:quantile_regression:superquantile}. Additionally, we use  Gurobi's barrier method to solve the following linear reformulation of quantile regression \eqref{eq:quantile_regression:superquantile}:
\begin{align}\label{eq:quantile_regression:lp}
\begin{array}{cl}
     \displaystyle\operatornamewithlimits{minimize}_{x,(z^+,z^-)\in\R^{2m}} & \displaystyle \frac1m \sum_{i=1}^m (1-\tau) \cdot z_i^- + \tau\cdot z_i^+\\[0.1in]
        \text{subject to}  & z_i^+ - z_i^- = b_i - x^\top a^i,\; z_i^+\geq 0, \; z_i^- \geq 0, \;  i\in\{1,\ldots,m\},
\end{array}
\end{align}
which we refer to as  G-QR.
For all of the above mentioned methods, we set a time limit of $3,\!600$s.

We  conduct  two experiments.
For the first experiment, we predict the $\tau$-th top-quantile for $\tau\in\{1\%,10\%,50\%\}$ using all methods at a precision of $\epsilon=10^{-4}$.
Notably,  $\tau=50\%$ presents the most computational challenge for ALM since for $\tau>0.5$, the superquantile formulation \eqref{eq:quantile_regression:superquantile} can be inverted to $\bar\tau = 1-\tau<0.5$ by reversing the sign of $\{b_i\}_{i=1}^m$ and $\{a^i\}_{i=1}^m$.
Secondly, we employ ALM to approximate the full distribution over a region of interest by solving prediction problems for the following 22 quantiles: 
$\{10^{-3},10^{-2},0.025,0.05,\ldots,0.475,0.5\}$,
each with a tolerance of $\epsilon=10^{-4}$. With the exception of the first quantile, each subsequent instance is warm-started using the solution from the preceding smaller  quantile.
Experiments are conducted in a Linux environment on a Dell workstation with eight physical Intel(R) Xeon(R) W-2145 CPUs @ 3.70GHz with two threads per core and 125GB of RAM.
Gurobi-based methods required about 85GB of RAM, and limited the number of features used for prediction, while ALM used  about 40GB of RAM.

Numerical experiments are summarized in \cref{tab:cold,tab:warm}.
\cref{tab:cold} shows that ALM is up to $\approx 20$x faster than all three of the Gurobi-based methods for small (or large, by using $\bar\tau$) quantile problems.
For the most difficult case when $\tau=0.5$, ALM required 1.5x the time needed by the most efficient Gurobi method (G-QR).
When solving instances from \cref{tab:cold}, the Gurobi methods were able to find feasible primal and dual solutions to numerical precision, but  the error $\eta$ is determined by the duality gap.
In large-scale problems, the linear systems solved by the barrier method can become prohibitively ill-conditioned. Even the G-OA method encountered this problem, as the set of active scenarios could change significantly from one iteration to the next.
In contrast, ALM has handled these challenging instances more reliably.

\cref{tab:warm} summarizes the results from the second experiment, demonstrating that warm-starting can substantially decrease the computational time for challenging larger-$\tau$ cases.
For instance, when warm-started, ALM can solve the $\tau = 0.5$
case approximately $1.4$ times faster than the best-performing Gurobi method.
Overall, ALM can obtain the solution path for $\approx40$ quantiles ranging from $\tau = 10^{-3}$ to $\tau = 1-10^{-3}$ with a precision of $\epsilon=10^{-4}$ in roughly the same duration that the most efficient Gurobi method takes to solve two problems in the $\tau\approx 1\%$ range.

\begin{table}[ht]
\centering
{\footnotesize
\begin{tabular}{cccccccccccc}
\toprule 
$k = \tau\%m$ &  & \multicolumn{4}{c}{Time [s]} & \multicolumn{4}{c}{$\eta$} \\
\cmidrule{3-6} \cmidrule{8-11}  &  & ALM & GRB & G-OA & G-QR &  & ALM & GRB & G-OA & G-QR \\
\midrule 
1 &  & 1.3e2 & 3.6e3 & 6.2e3 & 3.6e3 &  & 4.4e-9 & 2.5e-1 & 2.9e-1 & 2.3e1 \\
10 &  & 2.1e2 & 3.6e3 & 4.2e3 & 3.6e3 &  & 2.5e-9 & 2.7e-5 & 7.3e-3 & 1.1e-2 \\
25 &  & 6.4e2 & 7.6e2 & 4.2e3 & 2.5e3 &  & 2.5e-9 & 4.3e-9 & 1.6e-5 & 8.4e-9 \\
50 &  & 3.6e3 & 8.2e2 & 3.6e3 & 1.7e3 &  & 3.6e-4 & 1.4e-9 & 6.9e-4 & 9.6e-8 \\
\midrule
\end{tabular}
}%
\caption{
Results for the computation of the solution path for the cold-start quantile regression with an accuracy of  $\epsilon=10^{-8}$.  
The time limit is one hour.}
\label{tab:cold}
\end{table}

\begin{table}[ht]
\captionsetup{width=0.99\textwidth}
\centering
\footnotesize{
\begin{tabular}{ccccccccc}
\toprule
$k = \tau\%m$ & T $\coloneqq$ Time [s] & $\eta$ & ALM iter & \% T: $V^{-1}v$ & \% T: $\nabla\varphi$ & \% T: sort & \% T: $\mathopfont{proj}_{\overline{{\cal B}}}$ \\
\midrule
1.00e-3 & 1.3e2 & 6.5e-5 & 10 & 5 & 10 & 14 & 11 \\
1.00e-2 & 1.3e2 & 3.2e-5 & 4 & 1 & 5 & 53 & 12 \\
2.50e-2 & 1.4e2 & 3.7e-5 & 4 & 2 & 4 & 55 & 12 \\
5.00e-2 & 1.7e2 & 4.2e-5 & 5 & 4 & 4 & 54 & 11 \\
7.50e-2 & 1.5e2 & 8.5e-5 & 3 & 4 & 3 & 61 & 11 \\
1.00e-1 & 1.5e2 & 4.6e-5 & 3 & 5 & 3 & 60 & 10 \\
1.25e-1 & 1.5e2 & 6.7e-5 & 2 & 4 & 2 & 64 & 11 \\
1.50e-1 & 1.9e2 & 6.0e-5 & 4 & 7 & 3 & 56 & 10 \\
1.75e-1 & 2.0e2 & 5.4e-5 & 4 & 7 & 3 & 57 & 10 \\
2.00e-1 & 1.8e2 & 8.1e-5 & 3 & 7 & 2 & 59 & 9 \\
2.25e-1 & 1.8e2 & 8.7e-5 & 3 & 7 & 2 & 60 & 9 \\
2.50e-1 & 1.9e2 & 8.2e-5 & 3 & 8 & 2 & 59 & 9 \\
2.75e-1 & 2.2e2 & 5.7e-5 & 3 & 8 & 2 & 60 & 9 \\
3.00e-1 & 2.2e2 & 3.4e-5 & 3 & 8 & 2 & 60 & 9 \\
3.25e-1 & 2.0e2 & 5.5e-5 & 2 & 7 & 1 & 64 & 9 \\
3.50e-1 & 1.9e2 & 5.1e-5 & 2 & 8 & 2 & 62 & 9 \\
3.75e-1 & 2.3e2 & 3.9e-5 & 3 & 9 & 2 & 59 & 8 \\
4.00e-1 & 2.3e2 & 4.7e-5 & 3 & 10 & 2 & 58 & 8 \\
4.25e-1 & 2.2e2 & 8.7e-5 & 2 & 7 & 1 & 64 & 9 \\
4.50e-1 & 2.1e2 & 6.2e-5 & 2 & 8 & 1 & 62 & 8 \\
4.75e-1 & 2.7e2 & 3.3e-5 & 3 & 10 & 1 & 60 & 8 \\
5.00e-1 & 3.3e2 & 7.3e-5 & 4 & 12 & 1 & 57 & 8 \\
\midrule
\footnotesize{Total time:} & 4.3e3 & $-$ & $-$ &  \footnotesize{(avg)} 7 & \footnotesize{(avg)}  2 & \footnotesize{(avg)}  58 & \footnotesize{(avg)}  9 \\
\bottomrule
\end{tabular}
}%
\caption{Results for the computation of the solution path with an accuracy of  $\epsilon=10^{-4}$.
The last four columns  show the proportion of time spent  computing $V^{-1}v$, $\nabla\varphi$, sorting, and $\mathopfont{proj}_{\overline{{\cal B}}}$
relative to the total runtime. The final row presents the average  proportion across all instances.}
\label{tab:warm}
\end{table}

\section{Conclusions}
\label{sec:conclusions}

In this paper, we introduced a semismooth-Newton-based augmented Lagrangian method for the  superquantile-constrained optimization problem \eqref{eq:source_problem} when the uncertain parameter $\omega$ is supported on a finite set of scenarios with equal probability.
The proposed approach ensures that the generalized Jacobian in the semismooth Newton method can be obtained cheaply via the structured sparsity inherent to the superquantile operator.
In numerical experiments conducted on synthetic and realistic datasets where the number of scenarios greatly exceeds  the number of decision variables ($m\gg n$), we find that the ALM method improves upon off-the-shelf solvers such as Gurobi and OSQP as well as manually-tuned active set outer-approximation methods,
and can be faster and more reliable than the specialized software suite PSG.

Finally, when the uncertain function $g^{\,\ell}(\,\cdot\,;\omega^j)$ is nonlinear, the ALM method is still expected to exhibit improved performance relative to the generic nonlinear solvers.
Underlying this anticipation is the form of the generalized Jacobian 
\[
\sum_{i=1}^{mL}\nabla^2 G_i(x) \cdot\bigl(G(x)-\proj_{\overline{\mathcal{B}}}[\big]{G(x)}\bigr)_i+ JG(x)^\top\cdot (I-J_{\overline{\mathcal{B}}})\cdot JG(x), 
\]
which only depends on the Hessian of the effective tail scenarios.
Second-order solvers that are not able to exploit this sparsity must evaluate the Hessian on the entire scenario set, which could become prohibitively expensive.
Further improvements are expected when the function $G$ is  nonlinear but with a simple structure,
since the corresponding Hessian $\nabla^2 G_i(x)$ in the first term may be computed cheaply and the second term is always cheap to compute due to the special sparsity structure of $(I-J_{\overline{\mathcal{B}}})$, whereas standard convex solvers may have to resort to extended reformulations.

\section*{Acknowledgements}
The authors are partially supported by NSF grants CCF-2416172 and DMS-2416250, and NIH grant R01CA287413.  The authors thank Eric Sager Luxenberg and Johannes Royset for helpful discussion of the algorithms presented in the manuscript.

\bibliography{refs_abbrev}

\appendix

\section{Proofs of \texorpdfstring{\cref{prop:explicit_gen_jac}}{Lem2} and \texorpdfstring{\cref{lem:ineqalm}}{Lem3}}
\label{apx:proof}
 
{\bf Proof of \cref{prop:explicit_gen_jac}}. 
    Case 1 is clear.
    Case 2
    relies on the characterization of the critical cone given in \eqref{eq:criticalcone_analytic}.
    The corresponding projection problem can be stated as follows.
    For any $y^0\in\R^m$ with $\maxsum_k{y^0}>r$, the directional derivative of $\projection_{\mathcal{B}_k}$ along direction $(u^0,v^0,w^0)$ is given by projecting $(u^0,v^0,w^0)$ onto the critical cone $\bar{\mathcal{D}}$ at $\bar{y} = \proj_{\mathcal{B}_k}[\big]{\sorth{y}^0}$.
    In Case 2,
    the projection is onto a linear subspace, and hence the linear operator can be recovered by inspecting the projection onto the critical cone.
    It can be shown \cite[Section 4.3.4]{wu2014moreau} that for $u\in\R^{\abs{\bar\alpha}}$, $v\in\R^{\abs{\bar\beta}}$, $w\in\R^{\abs{\bar\gamma}}$, this problem can be expressed as

    \vspace{-5mm}
    \begin{align}
        \label{eq:criticalcone_projection}
        \renewcommand{\arraystretch}{1.5}
        \begin{array}{cl}
         \displaystyle \min_{u,v} & \tfrac12 \norm{u-u^0}_2^2 + \tfrac12 \norm{v-v^0_\sigma}_2^2 + \tfrac12\norm{w-w^0}_2^2\\
         \text{s.t.} & v_i \geq v_{[k-\bar k_0]},\;\forall i \in \bigl\{1,\ldots,\abs{\bar\beta_1}\bigr\},\\
         & v_i = v_{[k-\bar k_0]},\;\forall i \in\bigl\{\abs{\bar\beta_1}+1,\ldots,\abs{\bar\beta_1}+\abs{\bar\beta_2}\bigr\},\\
         & v_i \leq v_{[k-\bar k_0]},\;\forall i \in\bigl\{\abs{\bar\beta_1}+\abs{\bar\beta_2},\ldots,\abs{\bar\beta_1}+\abs{\bar\beta_2}+\abs{\bar\beta_3}\bigr\},\\
         & \ind_{\bar\alpha}^\top u + \bar{\mu}_{\bar\beta}^\top v = 0
        \end{array}
    \end{align}
    \vspace{-5mm}
    
    where $\bar\beta_1 = \{i\in\bar\beta: (\bar\mu_{\bar\beta})_i=1\}$, $\bar\beta_2 = \{i\in\bar\beta: (\bar\mu_{\bar\beta})_i\in(0,1)\}$, and $\bar\beta_3 = \{i\in\bar\beta: (\bar\mu_{\bar\beta})_i=0\}$.
    Problem \eqref{eq:criticalcone_projection} is obtained by using {Lemma 2.5 of \citet{wu2014moreau}}, where $\sigma$ is a permutation of $\bar\beta$ that sorts the $\bar\beta_1$ components, preserves the $\bar\beta_2$ components, and sorts the $\bar\beta_3$ components so that $(v_\sigma)_{\bar\beta_1} = (\sorth v_{\bar\beta_1})$, $(v_\sigma)_{\bar\beta_1} = (\sorth v_{\bar\beta_1})$, and $(v_\sigma)_{\bar\beta_3} = (\sorth v_{\bar\beta_3})$.
    The condition $\sorth y^0_{\bar{k}_0+1} - \bar\lambda < \bar{y}_k < y^0_{\bar{k}_0+1}$ implies that $0 < \bar\mu_{\bar\beta} = (\sorth y^0_{\bar\beta} - \bar y_{\bar\beta})/\bar\lambda < \ind_{\bar\beta}$ since
        \begin{align*}
            \sorth y_{\bar k_1}^0 - \bar\theta > 0 \implies\cdots\implies
            \sorth y_{\bar k_0+1}^0 - \bar\theta>0,\quad
            \sorth y_{\bar k_0+1}^0 - \bar\theta < \bar\lambda \implies\cdots\implies
            \sorth y_{\bar k_1}^0 - \bar\theta < \bar\lambda.
        \end{align*}
        Thus $\bar\beta_1\cup\bar\beta_3=\emptyset$.
        Since $\bar\beta_1\cup\bar\beta_3=\emptyset$, projection onto the critical cone simplifies in \eqref{eq:criticalcone_projection}, and the permutation $\sigma$ may be disregarded.
        The linear subspace is then easily verified as that defined in \eqref{eq:criticalcone_projection_case2} by noting that the equality constraints can be arranged successively.

        Case 3 follows the same argument as Case 2, where the constraints on $v_i$ for $i\in\{1,\ldots,\abs{\bar\beta_1}$ in \eqref{eq:criticalcone_projection} are automatically satisfied by the sorting. \hfill $\square$

\vskip 0.2in

{\bf Proof of \cref{lem:ineqalm}}.
    Let $h(y)\coloneqq\tfrac12\|\max(\sorth g-y,0)\|_2^2$.
    Then $h$ is continuously differentiable with gradient $\nabla h(y) = -\max(\sorth g-y,0)$.
    Since both \eqref{eq:ineqalm_proj} and \eqref{eq:eqalm_proj} are convex, in order to prove this lemma, it suffices to show that
    $\bar y - \sorth g\in \bigl(\mathcal{T}_{\mathcal{B}_k}(\bar y)\bigr)^* \implies -\max(\sorth g-\bar y,0) \in \bigl(\mathcal{T}_{\mathcal{B}_k}(\bar y)\bigr)^*$ is an element in the dual of the tangent cone of $\mathcal{B}_k$ at $\bar{y}$.
    Using the form of the KKT conditions which ensures that $\bar{y}\leq\sorth g$, the desired conclusion holds trivially since
    \begin{flalign*}
    && \bar{y} - \sorth{g} &= \max(\bar{y} - \sorth{g}, 0) + \min(\bar{y} - \sorth{g}, 0) = \max(\bar{y} - \sorth{g}, 0) - \max(\sorth{g} - \bar{y}, 0) = -\max(\sorth{g} - \bar{y}, 0). && \hfill \square
    \end{flalign*}

\section{Constraint non-degeneracy}
\label{sec:cq_non-degeneracy}
Here we describe the constraint non-degeneracy conditions of \eqref{eq:source_problem} with a single superquantile term in the objective function.
Such conditions are used in \cref{sec:p_alm_framework} when discussing convergence guarantees for the ALM method.
The constraint non-degeneracy follows the approach outlined in \cite{ding2016variational}.

Since the top-$k$-sum operator $\mathopfont{T}_k$  is positive homogeneous, we can consider the conic problem
\begin{equation}\label{eq:source_problem_conic}
    \begin{array}{cl}
         \displaystyle \min_{x\in X} & f(x) + t\\
         \text{s.t.} & \bigl(G(x), t\bigr) \in \mathcal{K}_{(k)}\coloneqq\mathsf{epi}\bigl(\maxsum_k{\cdot}\bigr),
    \end{array}
\end{equation}
where $\mathsf{epi}$ denotes the epigraph.
Let $\mathcal{G}(x,t)\coloneqq\bigl(G(x),t\bigr)$.
Then at any $(\bar x,\bar t)\in X\times\R$, Robinson's constraint qualification
\begin{align*}
    J\mathcal{G}(\bar x,\bar t) (X\times\R) + \tangentcone[\big]{\mathcal{K}_{(k)}}{\mathcal{G}(\bar x, \bar t)} = \R^m\times\R
\end{align*}
holds since there exists an interior to the feasible region by taking $t$ large enough.
In particular, this means that at an optimal solution $x^*\in X$, the corresponding set of multipliers are nonempty, convex, and compact \citep{zowe1979regularity}.
By \cite{robinson1984local,robinson1987local}, the constraint non-degeneracy holds at $\bigl(\bar x, \maxsum_k{G(\bar x)}\bigr)$ if and only if
\begin{align}
\label{eq:constraint_non-degeneracy}
    J G(\bar x) X + \mathcal{T}^{\mathsf{lin}}\bigl(G(\bar x)\bigr) = \R^m,
\end{align}
where $\mathcal{T}^{\mathsf{lin}}(y)$ is the linear subspace of $\R^m$ defined by
\begin{align*}
    \mathcal{T}^{\mathsf{lin}}(y) &\coloneqq \{d\in\R^m:\maxsum'_k{y;d}=-\maxsum'_k{y;-d}\}.
\end{align*}

Using \citep[Theorem 23.4]{rockafellar1970convex} and the form of $\partial\maxsum_k{y}$ in \citet[Lemma 2.2]{wu2014moreau}, the directional derivative of the top-$k$-sum operator at a point $y\in\R^m$ along a direction $d\in\R^m$ is given by
\begin{align}
    \maxsum'_k{y;d} &= \sup_{\mu}\{\inner{\mu}{d} : \mu\in\partial\maxsum_k{y}\}
    = \ind_{\alpha}^\top (d_\pi)_{\alpha} + \maxsum_{k-k_0}[\big]{(d_\pi)_\beta},
\end{align}
where $\alpha,\beta,\gamma$ are the index-sets associated with $y_\pi = \sorth{y}$, and $k_0,k_1$ are the corresponding sorting indices.
That is, the directional derivative is linear in the $\alpha$ components.
Letting $\bar{y}\coloneqq G(\bar{x})$, we obtain
\begin{align}
    \mathcal{T}^{\mathsf{lin}}(\bar y) = \bigl\{d\in\R^m : \maxsum_k[\big]{(d_\pi)_\beta} = -\maxsum_k[\big]{-(d_\pi)_\beta} \bigr\} = \bigl\{d \in\R^m : (d_\pi)_\beta = s\ind_\beta,\;s\in\R\bigr\}.
\end{align}

In the case when $X=\R^n$, and we assume that $\bar y = G(\bar{x})$ is sorted, then system \eqref{eq:constraint_non-degeneracy} is equivalent to the existence of $v\in\R^n$ for any $w\in\R^m$ such that
\begin{align*}
    JG(\bar{x})v + \mathcal{T}^{\mathsf{lin}}(\bar y) = w &\iff
    \begin{bmatrix}
        JG(\bar{x})_{\alpha,:}\\JG(\bar{x})_{\beta,:}\\JG(\bar{x})_{\gamma,:}
    \end{bmatrix}
    v
    +
    \begin{bmatrix}
        u_{\alpha}\\u_{\beta}\\u_{\gamma}
    \end{bmatrix}
    =
    \begin{bmatrix}
        w_{\alpha}\\w_{\beta}\\w_{\gamma}
    \end{bmatrix},\;\;
    u\in\mathcal{T}^{\mathsf{lin}}(\bar y)\\*
    &\iff
    JG(\bar{x})_{\beta,:}v + s\ind_\beta = w_\beta,\;s\in\R,
\end{align*}
where the second equivalence holds due to the form of $\mathcal{T}^{\text{lin}}(\bar y)$ because $u_\alpha$ and $u_\gamma$ are free.
Therefore, the constraint non-degeneracy condition \eqref{eq:constraint_non-degeneracy} holds if and only if $[G'(\bar{x})_{\beta,:} \;\; \ind_\beta]$ has full row rank, which requires $\abs{\beta}\leq n+1$.

For comparison, consider the linear independence constraint qualification of the nonlinear reformulated problem  in  \eqref{eq:nlp}.
The following conditions can be verified to hold at a solution $(x,t,z)\in\R^{n+1+m}$ with an active top-$k$-sum constraint (left) leading to the following active constraint matrix (right):

\vspace*{-15pt}\noindent
\begin{minipage}[t]{0.5\textwidth}
\vspace*{-10pt} %
\begin{subequations}
    \begin{alignat}{5}
    0 &= t + \tfrac1k \ind^\top z\\
    z_\alpha &= G(x)_\alpha - t\ind_\alpha &\qquad& z_\alpha &&> 0\\
    z_\beta  &= G(x)_\beta - t\ind_\beta  &\qquad& z_\beta  &&= 0\\
    z_\gamma &> G(x)_\gamma - t\ind_\gamma &\qquad& z_\gamma &&= 0
    \end{alignat}
\end{subequations}
\end{minipage}%
\begin{minipage}[t]{0.5\textwidth}
\vspace*{0pt} %
\begin{equation}
    \begin{aligned}
    C\coloneqq
    \begin{bmatrix}
        0 & 1 & \tfrac1k\ind^\top\\
        -G'(x)_{\alpha,:} & \ind_\alpha & I_{\alpha,:}\\
        -G'(x)_{\beta,:} & \ind_\beta & I_{\beta,:}\\
        0 & 0 & I_{\beta,:}\\
        0 & 0 & I_{\gamma,:}
    \end{bmatrix}.
\end{aligned}
\end{equation}
\end{minipage}

It is clear that $C$ has full row rank if and only if two conditions hold: (i) the submatrix $[JG(x)_{\beta,:}\;\;\ind_\beta]$ has full row rank, \ie $\rank{JG(x)_{\beta,:}\;\;\ind_\beta}=\abs{\beta}$ and $\abs{\beta}\leq n+1$; and (ii) $\abs{\beta}\leq n$.
Condition (i) arises from observing that the second and fifth block rows of $C$ are always linearly independent from each other and from the remaining block rows because the index sets $\alpha$, $\beta$, and $\gamma$ are distinct.
Therefore, the only possibility of linear dependence is due to blocks three and four.
Finally, given a full row rank matrix $E$, any block matrix $\begin{bmatrix} D & E \\ & E\end{bmatrix}$ has full row rank if and only if $D$ has full row rank \citep{horn2013matrix}.
Thus we recover condition (i).
Condition (ii) is obtained from the requirement that $1+m+\abs{\beta}\leq n+m+1$, \ie the number of rows of $C$ must not exceed the number of columns.
Therefore, necessary and sufficient conditions for satisfying the constraint non-degeneracy condition under the standard reformulation are slightly stronger than the top-$k$-sum operator formulation.

\vspace{-2mm}
\section{Reformulation of the quantile regression}
\label{apx:qr}
\vspace{-3mm}
Given input-output data $\{(a^i,b_i)\}_{i=1}^m \subseteq \R^n\times\R$ and some $\tau \in (0,1)$ with integer $(1-\tau)m$, consider the empirical quantile regression problem given by
\begin{align}\label{eq:quantile}
\begin{array}{rl}
     \displaystyle  \displaystyle\operatornamewithlimits{minimize}_{x\in \mathbb{R}^n, x_0\in \mathbb{R}} &  \displaystyle\frac1m \sum_{i=1}^m \rho_\tau \bigl(b_i - x^\top {a^i} -x_0 \bigr), 
\end{array}
\end{align}
where $z=z^+ - z^-$ with $z^+ = \max(z,0)$ and $z^- = \max(-z,0)$, and $\rho_\tau(z) = (1-\tau)\cdot z^- + \tau \cdot z^+$ is the so called ``check loss'' of \cite{koenker1978regression}.
By simple manipulation, we see that
\begin{align*}
    \rho_\tau(z) = (1-\tau)\cdot z^- + \tau\cdot z^+= (1-\tau)z^- + \tau z^+ + (z^+-z^+) = (1-\tau)\cdot(-z) + z^+.
\end{align*}
It is easy to see that one can scale the objective in \eqref{eq:quantile} by $(1-\tau)^{-1}$ and rewrite the problem as
\begin{align*}
\begin{array}{ll}
  \quad  \displaystyle\operatornamewithlimits{minimize}_{x\in \mathbb{R}^n, x_0\in \mathbb{R}} \quad  \displaystyle \frac1{(1-\tau)\cdot m} \sum_{i=1}^m \max\bigl(0,b_i - x^\top a^i-x_0\bigr) - \frac1{m} \sum_{i=1}^m \bigl(b_i - x^\top a^i - x_0\bigr)\\
=
\displaystyle\operatornamewithlimits{minimize}_{x\in \mathbb{R}^n} \; \displaystyle\operatornamewithlimits{minimize}_{x_0\in \mathbb{R}} \quad  \displaystyle x_0 + \frac1{(1-\tau)\cdot m} \sum_{i=1}^m \max\bigl(0,b_i - x^\top a^i-x_0\bigr) - \frac1{m} \sum_{i=1}^m \bigl(b_i - x^\top a^i\bigr)\\[0.2in]
=
\begin{array}{rl}
     \displaystyle \displaystyle\operatornamewithlimits{minimize}_{x\in \mathbb{R}^n} \quad \displaystyle \superquantile_{\tau}[\big]{b_i - x^\top a^i} - \frac1{m} \sum_{i=1}^m \bigl(b_i - x^\top a^i\bigr)
\end{array} \\[0.2in]
=
\begin{array}{rl}
     \displaystyle \displaystyle\operatornamewithlimits{minimize}_{x,t} & \displaystyle t - \frac1{m} \sum_{i=1}^m \bigl(b_i - x^\top a^i\bigr)\;\;
     \text{subject to}\;\; \superquantile_{\tau}[\big]{\{b_i - x^\top a^i\}_{i=1}^m} \leq t
\end{array}.
\end{array}
\end{align*}
By the translation equivariance property of the superquantile operator, we obtain formulation \eqref{eq:quantile_regression:superquantile}.

\end{document}